%% file: V3--An-overview-of-horospherical-varieties-and-coloured-fans.tex
\pdfoutput=1 

\documentclass[12pt, letterpaper, leqno]{amsart}
\usepackage[letterpaper, margin=1in]{geometry}
\usepackage[utf8]{inputenc}
\usepackage[style=alphabetic, maxnames=99]{biblatex}
\AtBeginBibliography{\small}
\usepackage{microtype}
\microtypesetup{protrusion=false}
\usepackage{xpatch}
\usepackage[english]{babel}
\usepackage{amsfonts}
\usepackage{amsmath}
\usepackage{amssymb}
\usepackage{mathrsfs}
\usepackage{stmaryrd}
\usepackage{nicematrix}

\usepackage{amsthm}
\xpatchcmd{\proof}{\itshape}{\bfseries}{}{}
\usepackage{thmtools}
\usepackage{thm-restate}
\usepackage{tikz}
\usepackage{tikz-cd}
\usetikzlibrary{arrows.meta,automata,positioning}
\usetikzlibrary{svg.path}
\usetikzlibrary{decorations.pathreplacing}
\usepackage{amscd}
\usepackage{quiver}
\usepackage{tikzit}
\usepackage{pdfpages}
\usepackage{float}
\usepackage[shortlabels]{enumitem}
\setlist[enumerate,1]{label={(\arabic*)}}
\usepackage[labelfont=bf]{caption}
\usepackage[linewidth=0pt]{mdframed}
\usepackage{subfig}
\usepackage{makecell}
\usepackage{multicol}
\usepackage{listings}
\usepackage{fancyhdr}
\usepackage[bottom, perpage]{footmisc}
\usepackage{etoolbox}

\patchcmd{\section}{\scshape}{\bfseries}{}{}
\makeatletter
\renewcommand{\@secnumfont}{\bfseries}
\makeatother

\allowdisplaybreaks

\usepackage{xr}
\usepackage{subfiles}

\usepackage{setspace}
\usepackage[shortcuts]{extdash}

\makeatletter   
\xpatchcmd{\@tocline}
{\hfil\hbox to\@pnumwidth{\@tocpagenum{#7}}\par}
{\ifnum#1<0\hfill\else\dotfill\fi\hbox to\@pnumwidth{\@tocpagenum{#7}}\par}
{}{}
\makeatother    
\makeatletter
\def\l@subsection{\@tocline{2}{0pt}{4pc}{6pc}{}}
\def\l@subsubsection{\@tocline{3}{0pt}{8pc}{8pc}{}}
\makeatother
\usepackage{xstring}
\usepackage{hyperref}
\usepackage[nameinlink,capitalize]{cleveref}
\usepackage{url}
\hypersetup{
	colorlinks   = true, 
	urlcolor     = cyan, 
	linkcolor    = blue, 
	citecolor   = blue,
	breaklinks = true
}

\numberwithin{equation}{section}
\numberwithin{figure}{section}

\makeatletter
\newcommand{\labeltext}[3][]{%
	\@bsphack%
	\csname phantomsection\endcsname
	\def\tst{#1}%
	\def\labelmarkup{}
	\def\refmarkup{}%
	\ifx\tst\empty\def\@currentlabel{\refmarkup{#2}}{\label{#3}}%
	\else\def\@currentlabel{\refmarkup{#1}}{\label{#3}}\fi%
	\@esphack%
	\labelmarkup{#2}
}
\makeatother

\declaretheorem[numberwithin=section,numberlike=equation,style=plain]{theorem}

\declaretheorem[numbered=no,style=plain,name=Theorem]{theorem*}
\declaretheorem[numberwithin=section,numberlike=equation,style=plain]{proposition}

\declaretheorem[numbered=no,style=plain,name=Proposition]{proposition*}
\declaretheorem[numberwithin=section,numberlike=equation,style=plain]{lemma}

\declaretheorem[numbered=no,style=plain,name=Lemma]{lemma*}
\declaretheorem[numberwithin=section,numberlike=equation,style=plain]{corollary}

\declaretheorem[numbered=no,style=plain,name=Corollary]{corollary*}

\declaretheorem[numbered=no,style=plain,name=Conjecture]{conjecture*}

\declaretheorem[numbered=no,style=plain,name=Question]{question*}

\declaretheorem[numberwithin=section,numberlike=equation,style=definition]{definition}

\declaretheorem[numbered=no,style=definition,name=Definition]{definition*}
\declaretheorem[numberwithin=section,numberlike=equation,style=definition]{remark}

\declaretheorem[numbered=no,style=definition,name=Remark]{remark*}
\declaretheorem[numberwithin=section,numberlike=equation,style=definition]{notation}

\declaretheorem[numbered=no,style=definition,name=Notation]{notation*}

\declaretheorem[numbered=no,style=definition,name=Axiom]{axiom*}

\declaretheorem[numbered=no,style=definition,name=Construction]{construction*}

\declaretheorem[numbered=no,style=definition,name=Algorithm]{algorithm*}

\declaretheorem[numbered=no,style=definition,name=Property]{property*}

\declaretheorem[numberwithin=section,numberlike=equation,style=definition,qed=$\diamondsuit$]{example}

\declaretheorem[numbered=no,style=definition,name=Example,qed=$\diamondsuit$]{example*}
\usepackage{mystyle}

\newcommand{\define}[1]{\textcolor{magenta}{\emph{#1}}}

\title{An overview of horospherical varieties and coloured fans}
\author{Sean Monahan}

\address{Sean Monahan, Department of Pure Mathematics, University of Waterloo}
\email{sean.monahan314@gmail.com}

\input{./Figures/sample.tikzstyles}
\newcommand{\tikzitfig}[1]{\ctikzfig{./Figures/#1}}

\begin{document}
	
\maketitle
%
\tableofcontents

\section{Introduction}\label{sec:introduction}

An essential step in mathematical research is the testing of one's proposed results on a small, well-understood class of examples. Such a class of examples should have an explicit, hands-on theory that lends itself to easy computations. In the context of varieties in algebraic geometry, a popular choice is the class of toric varieties. These objects admit a very explicit combinatorial description via polyhedral fans, and there is a rich dictionary between the algebro-geometric properties of the varieties and the simple combinatorial properties of the fans. Two key components of this theory are that toric varieties have an open torus orbit and have finitely many torus orbits. 

A downside of toric geometry is that it is fairly narrow. One way to broaden this class is to replace the torus with any reductive algebraic group $G$, while still keeping the key features: having an open $G$-orbit and finitely many $G$-orbits. This leads us to the class of horospherical varieties, which encompasses toric varieties, flag varieties (e.g. Grassmannians), and many other examples. Due to the work of Luna and Vust \cite{luna1983plongements}, horospherical varieties have an explicit combinatorial description using so-called \textit{coloured} fans, which generalizes the correspondence between toric varieties and fans. 

Horospherical varieties were first introduced by Vinberg and Popov \cite{vinberg1972certain}, with subsequent pioneering work by Pauer, Luna, Vust, and Brion \cite{pauer1981normale,luna1983plongements,brion1986espaces}. 
There has been significant interest and activity in horospherical geometry since its inception, particularly in the last 20 years. 
Much is known about the birational geometry of horospherical varieties, e.g. 
their minimal models \cite{brion1993mori,pasquier2018log}; many types of singularities (e.g. canonical, terminal, and log terminal) for horospherical varieties have been characterized \cite{langlois2017singularites,pasquier2016klt,pasquier2015survey}; 
the generalized Mukai conjecture has been proven for horospherical varieties \cite{pasquier2010pseudo}; 
smooth horospherical varieties with Picard number $1$ have been classified and studied \cite{pasquier2009smooth,gonzales2022geometry,hong2016smooth}; 
and the Cox ring and Cox GIT construction (as introduced in \cite{cox1995homogeneous} for toric varieties) were described for horospherical varieties \cite{brion2007total,gagliardi2014cox,gagliardi2019luna}. 
The aforementioned combinatorial theory plays a crucial role in all these works. 

The purpose of this paper is to consolidate many resources to provide a short, yet fairly comprehensive introduction to the combinatorial theory of horospherical varieties. In \cref{subsec:horospherical varieties as a subclass of spherical varieties}, we briefly discuss how horospherical varieties fit into a larger class of spherical varieties. Most of the literature in this area discusses this more general spherical class, so there is a lack of resources that focus solely on the horospherical subclass. We believe that the advantages of working with this subclass are significant from the perspective of seeking an easy-to-work-with class of examples, as discussed above. 

Note that, for the goals of this paper, we do not include proofs of (most) results, and instead we refer the reader to resources where they can find the appropriate proofs; we choose to focus more on examples and developing friendly notation. The inspiration for our new notation is the well-known book \cite{cox2011toric} on toric varieties.

\subsection{Summary of horospherical varieties}\label{subsec:summary of horospherical varieties}

Throughout this paper, let $k$ denote an algebraically closed field of characteristic $0$. For us, all varieties (assumed to be irreducible) and algebraic groups are assumed to be defined over $k$. Fix a connected reductive algebraic group $G$ (e.g. $G=\GL_n$, $G=\SL_n$, or $G=T$ an algebraic torus). We review the necessary algebraic group theory in \cref{sec:review of algebraic groups and roots}. 

We briefly review the definitions of the horospherical objects; see \cref{sec:horospherical objects} for more details. A \textit{horospherical homogeneous space} is a $G$-homogeneous space $G/H$ which has the structure of being a principal torus bundle over a flag variety, i.e. we have a principal bundle $G/H\to G/P$ where $G/P$ is a flag variety and the fibre $P/H$ is an algebraic torus. In this way, one can think of a horospherical homogeneous space as being an extension of a torus by a flag variety. On one extreme, when $G=T$ is a torus, the base $G/P=T/T$ is the trivial flag variety (i.e. a point) and the fibre is the torus $T/H$. On the other extreme, when $H=P$ is parabolic, the fibre $P/P$ is trivial and $G/H=G/P$ is a flag variety. 

These horospherical homogeneous spaces are meant to be the open $G$-orbit of a horospherical variety, and this principal bundle structure is important for the combinatorial description. The following are the key combinatorial ingredients that we obtain from $G/H$; see \cref{sec:colours} for details. First, the torus fibre $P/H$ has an associated lattice of one-parameter subgroups, which we denote by $N$; if $P/H$ is a rank $n$ torus, then $N\cong\Z^n$. Second, the flag variety base $G/P$ has finitely many Schubert divisors, i.e. prime divisors which are invariant under the action of a fixed Borel subgroup of $G$, and these pull back to Borel-invariant divisors $D_\alpha$ in $G/H$ indexed by $\alpha$ in some finite set $\calC$; we call $D_\alpha$ a \textit{colour divisor} and the index $\alpha$ a \textit{colour}. And lastly, each colour $\alpha\in\calC$ determines a \textit{colour point} $u_\alpha$ in $N$ (this uses the order of vanishing valuation $\nu_{D_\alpha}$ for $D_\alpha$; see \cref{subsec:eigenfunctions} for details). In total, we have this lattice $N$ equipped with a ``colour structure", meaning this finite set of colours $\calC$ and the finitely many distinguished colour points in $N$. Note that, when $G=T$ is a torus, there are no colour divisors or colour points since the flag variety base $G/P=T/T$ is trivial. This data of colours can be thought of as the main difference between toric varieties and horospherical varieties. 

A \textit{$G/H$-horospherical variety} (often called a \textit{$G/H$-embedding} in other resources) is a normal $G$-variety $X$ such that there is an open $G$-orbit which is $G$-equivariantly isomorphic to the horospherical homogeneous space $G/H$; so we may view $G/H\subseteq X$ as being open in $X$. When $G=T$ is a torus, we recover toric varieties; and when $H=P$ is parabolic, we recover the flag variety $G/P$. 

Now we review the combinatorial classification of horospherical varieties; see \cref{sec:classification of horospherical varieties} for more details. First, a \textit{coloured cone} on $N$ is a pair $\sigma^c=(\sigma,\calF)$ where $\sigma$ is a strongly convex polyhedral cone in the vector space $N_\R:=N\otimes_\Z \R$ and $\calF\subseteq\calC$ is a set of colours such that $u_\alpha\in\sigma\setminus\{0\}$ for each $\alpha\in\calF$. Then a \textit{coloured fan} $\Sigma^c$ on $N$ is a finite collection of coloured cones which is closed under taking coloured faces of coloured cones, and the intersection of any two coloured cones is a coloured face of both. Let $\calF(\Sigma^c)$ be the union over all colour sets $\calF$ as we range over all $\sigma^c=(\sigma,\calF)$ in $\Sigma^c$. When there are no colours, as in the case of $G=T$ a torus, these coloured cones and coloured fans are the same as ordinary polyhedral cones and fans. 

We can picture a coloured fan as follows: one simply draws an ordinary fan and then highlights the colour points corresponding to colours in $\calF(\Sigma^c)$. The following diagram is an example of a coloured fan $\Sigma^c$ which is drawn on the lattice $\Z^2$; there are two colour points $u_{\alpha_1}$ and $u_{\alpha_2}$, we have $\calF(\Sigma^c)=\{\alpha_1\}$, and there are three maximal coloured cones $\sigma_1^c$, $\sigma_2^c$, and $\sigma_3^c$. 

\tikzitfig{coloured-fan-of-non-toric-non-Q-factorial-SL-3-U-3-variety}

The fundamental theorem of horospherical varieties says that there is an explicit correspondence between horospherical varieties and coloured fans, generalizing the well-known correspondence between toric varieties and fans.

\begin{theorem*}[\cref{thm:classification of horospherical varieties}]
	There is a precise natural bijective correspondence
	\begin{align*}
		\{\text{$G/H$-horospherical varieties}\}/\cong ~&\longleftrightarrow~ \{\text{Coloured fans on } N\}
		\\ X_{\Sigma^c} ~~~&\longleftrightarrow~~~ \Sigma^c
	\end{align*}
	(Here, $\cong$ refers to $G$-equivariant isomorphisms).
\end{theorem*}

One way to think about this correspondence is as follows. Let $X$ be a $G/H$-horospherical variety. Recall that $G/H\to G/P$ is a principal torus bundle over the flag variety $G/P$; the fibre $P/H$ is a torus. We can think of $X$ as a ``compactification" of $G/H$, and \textit{if we can extend} the projection $G/H\to G/P$ to a morphism $X\to G/P$ (see \cref{def:toroidal}), then the fibre $Z$ over $eP$ is a ``compactification" of the torus fibre $P/H$. In this way, $Z$ becomes a toric variety, which has an associated fan $\Sigma$ on the lattice $N$ (since $N$ is the one-parameter subgroup lattice of $P/H$). However, we cannot always extend this map, and the obstruction to extending this map comes from the interaction between the $G$-orbits in $X$ and the colour divisors. It turns out that we can always obtain a canonical resolution $\wt{X}\to X$, where $\wt{X}$ is a $G/H$-horospherical variety which removes this obstruction (this is exactly the decolouration in \cref{ex:decolouration}), and $\wt{X}$ is determined by its toric fibre $\wt{X}\cong G\times^P Z$ (see \cref{subsec:toroidal horospherical varieties}). In total we have the following diagram:
\begin{equation*}
\begin{tikzcd}
	{P/H} &&& Z=\text{toric with fan } \Sigma \\
	{G/H} & X && {\wt{X}=G\times^P Z} \\
	{G/P}
	\arrow[hook, from=1-1, to=1-4]
	\arrow[hook, from=1-1, to=2-1]
	\arrow[hook, from=1-4, to=2-4]
	\arrow[hook, from=2-1, to=2-2]
	\arrow[from=2-1, to=3-1]
	\arrow["\text{resolution}"', from=2-4, to=2-2]
	\arrow[from=2-4, to=3-1]
\end{tikzcd}
\end{equation*}
This resolution is obtained by canonically resolving (in nice situations, blowing up) certain $G$-invariant closed subvarieties of $X$ which are always contained in some of the colour divisors $D_\alpha$. Then the set $\calF(\Sigma^c)$ precisely consists of those $\alpha$'s which determine the resolution. For example, if $\calF(\Sigma^c)=\varnothing$ (e.g. if $G=T$ is a torus and $X$ is a toric variety), then there is no obstruction and the resolution is simply the identity, so $X$ is combinatorially determined by just the fan $\Sigma$ of its toric fibre $Z$. 


The above correspondence gives rise to a comprehensive dictionary between horospherical varieties and coloured fans, through which one can study properties of the varieties by performing simple calculations with the coloured fans. This dictionary has a very familiar feel to the one in toric geometry. Given a $G/H$-horospherical variety $X$ with associated coloured fan $\Sigma^c$, the following are some of the main aspects of this dictionary that we explore: 
\begin{itemize}[leftmargin=2.5em]
	\item There is a bijective correspondence between coloured cones in $\Sigma^c$ and $G$-orbits in $X$, which is inclusion reversing with respect to orbit closures; see \cref{subsec:G-orbits}. 
	\item Horospherical morphisms correspond to linear maps of lattices which are compatible with the structure of the associated coloured fans; see \cref{subsec:horospherical morphisms and maps of coloured fans}. 
	\item $X$ is affine if and only if $\Sigma^c$ is generated by a single coloured cone and $\calF(\Sigma^c)=\calC$, and in this case $k[X]$ can be described explicitly by this coloured cone and certain regular functions on $G$; see \cref{subsec:affine horospherical varieties}.
	\item There is a local breakdown for $X$ using affine horospherical varieties; see \cref{subsec:affine local structure}. 
	\item Local properties such as factoriality, $\Q$-factoriality, and smoothness can be checked using relations between generators of the coloured cones in $\Sigma^c$ and properties of the Dynkin diagram of $G$; see \cref{subsec:factoriality and smoothness}. 
	\item Lastly, we can completely describe the class group and Picard group of $X$ using $\Sigma^c$; see \cref{sec:Weil and Cartier divisors}. In particular, Cartier divisors on $X$ correspond to piecewise linear functions on $\Sigma^c$.
\end{itemize}

\subsection{Horospherical varieties as a subclass of spherical varieties}\label{subsec:horospherical varieties as a subclass of spherical varieties}

We feel that it is important to briefly explain how horospherical varieties fit into the larger world of spherical varieties. If nothing else, this should help the reader navigate other resources on the topic, since most deal with the general spherical class. The following theorem gives a few equivalent definitions of spherical varieties.

\begin{theorem}[Characterizations of spherical]\label{thm:characterizations of spherical}
	Let $X$ be a normal $G$-variety. We say that $X$ is \define{spherical} if any of the following equivalent conditions are satisfied:
	\begin{enumerate}
		\item $X$ has an open orbit under a (or every) Borel subgroup.
		\item $X$ has finitely many orbits under a (or every) Borel subgroup.
		\item Every $G$-equivariant birational model of $X$ has finitely many $G$-orbits. 
	\end{enumerate}

	\begin{proof}
		See \cite[Theorem 2.1.2]{perrin2014geometry}. 
	\end{proof}
\end{theorem}

Horospherical varieties have an open Borel orbit, so they are spherical (see \cref{rmk:horospherical varieties are spherical}). Other examples of spherical varieties include toric varieties, flag varieties, and wonderful varieties. The former two are included in the class of horospherical varieties, but wonderful varieties are not. 

Spherical varieties also enjoy a combinatorial description using coloured fans; this also comes from the paper by Luna and Vust \cite{luna1983plongements} on embeddings of homogeneous spaces. However, this combinatorial theory is more complicated for spherical varieties, and is greatly simplified in the horospherical setting.

One way to characterize the subclass of horospherical varieties within the class of spherical varieties is the following; see \cref{subsec:eigenfunctions} for more details. The group of rational eigenfunctions with respect to the Borel action for a spherical variety $X$ forms a lattice denoted $\Lambda(X)$. Each valuation $\nu$ of $k(X)$ determines a functional on $\Lambda(X)$ by evaluating $\nu$ at eigenfunctions, so $\nu$ yields an element of the dual lattice $\Lambda(X)^\vee$. In fact, the set $\calV(X)$ of $G$-invariant valuations of $k(X)$ embeds into the vector space $\Lambda(X)^\vee_\R$ as a polyhedral cone called the \textit{valuation cone}. 

In the general combinatorial theory, the coloured fans associated to spherical varieties live on this lattice $\Lambda(X)^\vee$, and coloured fans are required to have a nonempty intersection with the valuation cone $\calV(X)$ (see \cite[Section 3]{knop1991luna} for a precise definition of a coloured fan in the general setting). When $X$ is horospherical, this lattice $\Lambda(X)^\vee$ is the lattice $N$ discussed above, and so the coloured fans for horospherical varieties are on this lattice $N$. Moreover, the following result says that, precisely in the horospherical setting, we do not need to keep track of the valuation cone, which simplifies the definition of a coloured fan to the one that we present in this paper. 

\begin{proposition*}[\cref{prop:valuation cone}]
	A spherical $G$-variety $X$ is horospherical if and only if $\calV(X)=\Lambda(X)^\vee_\R$. 
\end{proposition*}

Not only does this simplify the definition of a coloured fan, but it also greatly simplifies the entire combinatorial theory. One can make stronger statements in the dictionary between horospherical varieties and coloured fans than in the general spherical setting (e.g. combinatorial smoothness criteria and local structure theorems).

\subsection{Resources}\label{subsec:resources}

Below are the primary resources used to create this paper:
\begin{itemize}[leftmargin=2.5em]
	\item \cite{knop1991luna} is an introduction to the Luna-Vust theory of spherical varieties. 
	\item \cite{pasquier2009introduction} is an introduction to this topic, which briefly highlights the subclass of horospherical varieties in Section 2.1. 
	\item \cite{pasquier2006thesis} Pasquier's thesis contains much of the basics of horospherical varieties (particularly in Section 2). 
	\item \cite{perrin2014geometry,perrin2018sanya} are detailed survey articles/notes on spherical varieties. 
	\item \cite{timashev2011homogeneous} is an extensive resource on homogeneous spaces and equivariant embeddings, including spherical varieties. Section 28 focuses on horospherical varieties (also known as ``$S$-varieties").
	\item \cite{cox2011toric} is a comprehensive, detailed resource on toric varieties. We try to make the notation in this paper have a familiar feel to their notation. 
\end{itemize}

\subsection{Fixed notation}\label{subsec:fixed notation}

Throughout this paper, $k$ is a fixed algebraically closed field of characteristic $0$. All varieties and algebraic groups are assumed to be defined over this base field $k$. Note that a \textit{variety} for us is an integral separated scheme of finite type over $k$, so in particular varieties are irreducible. 

Furthermore, we fix a connected reductive (linear) algebraic group $G$, a Borel subgroup $B\subseteq G$, and a maximal torus $T\subseteq B$. Let $U$ denote the subgroup of all unipotent elements in $B$, and let $S$ be the set of simple roots of $G$ (relative to $B$ and $T$). We discuss these items in the next section. 

In this paper, instead of considering the $B$-action on a horospherical variety, we choose to consider the opposite $B^-$-action; see \cref{rmk:horospherical varieties are spherical} and the preamble of \cref{sec:colours}. Most sources consider the $B$-action and defined colour divisors as being $B$-invariant, but we choose to define colour divisors as the $B^-$-invariant ones. If the reader prefers to think about $B$-invariant divisors, then they should take $H$ to contain $U^-\subseteq B^-$ and just essentially swap the roles of $B$ and $B^-$ throughout.

\section{Review of algebraic groups and roots}\label{sec:review of algebraic groups and roots}

We briefly review the algebraic group theory that will be used throughout the paper. For more details on algebraic groups, we recommend \cite{humphreys1975linear}. Recall that all varieties and algebraic groups are defined over $k$, an algebraically closed field of characteristic $0$. 

Note that, by ``algebraic group" we always mean ``\textit{linear} algebraic group". For us, $\G_m$ denotes the multiplicative group of the field $k$. Recall that a \define{torus} is an algebraic group of the form $\G_m^n$ for some $n$. 

Throughout this paper, fix a connected reductive algebraic group $G$. Let us recall what it means for $G$ to be reductive. The radical $\Rad(G)$ of $G$ is the identity component of the maximal normal solvable subgroup of $G$, and the unipotent radical $\Rad_u(G)$ of $G$ is the subgroup of $\Rad(G)$ consisting of all unipotent elements; recall that an element of an algebraic group is \define{unipotent} if all of its eigenvalues are $1$. Then $G$ being \define{reductive} means that $\Rad_u(G)$ is trivial. For example, the radical $\Rad(\GL_n)\cong \G_m$ of $\GL_n$ is the subgroup of diagonal matrices where all diagonal entries are equal, so the unipotent radical $\Rad_u(\GL_n)$ is trivial, and hence $\GL_n$ is reductive.

\subsection{Borel subgroups}\label{subsec:Borel subgroups}

Fix a Borel subgroup $B\subseteq G$ (i.e. a maximal solvable closed connected subgroup) and a maximal torus $T\subseteq G$ which is contained in $B$. There are two main examples to keep in mind: (1) when $G=T$ is a torus, in which case $G=B=T$; and (2) when $G=\SL_n$ (or $\GL_n$), in which case we can take $B=B_n$ to be the subgroup of upper triangular matrices and $T=T_n$ to be the subgroup of diagonal matrices. 

Let $U\subseteq B$ denote the subgroup consisting of all unipotent elements in $B$. We can describe $U$ for our two main examples: (1) when $G=T$ is a torus, we have $U=\{1\}$; and (2) when $G=\SL_n$, we have $U=U_n\subseteq B_n$, the subgroup of upper triangular matrices with all $1$'s on the diagonal. This $U$ is a maximal unipotent subgroup of $G$, so $U$ is a solvable closed connected subgroup of $G$. 

Note that all Borel subgroups of $G$ are conjugate to $B$. Moreover, all maximal tori (resp. maximal unipotent subgroups) of $G$ are those of the Borel subgroups of $G$, and all maximal tori (resp. maximal unipotent subgroups) are conjugate, both within $G$ and within their Borel. In particular, the \define{Weyl group} $W=W(G,T):=N_G(T)/T$ acts freely and transitively by conjugation on the set of Borel subgroups of $G$ which contain $T$.

There exists a unique Borel subgroup $B^-\subseteq G$, called the \define{opposite} of $B$ (relative to $T$), such that $B\cap B^-=T$. Note that $B^-$ is conjugate to $B$ by an element of $W$. In our example with $G=\SL_n$, the opposite Borel $B_n^-$ (relative to $T_n$) is the subgroup of lower triangular matrices.

\subsection{Characters}\label{subsec:characters}

For characters and one-parameter subgroups, we use much of the notation and conventions from \cite[Section 1.1]{cox2011toric}. Recall that a \define{character} of an algebraic group $G'$ is an algebraic group homomorphism $\chi:G'\to \G_m$. The dual objects, i.e. algebraic group homomorphisms $\lambda:\G_m\to G'$, are called \define{one-parameter subgroups}. We denote the abelian groups of characters and one-parameter subgroups of $G'$ by $\frakX(G')$ and $\frakX^*(G')$, respectively. We use additive notation for $\frakX(G')$ and $\frakX^*(G')$. However, the group operation on characters (and one-parameter subgroups) is naturally multiplication, so for $m\in\frakX(G')$ we use $\chi^m$ to denote the corresponding multiplicative map $\chi^m:G'\to\G_m$. Note that $m_1+m_2$ corresponds to $\chi^{m_1+m_2}=\chi^{m_1}\chi^{m_2}$. 

We are mostly interested in characters and one-parameter subgroups of diagonalizable (algebraic) groups, i.e. closed subgroups of tori. There is a contravariant equivalence of categories $K\mapsto \frakX(K)$ between diagonalizable groups and finitely generated abelian groups. In particular, this restricts to a contravariant equivalence of categories between tori and lattices (i.e. free abelian groups of finite rank). If $K$ is a torus, then the groups $\frakX(K)$ and $\frakX^*(K)$ are dual lattices. After identifying $\frakX(K)$ and $\frakX^*(K)$ with $\Z^n$ for some $n$, we can view the dual pairing as the usual dot product on $\Z^n$; note that the $\R$-vector spaces $\frakX(K)_\R:=\frakX(K)\otimes_\Z \R$ and $\frakX^*(K)_\R$ can be identified with $\R^n$ with the usual dot product. 

Given a lattice $N$, we let $T_N\cong \G_m^{\rank(N)}$ denote the torus whose one-parameter subgroup lattice is $N$; there is a natural isomorphism $N\otimes_\Z \G_m\isomap T_N$ via $\lambda\otimes t\mapsto \lambda(t)$, where $\lambda:\G_m\to T_N$ denotes a one-parameter subgroup. If $\Phi:N_1\to N_2$ is a $\Z$-linear map of lattices, then the dual map $\Phi^\vee:N_2^\vee\to N_1^\vee$ induces a map of tori $T_{N_1}\to T_{N_2}$ via the contravariant equivalence mentioned above (viewing $N_i^\vee$ as the character lattice of $T_{N_i}$); we denote this map by $\Phi_\Gm:T_{N_1}\to T_{N_2}$. It is easy to verify that $\Phi_\Gm$ is surjective if and only if $\Phi^\vee$ is injective, or equivalently the $\R$-linear extension $\Phi_\R:(N_1)_\R\to (N_2)_\R$ is surjective.

If $H'=\cap_{\chi\in\frakX(G')} \ker(\chi)$, then $G'/H'$ is diagonalizable and $\frakX(G')=\frakX(G'/H')$. Thus $\frakX(G')$ is finitely generated, and if $G'$ is connected, then $\frakX(G')$ is a lattice.

\begin{remark}\label{rmk:character lattice of Borel}
	For our Borel subgroup $B$ and maximal unipotent $U\subseteq B$, we have $\cap_{\chi\in\frakX(B)} \ker(\chi)=U$ and $B/U\cong T$, so $\frakX(B)=\frakX(B/U)=\frakX(T)$.
\end{remark}

\begin{example}[Characters of $B_n\subseteq \SL_n$]\label{ex:characters for SL_n maximal torus}
	Consider $G=\SL_n$ with the Borel subgroup $B_n$ and maximal torus $T_n$. We have $\frakX(B_n)=\frakX(T_n)\cong\Z^{n-1}$. Throughout this paper, we implicitly use the basis $\{e_1,\ldots,e_{n-1}\}$ for $\frakX(T_n)$ where $e_i$ corresponds to the character $\chi_i:T_n\to\G_m$ sending a diagonal matrix to the product of the first $i$ diagonal entries. Said another way, $\chi_i$ is the character $B_n\to\Gm$ sending an upper triangular matrix to the determinant of the upper-left $i\times i$-submatrix. Note that, if we allowed for $i=n$, then $\chi_n$ would be the full determinant, which is trivial in $\SL_n$. This basis will be particularly nice in \cref{sec:colours}.
\end{example}

\subsection{Roots}\label{subsec:roots}

Throughout this subsection we assume that $G$ is semisimple, i.e. $\Rad(G)$ is trivial. For example, $G=\SL_n$ is semisimple, but $G=\GL_n$ or $G=T$ a torus are not semisimple. We explain how to treat the case where $G$ is reductive at the end in \cref{rmk:roots of reductive group}. Let $\frakg$, $\frakb$, and $\frakt$ denote the Lie algebras of $G$, $B$, and $T$, respectively. 

The adjoint representation of $\frakg$ in $GL(\frakg)$ is defined by $\opn{ad}(x):y\mapsto [x,y]$ for $x,y\in\frakg$, where $[\cdot,\cdot]$ is the Lie bracket. An nonzero element $\alpha\in\frakX(T)$ is called a \define{root} of $(G,T)$ if the space of $\alpha$-eigenvectors for the adjoint representation
\begin{align*}
	\frakg_\alpha := \{x\in\frakg : [t,x]=\chi^\alpha(t)x, ~\forall t\in\frakt\}
\end{align*}
is nonzero. Let $R=R(G,T)$ denote the set of roots of $(G,T)$, which is a finite set. Also let $R^+:=\{\alpha\in R:\frakg_\alpha\subseteq\frakb\}$ denote the subset of \define{positive roots} (relative to $B$), and let $R^-:=R\setminus R^+$ denote the subset of \define{negative roots}. Note that $R^-=\{-\alpha:\alpha\in R^+\}$. Then
\begin{align*}
	\frakg = \frakt \oplus\bigoplus_{\alpha\in R} \frakg_\alpha = \underbrace{\frakt \oplus\bigoplus_{\alpha\in R^+} \frakg_\alpha}_{\frakb} \oplus\bigoplus_{\alpha\in R^-} \frakg_\alpha
\end{align*}
and $\dim_k(\frakg_\alpha)=1$ for each $\alpha\in R$. 

Let $S\subseteq R^+$ denote the set of \define{simple roots}, i.e. positive roots which cannot be written as a sum of two positive roots. Every positive (resp. negative) root is a unique $\Z_{\geq 0}$-linear (resp. $\Z_{\leq 0}$-linear) combination of elements from $S$; although, not \textit{every} $\Z$-linear combination of elements from $S$ yields a root. 
Since $G$ is semisimple, one can show that $S$ is an $\R$-basis for $\frakX(T)_\R=\frakX(T)\otimes_\Z \R$, and in particular $\# S = \dim_\R(\frakX(T)_\R) = \dim(T)$. For $I\subseteq S$, we let $R_I$ denote the subset of $R$ generated by $I$, i.e. $R_I$ consists of the roots which are $\Z$-linear combinations of elements in $I$. Also let $R_I^+:=R_I\cap R^+$ and $R_I^-:=R_I\cap R^-$. 

The Weyl group $W=N_G(T)/T$ acts on $\frakX(T)$ by $w\cdot m$ being the character defined as $(w\cdot \chi^m):=\chi^m(w^{-1}tw)$ for all $w\in W$, $m\in \frakX(T)$, and $t\in T$. This action faithfully permutes the roots, so $W$ is a subgroup of the symmetric group of the set $R$. In particular, $W$ is finite. 

Given $\alpha\in R$, let $G_\alpha$ denote the subgroup of $G$ whose Lie algebra is $\frakt\oplus\frakg_\alpha\oplus\frakg_{-\alpha}$. Each $G_\alpha$ contains $T$, and the Weyl group of $G_\alpha$, i.e. $N_{G_\alpha}(T)/T$, has order $2$. Therefore, each $\alpha\in R$ determines a unique nontrivial element of $W$ which has a representative in $N_{G_\alpha}(T)$; we call this element the \define{reflection} associated to $\alpha$, and we denote it $r_\alpha\in W$. Note that $r_\alpha=r_{-\alpha}$ and $r_\alpha\cdot\alpha=-\alpha$. 

The Weyl group $W$ is generated by $\{r_\alpha:\alpha\in R\}$, and is even generated by $\{r_\alpha:\alpha\in S\}$. For any $I\subseteq S$, we let $W_I$ denote the subgroup of $W$ generated by $\{r_\alpha:\alpha\in I\}$. In particular, $W_\varnothing$ is trivial, and $W_S=W$. 

\begin{example}[Roots of $\SL_n$]\label{ex:roots of SL_n}
	Consider $G=\SL_n$. Then $\mathfrak g=\mathfrak {sl}_n$ is the traceless $n\times n$ matrices (over $k$), and $[x,y]=xy-yx$. Consider $T=T_n$ and $B=B_n$. Throughout this paper, we implicitly use the following notation for the roots and simple roots of $\SL_n$. 
	
	If $\alpha\in R$ is a root, then there exists $0\neq x\in\mathfrak g$ such that $tx-xt=\chi^\alpha(t)x$ for all $t\in\mathfrak t$. Since $t\in\mathfrak t$ is diagonal, we can write $t=\diag(t_1,\ldots,t_n)$ for some $t_1,\ldots,t_n\in k$. Hence, for all $1\leq i,j\leq n$ we have
	\begin{align*}
		\chi^\alpha(t)x_{i,j} = (tx-xt)_{i,j} = t_ix_{i,j}-t_jx_{i,j} = (t_i-t_j)x_{i,j}
	\end{align*}
	(the subscript $i,j$ is denoting the $(i,j)$-entry). Thus $\chi^\alpha(t)=t_i-t_j$ as a map $\frakt\to k$; or if we view $\chi^\alpha$ as a character of $T$ then $\chi^\alpha(t)=t_it_j^{-1}$. Therefore, the roots are
	\begin{align*}
		R = \{\alpha_{i,j} : 1\leq i\neq j\leq n\}, \qquad \chi^{\alpha_{i,j}}:T\to \Gm ~~~ \diag(t_1,\ldots,t_n)\mapsto t_it_j^{-1}.
	\end{align*}
	
	It is easy to check that $\mathfrak g_{\alpha_{i,j}}$ is the subalgebra generated by $E_{i,j}$ (i.e. the standard basis matrix with $1$ in the $(i,j)$ entry and $0$ elsewhere). So we have $\mathfrak g_{\alpha_{i,j}}\subseteq\mathfrak b$ if and only if $j>i$. Therefore, the positive roots are $R^+=\{\alpha_{i,j}:1\leq i<j\leq n\}$, and it is easy to see that the simple roots are $S=\{\alpha_i:1\leq i<n\}$ where $\alpha_i:=\alpha_{i,i+1}$. 
	
	Note that $N_{\SL_n}(T)$ is the group of generalized permutation matrices in $\SL_n$, so the Weyl group $W(\SL_n,T)=N_{\SL_n}(T)/T$ is isomorphic to $\frakS_n$ (the symmetric group). The reflection $r_{\alpha_i}$ for $\alpha_i\in S$ corresponds to the transposition $(i ~~ i+1)\in\frakS_n$. 
\end{example}

\begin{remark}[Roots of reductive $G$]\label{rmk:roots of reductive group}
	If $G$ is reductive, then $[G,G]$ (the commutator) and $G/\Rad(G)$ are semisimple, and there is a natural bijective correspondence between the roots of $G$, $[G,G]$, and $G/\Rad(G)$. The key difference in the root data between $G$ and these associated semisimple groups is in the rank of the maximal torus. For example, the roots of $G=\GL_n$ are the same as the roots of $[G,G]=G/\Rad(G)=\SL_n$, but the maximal torus of $\GL_n$ has rank $n$ and the maximal torus of $\SL_n$ has rank $n-1$. As a result, the simple roots $S$ may not span $\frakX(T)_\R$; in $G=\GL_n$ we have $\frakX(T)_\R\cong\R^n$ but $S$ spans an $(n-1)$-dimensional subspace. 
\end{remark}

\begin{example}[Roots of torus]\label{ex:roots of torus}
	Using \cref{rmk:roots of reductive group}, we see that a torus has a trivial root system, i.e. the root set is empty and the Weyl group is trivial. 
\end{example}

\begin{remark}\label{rmk:ss and torus isogeny}
	If $G$ is reductive, then there always exists an isogeny (i.e. a surjective algebraic group homomorphism with finite kernel) $G^{ss}\times T_0\to G$ where $G^{ss}$ is semisimple simply connected (i.e. every isogeny from a connected algebraic group to $G^{ss}$ is an isomorphism) and $T_0$ is a torus. We can choose a Borel subgroup and a maximal torus of $G^{ss}\times T_0$ that map onto $B$ and $T$, respectively. Then the roots of $G^{ss}\times T_0$ (relative to this Borel and maximal torus) are the same as the roots of $G$, and the maximal tori have the same rank. 
\end{remark}

\subsection{Parabolic subgroups}\label{subsec:parabolic subgroups}

We no longer restrict to the assumption that $G$ is semisimple, i.e. $G$ can be reductive. A \define{parabolic subgroup} of $G$ is a subgroup of $G$ which contains a Borel subgroup of $G$. We say that a parabolic subgroup $P\subseteq G$ is \define{standard} if it contains our fixed $B$. If $P$ is a parabolic subgroup of $G$, then $P$ is closed and connected, we have $N_G(P)=P$, and the homogeneous space $G/P$ is a smooth projective variety. Any variety of the form $G/P$ with $P$ parabolic is called a \define{flag variety}. 

For $I\subseteq S$, set $P_I:=BW_IB\subseteq G$. For example, $P_\varnothing=B$ and $P_S=G$. Each $P_I$ is a standard parabolic subgroup of $G$ whose Lie algebra is
\begin{align*}
	\frakp_I = \underbrace{\frakt \oplus\bigoplus_{\alpha\in R^+} \frakg_\alpha}_{\frakb} \oplus\bigoplus_{\alpha\in R_I^-} \frakg_\alpha.
\end{align*}

A \define{Levi subgroup} of an algebraic group $G'$ is a closed subgroup $L'\subseteq G'$ such that the natural map $L'\to G'/\Rad_u(G')$ is an isomorphism. Note that, if $G'$ is connected, then $L'$ is a connected reductive algebraic group. Every standard parabolic subgroup $P\subseteq G$ has a unique Levi subgroup which contains $T$; we call this the \define{standard Levi subgroup} of $P$. If $P=P_I$ for some $I\subseteq S$, then we denote this unique Levi subgroup by $L_I$. For example, $L_\varnothing=T=B/U$ and $L_S=G=G/\Rad_u(G)$. 

Since $L_I$ is connected and reductive, we can consider its roots. Relative to the Borel subgroup $B\cap L_I\subseteq L_I$ and the maximal torus $T\subseteq L_I$, the set of roots for $L_I$ is $R_I$, the set of simple roots is $I$, and the Weyl group of $L_I$ is $W_I$. The Lie algebra of $L_I$ is
\begin{align*}
	\frakl_I = \frakt \oplus\bigoplus_{\alpha\in R_I} \frakg_\alpha = \underbrace{\frakt \oplus\bigoplus_{\alpha\in R_I^+} \frakg_\alpha}_{\frakb\cap\frakl_I} \oplus\bigoplus_{\alpha\in R_I^-} \frakg_\alpha.
\end{align*}

\begin{proposition}[Classification of parabolic subgroups]\label{prop:classification of parabolic subgroups}
	The map $I\mapsto P_I$ is a bijection between subsets $I\subseteq S$ and standard parabolic subgroups of $G$. This bijection satisfies the following properties:
	\begin{enumerate}
		\item If $P_I$ and $P_J$ are conjugate, then $P_I=P_J$.
		\item $P_I\subseteq P_J$ (resp. $L_I\subseteq L_J$) if and only if $I\subseteq J$.
		\item $P_{I\cap J}=P_I\cap P_J$ and $L_{I\cap J}=L_I\cap L_J$. 
	\end{enumerate}
	
	\begin{proof}
		See \cite[Theorem 29.3]{humphreys1975linear}. 
	\end{proof}
\end{proposition}

\begin{example}[Parabolics of $\SL_n$]\label{ex:parabolic and Levi subgroups of SL_n}
	Consider $G=\SL_n$ with the Borel subgroup $B_n$. We describe the standard parabolic and Levi subgroups of $\SL_n$ using \cref{prop:classification of parabolic subgroups}. Recall that $S=\{\alpha_1,\ldots,\alpha_{n-1}\}$. 
	
	Since $B_n$ is the subgroup of upper triangular matrices, the parabolics of $\SL_n$ containing $B_n$ are subgroups consisting of block upper triangular matrices. For the parabolic $P_I$, the set $I\subseteq S$ describes the blocks on the diagonal. Specifically, if $S\setminus I=\{\alpha_{i_1},\ldots,\alpha_{i_s}\}$ with $i_1<i_2<\cdots<i_s$, and we set $i_0:=0$ and $i_{s+1}:=n$, then $P_I$ is the subgroup of $\SL_n$ consisting of block upper triangular matrices where the blocks on the diagonal have sizes $i_1-i_0,i_2-i_1,\ldots,i_{s+1}-i_s$. 
	
	The standard Levi subgroups consist of block diagonal matrices: the set $I$ describes the blocks on the diagonal for $L_I$. Specifically, $L_I$ is the subgroup of block diagonal matrices where the blocks have sizes $i_1-i_0,i_2-i_1,\ldots,i_{s+1}-i_s$. 
	
	As some specific examples in $\SL_4$, we have
	\begin{align*}
		P_{\{\alpha_1,\alpha_2\}} = \begin{bmatrix}
			* & * & * & * \\ * & * & * & * \\ * & * & * & * \\ 0 & 0 & 0 & *
		\end{bmatrix},
		\quad P_{\{\alpha_1,\alpha_3\}} = \begin{bmatrix}
			* & * & * & * \\ * & * & * & * \\ 0 & 0 & * & * \\ 0 & 0 & * & *
		\end{bmatrix},
		\quad P_{\{\alpha_2\}} = \begin{bmatrix}
			* & * & * & * \\ 0 & * & * & * \\ 0 & * & * & * \\ 0 & 0 & 0 & *
		\end{bmatrix}
	\end{align*}
	and
	\begin{align*}
		L_{\{\alpha_1,\alpha_2\}} = \begin{bmatrix}
			* & * & * & 0 \\ * & * & * & 0 \\ * & * & * & 0 \\ 0 & 0 & 0 & *
		\end{bmatrix},
		\quad L_{\{\alpha_1,\alpha_3\}} = \begin{bmatrix}
			* & * & 0 & 0 \\ * & * & 0 & 0 \\ 0 & 0 & * & * \\ 0 & 0 & * & *
		\end{bmatrix},
		\quad L_{\{\alpha_2\}} = \begin{bmatrix}
			* & 0 & 0 & 0 \\ 0 & * & * & 0 \\ 0 & * & * & 0 \\ 0 & 0 & 0 & *
		\end{bmatrix}.
	\end{align*}
	
	With this description, one can understand the flag varieties $\SL_n/P_I$. In particular, when $I=\varnothing$, we have the prototypical flag variety $\SL_n/B$ which is the space of complete flags in the vector space $k^n$. Another key example is the Grassmannian $\opn{Gr}(\ell,k^n)$ of $\ell$-dimensional subspaces of $k^n$, which we can write as $\SL_n/P_I$ where $S\setminus I=\{\alpha_\ell\}$. For example, $\SL_4/P_{\{\alpha_1,\alpha_3\}}$ is $\opn{Gr}(2,k^4)$. 
\end{example}

\begin{remark}\label{rmk:characters of parabolic}
	If $P$ is a standard parabolic subgroup of $G$, then we have a surjection of tori $T=B/U\to P/(\cap_{\chi\in\frakX(P)} \ker(\chi))$, so we always have an inclusion of character lattices $\frakX(P)\subseteq \frakX(T)$; compare with \cref{rmk:character lattice of Borel}. 
\end{remark}

A final word about parabolic subgroups is that, for any parabolic $P\subseteq G$, we can consider its \define{opposite} parabolic subgroup $P^-\subseteq G$ containing $B^-$: if $P=P_I$ is standard, then $P^-=P_I^-$ is the unique parabolic subgroup of $G$ (standard relative to $B^-$) which has the same standard Levi subgroup $L_I$ as $P_I$ and satisfies $P_I\cap P_I^-=L_I$. In \cref{ex:parabolic and Levi subgroups of SL_n} for $G=\SL_n$, the opposite parabolic $P_I^-$ is the appropriate subgroup of block \textit{lower} triangular matrices, where the blocks are described by $I$ as discussed in this example.

\subsection{Bruhat decomposition}\label{subsec:Bruhat decomposition}

To finish this section, we discuss the Borel orbits of flag varieties using the well-known Bruhat decomposition. This is very useful in describing the Borel orbits of horospherical homogeneous spaces, which is a key ingredient in the study of horospherical (and spherical) varieties. 

\begin{proposition}[Bruhat decomposition]\label{prop:Bruhat decomposition}
	If $P_I\subseteq G$ is a standard parabolic subgroup (with $I\subseteq S$), then we have
	\begin{align*}
		G = \bigsqcup_{w\in W/W_I} B^-wP_I
	\end{align*}
	(here $w$ is an element of $W/W_I$ and is also a representative in $G$). Furthermore, $B^-wP_I=B^-w'P_I$ if and only if $w=w'$ in $W/W_I$. 
	
	\begin{proof}
		See \cite[Theorem 28.3]{humphreys1975linear}. 
	\end{proof}
\end{proposition}

This proposition tells us the $B^-$-orbit structure of flag varieties because we have
\begin{align}\label{eq:Schubert cell decomposition}
	G/P_I = \bigsqcup_{w\in W/W_I} B^-wP_I/P_I.
\end{align}
The $B^-wP_I/P_I$ appearing above are called the \define{Schubert cells} of the flag variety $G/P_I$,\footnote{Based on our convention of using the opposite Borel $B^-$ in the cell $B^-wP_I$, these are typically known as the \textit{opposite} cells; e.g. see \cite[Definition 1.2.4]{brion2004lectures}} which are exactly the $B^-$-orbits of $G/P_I$. The closures of the Schubert cells are called \define{Schubert varieties}; these are closed subvarieties of $G/P_I$. 

If $w\in W/W_I$ is written as a product of generators $w=r_{\alpha_1}\cdots r_{\alpha_\ell}$ with $\alpha_i\in S$ and $\ell$ is minimal, then we call $\ell=\ell(w)$ the \define{length} of $w$; see \cite[Section 29.2]{humphreys1975linear} or \cite[pg. 8]{brion2004lectures} for more details. The identity $w=e$ in $W/W_I$ is the unique shortest element, i.e. $\ell(e)=0$. The codimension of the Schubert cell $B^-wP_I/P_I$ (or of the corresponding Schubert variety) in $G/P_I$ is equal to $\ell(w)$. Thus, the only Schubert cell of full dimension (i.e. dimension $\dim(G/P_I)$) is when $w=e$: $B^-P_I/P_I$; we call this the \define{big Schubert cell}. This big Schubert cell is an open $B^-$-orbit in $G/P_I$, and the corresponding Schubert variety is $G/P_I$.

\section{Horospherical objects}\label{sec:horospherical objects}

Before studying the combinatorial description of horospherical varieties, we use this section to define the various horospherical terms: horospherical subgroups, horospherical homogeneous spaces, horospherical varieties, and horospherical morphisms.

\subsection{Horospherical subgroups}\label{subsec:horospherical subgroups}

We first define horospherical subgroups, which act as the general stabilizer groups for horospherical varieties. 

\begin{definition}[Horospherical subgroup]\label{def:horospherical subgroup}
	A closed subgroup $H\subseteq G$ is \define{horospherical} if it contains a maximal unipotent subgroup of $G$. If $H$ contains our fixed $U$, then we say that $H$ is \define{standard}. 
\end{definition}

\begin{example}[Basic examples of horospherical subgroups]\label{ex:basic examples of horospherical subgroups}
	\mbox{}\par
	\begin{itemize}[leftmargin=2.5em]
		\item If $G=T$ is a torus, then $U=\{1\}$, so any closed subgroup of $T$ is horospherical. 
		\item For any $G$, maximal unipotent subgroups are horospherical.
		\item Any parabolic subgroup $P\subseteq G$ is horospherical. 
		\item If $G_i$ are two connected reductive algebraic groups and $H_i\subseteq G_i$ are horospherical subgroups for $i=1,2$, then $H_1\times H_2$ is a horospherical subgroup of $G_1\times G_2$. \qedhere
	\end{itemize}
\end{example}

The following proposition describes the normalizer of a horospherical subgroup, which is used heavily in our study of horospherical varieties. 

\begin{proposition}[Normalizer of horospherical subgroup]\label{prop:normalizer of horospherical is parabolic}
	Let $H\subseteq G$ be a (resp. standard) horospherical subgroup. Then $N_G(H)$ is a (resp. standard) parabolic subgroup of $G$ and $N_G(H)/H$ is a torus. 
\end{proposition}

For a horospherical subgroup $H\subseteq G$, we call $P=N_G(H)$ the \define{associated parabolic} to $H$, and we call $P/H$ the \define{associated torus} (to $H$ or $G/H$). Let $N=N(G/H):=\frakX^*(P/H)$ be the one-parameter subgroup lattice of the torus $P/H$ associated to $G/H$; so $N^\vee=\frakX(P/H)$ is the dual lattice of characters. As in \cref{subsec:characters}, we can write $P/H=T_N$.

Assuming $H$ is standard, one can quickly verify that $P=TH$ and thus $T_N=T/(T\cap H)$; recall that $T\subseteq B$ is the maximal torus. So we have a surjection of tori and a corresponding inclusion of character lattices:
\begin{align}\label{eq:surjection from maximal torus to associated torus}
	T \surjmap T_N=T/(T\cap H), \qquad N^\vee=\frakX(T_N) \hookrightarrow \frakX(T).
\end{align}

\begin{example}[Basic examples of associated parabolics and tori]\label{ex:basic examples of associated parabolics and tori}
	Let $H\subseteq G$ be a horospherical subgroup. For the examples in \cref{ex:basic examples of horospherical subgroups}, we have:
	\begin{itemize}[leftmargin=2.5em]
		\item If $G=T$ is a torus, then the associated parabolic to $H$ is $T$, and the associated torus is $T/H$. 
		\item The associated parabolic to the maximal unipotent $U$ is the Borel $B$, and the associated torus is the maximal torus $B/U=T$.
		\item If $H=P\subseteq G$ is parabolic, then the associated parabolic is $P$ itself, and the associated torus is the trivial group $P/P$. 
		\item If $G_i$ are two connected reductive algebraic groups and $H_i\subseteq G_i$ are horospherical subgroups for $i=1,2$, then the associated parabolic to $H_1\times H_2$ is $P_1\times P_2$, where $P_i$ is the associated parabolic to $H_i$, and the associated torus is $(P_1\times P_2)/(H_1\times H_2)=(P_1/H_1)\times(P_2/H_2)$. \qedhere
	\end{itemize}
\end{example}

Now we give a classification of horospherical subgroups. Consider a pair $(I,M)$ where $I\subseteq S$ and $M\subseteq\frakX(P_I)\subseteq\frakX(T)$. The pair $(I,M)$ determines a standard horospherical subgroup $H_{(I,M)}$ defined below whose associated parabolic is $P_I$:
\begin{align*}
	H_{(I,M)} := \bigcap_{\chi\in M} \ker(\chi).
\end{align*}
Indeed, this is a standard horospherical subgroup because every character $\chi:P_I\to\G_m$ is trivial on $U$, so $H_{(I,M)}\supseteq U$. Notice that, if $M$ and $M'$ generate the same sublattice of $\frakX(P_I)$, then $H_{(I,M)}=H_{(I,M')}$. So we can just look at pairs $(I,M)$ where $I\subseteq S$ and $M$ is a sublattice of $\frakX(P_I)\subseteq \frakX(T)$. The following proposition shows that every standard horospherical subgroup arises in this way.

\begin{proposition}[Classification of horospherical subgroups]\label{prop:classification of horospherical subgroups}
	There is a bijection
	\begin{align*}
		\{(I,M):I\subseteq S,~ M\subseteq\frakX(P_I) \text{ a sublattice}\} &\map \{\text{Standard horospherical subgroups of $G$}\}
		\\ (I,M) &\longmapsto H_{(I,M)}.
	\end{align*}
	In particular, for a fixed standard parabolic subgroup $P\subseteq G$, this restricts to an inclusion reversing bijection
	\begin{align*}
		\{\text{Sublattices } M\subseteq\frakX(P)\} &\map \{\text{Horospherical subgroups with associated parabolic $P$}\}.
	\end{align*}
	
	\begin{proof}
		See \cite[Proposition 2.4]{pasquier2008varietes}. 
	\end{proof}
\end{proposition}

\begin{remark}[Minimal horospherical subgroups]\label{rmk:minimal horospherical subgroups}
	Given a standard parabolic subgroup $P=P_I$ of $G$ for $I\subseteq S$, there is a unique minimal horospherical subgroup $H_{(I,\frakX(P))}$ whose associated parabolic is $P$. This is equal to the commutator $[P,P]$. For example, $[B,B]=U$. 
\end{remark}

\begin{example}\label{ex:horospherical subgroups of SL_2}
	Consider $G=\SL_2$. Then $S=\{\alpha\}$ and $\frakX(T_2)\cong\Z$. We use \cref{prop:classification of horospherical subgroups} to determine all standard horospherical subgroups of $\SL_2$.
	
	If $I=S$, then $\frakX(P_S)=\frakX(\SL_2)$ is trivial, so $H_{(S,\{0\})}=\SL_2$ is the only (standard) horospherical subgroup of $\SL_2$ with associated parabolic $P_S=\SL_2$. 
	
	The other possibility for the subset $I\subseteq S$ is $I=\varnothing$. In this case, $\frakX(P_\varnothing)=\frakX(B_2)=\frakX(T_2)\cong\Z$. Therefore, the standard horospherical subgroups of $\SL_2$ with associated parabolic $P_\varnothing=B_2$ are $H_{(\varnothing,M)}$ where $M\subseteq\Z$ is a sublattice. Since sublattices of $\Z$ are of the form $a\Z$ for $a\in\Z_{\geq 0}$, we can write these horospherical subgroups explicitly as 
	\begin{align*}
		H_{(\varnothing,a\Z)} = \left\{\begin{bmatrix}\zeta_a&*\\0&\zeta_a^{-1}\end{bmatrix} : \zeta_a^a=1\right\}.
	\end{align*}
	Note that, when $a=0$, this yields the Borel $P_\varnothing=B_2$. 
\end{example}

\subsection{Horospherical homogeneous spaces}\label{subsec:horospherical homogeneous spaces}

If $H\subseteq G$ is a horospherical subgroup, then we call the homogeneous space $G/H$ a \define{horospherical homogeneous space}. Since $G/H$ is the same (up to $G$-equivariant isomorphism) after conjugating $H$, we may assume that $H$ is standard; going forward, we always assume that this is the case for horospherical homogeneous spaces. 

\begin{remark}[$G$-equivariant automorphisms of $G/H$]\label{rmk:equivariant automorphisms of G/H}
	Let $G/H$ be a horospherical homogeneous space. The associated torus $P/H$, with $P=N_G(H)\supseteq B$, is naturally isomorphic to the group of $G$-equivariant automorphisms of $G/H$, denoted $\Aut^G(G/H)$. 
\end{remark}

Given a horospherical homogeneous space $G/H$, we have an \define{associated flag variety} $G/P$, where $P=N_G(H)$ is the associated parabolic, and we have a natural projection map $G/H\to G/P$. The fibre over $eP\in G/P$ is the associated torus $P/H$, and every other fibre is isomorphic to this one. That is, $G/H\to G/P$ is a principal $P/H$-bundle. The following proposition shows that this structure characterizes horospherical homogeneous spaces. 

\begin{proposition}[Characterization of horospherical homogeneous spaces]\label{prop:characterization of horospherical homogeneous spaces}
	A homogeneous space $G/H$ is horospherical if and only if it is a principal torus bundle over a flag variety. 
	
	\begin{proof}
		The forward direction is shown above. Conversely, if there exists a flag variety $G/P$ such that $G/H\to G/P$ is a principal torus bundle, then we may assume that $H\subseteq P$ and this map is the natural projection. The fibre $P/H$ is a torus, so $H$ is normal in $P$ and it follows that $H$ is horospherical with associated parabolic $P$. 
	\end{proof}
\end{proposition}

\begin{example}[Primary horospherical homogeneous space]\label{ex:primary horospherical homogeneous space}
	Since $U$ is the smallest standard horospherical subgroup of $G$, we see that $G/U$ is the ``biggest" horospherical $G$-homogeneous space. What we mean by this is that, for any other horospherical homogeneous space $G/H$, we have a $G$-equivariant surjection $G/U\to G/H$, i.e. the natural projection induced by $U\subseteq H$. We call $G/U$ the \define{primary horospherical homogeneous space} (with respect to $G$). 
	
	Now consider a standard parabolic subgroup $P\subseteq G$. As in \cref{rmk:minimal horospherical subgroups}, $[P,P]$ is the smallest horospherical subgroup with associated parabolic $P$, so $G/[P,P]$ is the ``biggest" horospherical $G$-homogeneous space with associated flag variety $G/P$. That is, for any other horospherical homogeneous space $G/H$ with associated flag variety $G/P$, we have a $G$-equivariant surjection $G/[P,P]\to G/H$, i.e. the natural projection induced by $[P,P]\subseteq H$. We call $G/[P,P]$ the \define{$P$-primary horospherical homogeneous space} (with respect to $G$). 
\end{example}

Again, consider a horospherical homogeneous space $G/H$ with associated flag variety $G/P$; so $P=N_G(H)\supseteq B$. Recall that the Bruhat decomposition gives the $B^-$-orbit structure of $G/P$; see \cref{eq:Schubert cell decomposition}. This can be used to obtain the $B^-$-orbit structure for $G/H$:
\begin{align}\label{eq:B-orbit decomposition of G/H}
	G/H = \bigsqcup_{w\in W/W_I} B^-wH/H
\end{align}
(note that $B^-wP/H=B^-wH/H$ since $P=TH$ and $wT=w$). That is, the $B^-$-orbits of $G/H$ are the preimages of the $B^-$-orbits of $G/P$ under the projection $G/H\to G/P$. 

\begin{remark}[Open $B^-$-orbit of $G/H$]\label{rmk:horospherical homogeneous space has open B-orbit}
	Every horospherical homogeneous space has an affine open $B^-$-orbit, namely the big cell $B^-H/H$, i.e. the $B^-$-orbit of the point $eH\in G/H$. One way to see that $B^-H/H$ is affine is using \cref{subsec:simple open cover}: we have $B^-H/H=X_{G/H,B^-}$. 
\end{remark}

\subsection{Horospherical varieties}\label{subsec:horospherical varieties}

Now we define the main algebro-geometric object of this paper: horospherical varieties. Let $G/H$ be a horospherical homogeneous space. 

\begin{definition}[Horospherical variety]\label{def:horospherical variety}
	A \define{$G/H$-horospherical variety} is a normal $G$-variety $X$ together with a base point $x\in X$ such that the orbit $G\cdot x\subseteq X$ is open and the stabilizer $G_x$ is the horospherical subgroup $H$. 
\end{definition}

We call $X$ a \define{horospherical $G$-variety} if it is a $G/H$-horospherical variety for some horospherical $H\subseteq G$. In a $G/H$-horospherical variety $X$, the open orbit $G\cdot x$ is $G$-equivariantly isomorphic to the horospherical homogeneous space $G/H$. So we regard $G/H\subseteq X$ as the open orbit; in this notation, the base point is simply $eH\in G/H\subseteq X$. 

\begin{remark}[Horospherical varieties are spherical]\label{rmk:horospherical varieties are spherical}
	By \cref{rmk:horospherical homogeneous space has open B-orbit}, every $G/H$-horospherical variety $X$ is spherical since the $B^-$-orbit $B^-H/H\subseteq X$ of the base point is open. 
	
	We briefly discuss spherical varieties in \cref{subsec:eigenfunctions}, and we see an important way of viewing horospherical varieties as a subclass of spherical varieties (\cref{prop:valuation cone}). 
\end{remark}

\begin{proposition}[Rationality of horospherical varieties]\label{prop:rationality of horospherical varieties}
	Every horospherical variety is rational and has rational singularities (hence is Cohen-Macaulay). 
	
	\begin{proof}
		See \cite[Corollary 2.1.3, Corollary 2.3.4]{perrin2014geometry} (the latter uses $\charac(k)=0$). 
	\end{proof}
\end{proposition}

\begin{example}[Basic examples of horospherical varieties]\label{ex:basic examples of horospherical varieties}
	\mbox{}\par
	\begin{itemize}[leftmargin=2.5em]
		\item Any toric variety with open torus $T$ is a $T/\{1\}$-horospherical variety. 
		\item $G/H$ is itself a $G/H$-horospherical variety.
		\item As a special case of the previous example, any flag variety $G/P$ is the \textit{unique} $G/P$-horospherical variety.  
		\item If $X_i$ are $G_i/H_i$-horospherical varieties for $i=1,2$, then $X_1\times X_2$ is a $(G_1/H_1)\times (G_2/H_2)=(G_1\times G_2)/(H_1\times H_2)$-horospherical variety. \qedhere
	\end{itemize}
\end{example}

\begin{example}\label{ex:A^2 and P^2 are SL_2/U_2-varieties} 
	Consider $G=\SL_2$. There is an $\SL_2$-action on $\A^2$ via matrix multiplication. Choosing the base point $(1,0)\in\A^2$, the orbit $\SL_2\cdot(1,0)=\A^2\setminus\{0\}$ is open in $\A^2$, and the stabilizer is $U_2$. Therefore, $\A^2$ is a $\SL_2/U_2$-horospherical variety. 
	
	We have an open embedding $\A^2\hookrightarrow\P^2$ using the affine chart $\{(1:y:z)\}\subseteq\P^2$. The $\SL_2$-action on $\A^2$ above yields an action on $\P^2$, and with the base point $(1:1:0)\in\P^2$ we can view $\P^2$ as a $\SL_2/U_2$-horospherical variety using the previous paragraph. 
\end{example}

The varieties in \cref{ex:A^2 and P^2 are SL_2/U_2-varieties} are viewed as $\SL_2/U_2$-horospherical varieties, but both are intrinsically toric varieties so they can be studied using toric geometry. To get a sense of which varieties are horospherical but not toric, we prove the following proposition; compare with \cite[Proposition 5.14]{monahan2025horospherical-stack}. 

\begin{proposition}[Non-toric condition]\label{prop:non-toric condition}
	Let $X$ be a $G/H$-horospherical variety. If the associated flag variety $G/P$ is not a product of projective spaces, then $X$ is not a toric variety. Moreover, when $X$ is affine or complete (e.g. projective), the converse is also true. 
	
	\begin{proof}
		Without loss of generality, we may assume that $\scrO_X^*(X)=k^*$. Indeed, this uses \cref{rmk:removing torus factors} in which we can write $X=X_1\times T_1$ where $\scrO_{X_1}^*(X_1)=k^*$ and $T_1$ is a torus, and it is clear that $X$ is a toric variety if and only if $X_1$ is a toric variety. By reducing to this case, we can consider the Cox ring of $X$ (see \cite{gagliardi2014cox} for the horospherical case). 
		
		We prove the contrapositive, so suppose that $X$ is a toric variety. Then the Cox ring $\Cox(X)$ (see \cite{gagliardi2014cox} for the horospherical case) is a polynomial ring over $k$. By \cite[Theorem 3.8]{gagliardi2014cox}, this implies that $\Cox(G/P)$ is a polynomial ring over $k$. Thus, $G/P$ is a toric variety by \cite[Proposition 6.1]{gagliardi2017generalized} (since $G/P$ is complete), so $G/P$ is a product of projective spaces by \cite[Theorem 1]{thomsen1997affinity}. 
		
		Now we tackle the converse statement in the case where $X$ is affine or complete. In general, each implication in the above argument is reversible except for the first one, i.e. it is \textit{not} known to be true that if $\Cox(X)$ is a polynomial ring over $k$, then $X$ is a toric variety. However, this is known to be true when $X$ is complete (see \cite[Exercise 3.9]{arzhantsev2015cox}), so the converse statement holds in this case. Below we explain why this also works when $X$ is affine. 
		
		If $X$ is affine and $\Cox(X)$ is a polynomial ring over $k$, then $\Spec(\Cox(X))$ is $\A^n$ for some $n$, and we have a good quotient map $\A^n\to X$ for the action of a diagonalizable group $K$ (this is using the Cox construction for horospherical varieties; for instance, see \cite[Section 2]{gagliardi2019luna}). Moreover, there is a connected reductive algebraic group $G'$ such that $\A^n$ is a horospherical $G'$-variety, and $K$ acts by $G'$-equivariant automorphisms. As in the argument of \cite[Proposition 6.1]{gagliardi2017generalized}, $G'$ acts linearly on $\A^n$, and therefore $K$ acts linearly on $\A^n$. It follows that the good quotient $X$ of $\A^n$ by $K$ is a toric variety. 
	\end{proof}
\end{proposition}

\begin{example}\label{ex:non-toric affine SL_3/U_3-variety}
	Consider $G=\SL_3$. There is an action of $\SL_3$ on $\A^6=\A^3_{\vec x}\times\A^3_{\vec y}$ by matrix multiplication on the first $\A^3_{\vec x}$ factor and contragredient (i.e. inverse transpose) matrix multiplication on the second $\A^3_{\vec y}$ factor; here $\vec x$ and $\vec y$ represent coordinate vectors in $3$-space. That is, $A\in \SL_3$ acts by $A\cdot (\vec x,\vec y)=(A\vec x,(A^T)^{-1}\vec y)$. Notice that, if $\vec x\cdot \vec y=0$ (usual dot product), then $(A\vec x)\cdot ((A^T)^{-1}\vec y)=0$. It follows that, choosing the base point $p:=(1,0,0,0,0,1)$, the orbit
	\begin{align*}
		\SL_3\cdot p = \{\vec x\cdot \vec y=0 : \vec x,\vec y\neq 0\}
	\end{align*}
	is open in $X:=\{\vec x\cdot \vec y=0\}\subseteq \A^3_{\vec x}\times\A^3_{\vec y}$, and the stabilizer is $U_3$. Therefore, $X$ is a $\SL_3/U_3$-horospherical variety. Notice that $X$ is not toric because the associated flag variety $\SL_3/B_3$ is not a product of projective spaces (see \cref{prop:non-toric condition}). 
\end{example}

\begin{example}\label{ex:non-toric class of examples for SL_n}
	Consider $G=\SL_n$ for $n\geq 3$ and a horospherical subgroup $H\subseteq \SL_n$ whose associated parabolic $P$ is not maximal, i.e. $P=P_I$ corresponds to a subset $I\subseteq S$ with $\# I<n-2$ (e.g. $P=B_n$). Then $\SL_n/P$ is not a (product of) projective space(s), so \cref{prop:non-toric condition} tells us that every $\SL_n/H$-horospherical variety is not a toric variety. 
\end{example}

\begin{remark}\label{rmk:ss and torus isogeny horospherical varieties}
	As in \cref{rmk:ss and torus isogeny}, we have an isogeny $\varepsilon:G^{ss}\times T_0\to G$ where $G^{ss}$ is semisimple simply connected and $T_0$ is a torus. A horospherical homogeneous space $G/H$ is $\varepsilon$-equivariantly isomorphic to the horospherical homogeneous space $G^{ss}\times T_0/\varepsilon^{-1}(H)$. Moreover, these horospherical homogeneous spaces have the same equivariant embeddings. Therefore, $G/H$-horospherical varieties are the same as $G^{ss}\times T_0/\varepsilon^{-1}(H)$-horospherical varieties. 
\end{remark}

We end this section by defining the morphisms which preserve the horospherical structure.

\begin{definition}[Horospherical morphism]\label{def:horospherical morphism}
	Let $X_i$ be $G/H_i$-horospherical varieties for $i=1,2$. A \define{horospherical morphism} is a morphism of varieties $X_1\to X_2$ which is $G$-equivariant and sends the base point of $X_1$ to the base point of $X_2$. 
\end{definition}

From the definition of horospherical morphism, the map $X_1\to X_2$ restricts to a map of open orbits $G/H_1\to G/H_2$ and sends $eH_1\mapsto eH_2$. This means that we must have $H_1\subseteq H_2$ and the restricted map $G/H_1\to G/H_2$ is the natural projection. 

Therefore, all horospherical morphisms are dominant. This seems like a somewhat restrictive condition, but there are many important types of morphisms that are horospherical, such as the following: open embeddings of $G$-invariant subsets, $G$-equivariant isomorphisms, $G$-equivariant surjections, and blow-ups of $G$-invariant closed subsets.

\begin{remark}\label{rmk:comparing horospherical varieties different groups}
	Consider two horospherical homogeneous spaces $G_i/H_i$ where these $G_i$ are a connected reductive algebraic groups for $i=1,2$, and suppose that there is a surjective algebraic group homomorphism $\varepsilon:G_1\to G_2$. If $X_i$ are $G_i/H_i$-horospherical varieties for $i=1,2$, then we can view $X_2$ as a $G_1/\varepsilon^{-1}(H_2)$-horospherical variety where $G_1$ acts via $\varepsilon$. This perspective is useful because now we can view both $X_1$ and $X_2$ as horospherical $G_1$-varieties, and in particular we can study horospherical morphisms $X_1\to X_2$. 
\end{remark}

\section{Colours}\label{sec:colours}

Throughout this section, we fix a horospherical homogeneous space $G/H$, where $H\supseteq U$ is horospherical with associated parabolic $P=N_G(H)\supseteq B$. Write $P=P_I$ for a unique $I\subseteq S$, and let $\calC=\calC(G/H):=S\setminus I$. Recall that we have the lattice $N=N(G/H)$ associated to the torus $P/H=T_N$. 

In this section we put a ``colour structure" on $N$ by connecting $\calC$ and $N$ via the $B^-$-invariant divisors of $G/H$. We should think of $\calC$ as an abstract set of \define{colours} for $G/H$. Note that we prefer to consider $B^-$-invariant divisors instead of $B$-invariant divisors, which is not common in the literature. We choose to do this because the open $B$-orbit in $G/H$ does not contain the base point $eH$, but the open $B^-$-orbit does contain $eH$, so it is more convenient to work with the open $B^-$-orbit, in which case one should study the complementary $B^-$-invariant divisors. One can simply swap the roles of $B$ and $B^-$ (since $(B^-)^-=B$ and opposites are conjugate) throughout this paper to instead consider $B$-invariant divisors, if desired; in this case, $H$ should contain the unipotent radical $U^-$ of $B^-$ instead of that of $B$.

\subsection{Eigenfunctions}\label{subsec:eigenfunctions}

We start with some technical details regarding eigenfunctions and valuations. For any algebraic group $G'$ and any $G'$-module $V$, the set of $G'$-eigenvectors in $V$ is denoted $V^{(G')}:=\{v\in V\setminus\{0\}:g\cdot v=\chi(g)v \text{ for some } \chi\in\frakX(G'),~\forall g\in G'\}$. 

In particular, if $X$ is a normal $G$-variety, then the \define{$B^-$-eigenfunctions} of $k(X)$ are
\begin{align*}
	k(X)^{(B^-)} = \{f\in k(X)^* : b\cdot f=\chi(b)f \text{ for some } \chi\in\frakX(T),~\forall b\in B^-\}.
\end{align*}
Note that we used characters of the maximal torus $T$ above since $\frakX(B)=\frakX(B^-)=\frakX(T)$. The $B^-$-eigenfunctions are $\frakX(T)$-graded:
\begin{align*}
	k(X)^{(B^-)}\cup\{0\} = \bigoplus_{m\in \frakX(T)} k(X)^{(B^-)}_{-m}
\end{align*}
where
\begin{align*}
	k(X)^{(B^-)}_{-m} &:= \{f\in k(X) : b\cdot f=\chi^{-m}(b)f ~\forall b\in B^-\}
	\\&= \{f\in k(X) : f(b\cdot x)=\chi^m(b)f(x) ~\forall b\in B^-,x\in X\}
\end{align*}
where $\chi^m$ denotes the character function associated to the lattice point $m\in \frakX(T)$. Recall that $(b\cdot f)(x):=f(b^{-1}\cdot x)$, which explains the minus sign used above. 

Given $f\in k(X)^{(B^-)}$, let $m_f$ denote the element of $\frakX(T)$ such that $f\in k(X)^{(B^-)}_{-m_f}$. The assignment $f\mapsto m_f$ is a group homomorphism $k(X)^{(B^-)}\to\frakX(T)$; let $\Lambda(X)\subseteq\frakX(T)$ denote the image of this map. If $X$ is spherical, then the kernel of this map is exactly the constant rational functions, i.e. we have an exact sequence of groups
\begin{align}\label{eq:exact sequence of eigenfunctions}
	\begin{split}
		1 \map k^* \map k(X)^{(B^-)} &\map \Lambda(X) \map 0
		\\ f &\longmapsto m_f.
	\end{split}
\end{align}

\begin{lemma}\label{lemma:image of eigenfunctions is N^vee for horospherical}
	If $X$ is a $G/H$-horospherical variety, then $\Lambda(X)=N^\vee$. Therefore, we have $N^\vee\cong k(X)^{(B^-)}/k^*$. 
	
	\begin{proof}
		See the paragraph immediately before Section 7.3 in \cite{timashev2011homogeneous}; this is a corollary of \cite[Theorem 4.7]{timashev2011homogeneous}. 
	\end{proof}
\end{lemma}

For a normal $G$-variety $X$, let $\Val(X)$ denote the set of valuations of $X$, i.e. group homomorphisms $\nu:k(X)^*\to\R$ such that the image is a discrete subgroup of $\R$ (i.e. isomorphic to $\Z$), $\nu(k^*)=0$, and $\nu(f_1+f_2)\geq\min\{\nu(f_1),\nu(f_2)\}$ for all $f_1,f_2\in k(X)^*$. Let $\calV(X)\subseteq\Val(X)$ denote the set of $G$-invariant valuations. 

If $X$ is spherical, then \cref{eq:exact sequence of eigenfunctions} and the fact that valuations are trivial on constant functions imply that we have a set map
\begin{align}\label{eq:valuation map}
	\begin{split}
		\xi_{\Val}:\Val(X) &\map \Lambda(X)^\vee_\R = \Hom_\Z(\Lambda(X),\R)
		\\ \nu &\longmapsto (m_f\mapsto \nu(f)).
	\end{split}
\end{align}
For simplicity, we denote $u_\nu:=\xi_{\Val}(\nu)$ for $\nu\in\Val(X)$. If $D\subseteq X$ is a prime divisor then it has a corresponding order of vanishing valuation $\nu_D:k(X)^*\to\Z$; when $D$ is $G$-invariant, we have $\nu_D\in\calV(X)$. In this case, we write $u_D$ to denote $u_{\nu_D}$, which is a lattice point in $\Lambda(X)^\vee\subseteq\Lambda(X)^\vee_\R$. 

\begin{example}\label{ex:valuation points for toric varieties}
	Let $X$ be a toric variety with open torus $T$, i.e. $X$ is a $T/\{1\}$-horospherical variety. Then $X$ corresponds to a fan $\Sigma$ on the lattice $N=N(T/\{1\})$. Each ray $\rho$ of $\Sigma$ corresponds to a $T$-invariant prime divisor $D_\rho\subseteq X$. In the notation above, the point $u_{D_\rho}$ is exactly the minimal ray generator for $\rho$; see \cite[Proposition 4.1.1]{cox2011toric}. A generalization of this to all horospherical varieties is stated in \cref{cor:non-coloured rays and G-invariant divisors}. 
\end{example}

The map $\xi_{\Val}$ is usually not injective. However, for a spherical variety $X$, we see below that this map is injective when restricted to $\calV(X)$. Moreover, horospherical varieties can be characterized as the subclass of spherical varieties for which $\res{\xi_{\Val}}{\calV(X)}$ is surjective. 

\begin{proposition}[Valuation cone]\label{prop:valuation cone}
	Let $X$ be a normal $G$-variety. If $X$ is spherical, then the restricted map $\res{\xi_{\Val}}{\calV(X)}:\calV(X) \to \Lambda(X)^\vee_\R$ is injective, and the image is a convex polyhedral cone called the \define{valuation cone}. Moreover, this restricted map is a bijection if and only if $X$ is horospherical. 
	
	\begin{proof}
		See \cite[Corollary 1.8]{knop1991luna}, and see \cite[Corollary 6.2]{knop1991luna} for the ``moreover" claim. 
	\end{proof}
\end{proposition}

Let $X$ be a $G/H$-horospherical variety. Combining \cref{lemma:image of eigenfunctions is N^vee for horospherical} and \cref{prop:valuation cone}, we see that the valuation cone $\calV(X)$ is in bijection with $N_\R$, so points in $N_\R$ are all of the form $u_\nu$ for some unique $\nu\in\calV(X)$. However, some points in $N_\R$ may correspond to multiple valuations when we use \textit{non-$G$-invariant} valuations. Also, recall that each prime divisor $D\subseteq X$ yields a lattice point $u_D\in N$. 

\begin{remark}\label{rmk:eigenfunctions and valuations on open subset}
	If $X$ is a $G/H$-horospherical variety, then $k(X)^{(B^-)}=k(G/H)^{(B^-)}$, $\Lambda(X)=\Lambda(G/H)$, $\Val(X)=\Val(G/H)$, and $\calV(X)=\calV(G/H)$. Thus, everything developed in this subsection is the same for every $G/H$-horospherical variety because it only depends on $G/H$. 
\end{remark}

\subsection{Colour divisors}\label{subsec:colour divisors}

In a $G/H$-horospherical variety $X$, the open $B^-$-orbit $B^-H/H\subseteq X$ of the base point is affine (see \cite[Theorem 3.5]{timashev2011homogeneous} because $B^-$ is solvable), and the boundary $X\setminus (B^-H/H)$ is equal to the union of the $B^-$-invariant prime (i.e. irreducible) divisors of $X$. These $B^-$-invariant prime divisors come in two flavours: either they are $G$-invariant or they are not. The divisors of the latter flavour are key in developing the combinatorics of horospherical varieties. 

\begin{definition}[Colour divisor]\label{def:colour divisor}
	Let $X$ be a $G/H$-horospherical variety. A \define{colour divisor} of $X$ is a prime divisor which is $B^-$-invariant but not $G$-invariant.
\end{definition}

For a $G/H$-horospherical variety $X$, let $\calD(X)$ denote the (finite) set of $B^-$-invariant prime divisors of $X$. This breaks down as $\calD(X)=\calD^G(X)\sqcup\calD^\calC(X)$ where $\calD^G(X)$ and $\calD^\calC(X)$ are the sets of $G$-invariant prime divisors and colour divisors of $X$, respectively.

\begin{example}[Basic examples of colour divisors]\label{ex:basic examples of colour divisors}
	\mbox{}\par
	\begin{itemize}[leftmargin=2.5em]
		\item Consider the case where $G=T$ is a torus and $X$ is a toric variety with open torus $T$, i.e. $X$ is a $T/\{1\}$-horospherical variety. Then $B^-=T$, so there are no colour divisors for $X$. Thus, the set $\calD(X)$ is just $\calD^G(X)=\calD^T(X)$.
		\item Consider a flag variety $G/P$. Since $G/P$ is a single $G$-orbit, there are no $G$-invariant prime divisors. Thus, the set $\calD(G/P)$ is just $\calD^\calC(G/P)$. Note that the colour divisors are exactly the Schubert divisors. 
		\item If $X_i$ are horospherical $G_i$-varieties for $i=1,2$, then $\calD^\calC(X_1\times X_2)$ and $\calD^{G_1\times G_2}(X_1\times X_2)$ can be naturally identified with $\calD^\calC(X_1)\sqcup\calD^\calC(X_2)$ and $\calD^{G_1}(X_1)\sqcup\calD^{G_2}(X_2)$, respectively. \qedhere
	\end{itemize}
\end{example}

In the rest of this subsection we describe the $B^-$-invariant prime divisors of horospherical varieties. First, we compute the $B^-$-invariant prime divisors of $G/H$. Since $G/H$ is a single $G$-orbit, there are no $G$-invariant divisors, i.e. $\calD^G(G/H)=\varnothing$. So we just need to compute the colour divisors, which are the $B^-$-orbit closures of codimension $1$ in $G/H$. 

Recall that from \cref{eq:B-orbit decomposition of G/H}, the $B^-$-orbits of $G/H$ are the preimages of the $B^-$-orbits of $G/P$ under the natural projection $G/H\to G/P$. In particular, the colour divisors of $G/H$ are the pullbacks of the colour divisors of $G/P$. The colour divisors of $G/P$ are exactly the codimension $1$ Schubert varieties, i.e. the $\ol{B^-wP/P}\subseteq G/P$ where the length of $w\in W/W_I$ is $1$, i.e. $w=r_\alpha$ for some $\alpha\in \calC=S\setminus I$. Therefore, the colour divisors of $G/H$ are
\begin{align}\label{eq:colour divisors of G/H}
	\calD^\calC(G/H) = \{\ol{B^-r_\alpha H/H}\subseteq G/H : \alpha\in\calC\}.
\end{align}

The following proposition tells us about the $B^-$-invariant divisors of any $G/H$-horospherical variety. In particular, this proposition says that the colour divisors are determined by the colour divisors of $G/H$ given above. 

\begin{proposition}\label{prop:colour divisors of horospherical variety}
	Let $X$ be a $G/H$-horospherical variety. Every colour divisor of $X$ is the closure of a colour divisor of $G/H$, and conversely every colour divisor of $G/H$ is the restriction of a colour divisor of $X$. That is, we have a bijective correspondence
	\begin{align*}
		\calD^\calC(G/H) &\longleftrightarrow \calD^\calC(X)
		\\ D &\longmapsto \ol{D}
		\\ D\cap G/H &\longmapsfrom D.
	\end{align*}
	
	\begin{proof}
		The map $\calD^\calC(G/H)\to \calD^\calC(X)$ is clearly well defined. To show that the reverse map is well defined, we just need to check that $D\cap G/H\neq \varnothing$ for each $D\in\calD^\calC(X)$. If $D\cap G/H=\varnothing$, then $D$ is an irreducible component of the boundary $X\setminus G/H$, so $D$ must be a $G$-invariant prime divisor, which is a contradiction. 
	\end{proof}
\end{proposition} 

Let $X$ be a $G/H$-horospherical variety. As a consequence of \cref{prop:colour divisors of horospherical variety}, we have a bijection between abstract colours and colour divisors:\footnote{Most sources do not make this distinction, but we find it helpful to distinguish between the combinatorial abstract colours and the geometric colour divisors.}
\begin{align}\label{eq:bijection between colours and colour divisors}
	\calC \map \calD^\calC(X) \qquad \alpha\longmapsto D_\alpha^X:=\ol{B^-r_\alpha H/H}\subseteq X.
\end{align}
When there is no confusion, we simply write $D_\alpha$ rather than $D_\alpha^X$ for the colour divisor in $X$ corresponding to $\alpha\in\calC$. 

From \cref{ex:basic examples of colour divisors} we know that toric varieties (viewed as $T/\{1\}$-horospherical varieties) have no colour divisors, so \cref{eq:bijection between colours and colour divisors} implies that they are ``colourless". The following corollary tells us that toric varieties are the \textit{only} ``colourless" horospherical varieties; so we can think of colours (or colour \textit{divisors}) as the main difference between horospherical varieties and toric varieties. 

\begin{corollary}\label{cor:torus if and only if colourless}
	The horospherical homogeneous space $G/H$ is a torus if and only if $\calC=\varnothing$ (i.e. $G/H$ is ``colourless"). 
	
	\begin{proof}
		The forward direction follows from \cref{ex:basic examples of colour divisors} and \cref{eq:bijection between colours and colour divisors}. For the reverse direction, $\calC=\varnothing$ means $I=S$, so $G/H=P/H$ is a torus.
	\end{proof}
\end{corollary}

\subsection{Coloured lattices}\label{subsec:colour structure on N=N(G/H)}

Now we use the algebro-geometric data developed in the previous subsections to put a ``colour structure" on $N=N(G/H)$. Recall that $\calC=\calC(G/H)=S\setminus I$ is in bijection with $\calD^\calC(G/H)$; namely, the colour $\alpha\in\calC$ corresponds to the colour divisor $D_\alpha\subseteq G/H$. Let $\xi:\calC\to N$ be the map sending 
\begin{align}\label{eq:colour map}
\alpha\mapsto \xi(\alpha):=\xi_{\Val}(D_\alpha)=u_{D_\alpha}.
\end{align}
We also denote $\xi(\alpha)$ by $u_\alpha\in N$, or sometimes $u_\alpha^N\in N$ to signify that $u_\alpha$ is in $N$; these points are called \define{colour points}. These colour points can be difficult to compute directly from $\xi_{\Val}$, but there are other ways to determine them quickly; see \cref{ex:primary coloured lattice} and \cref{rmk:mapping from primary coloured lattice}. 

Note that the lattice $N$ depends on the torus $P/H$, the abstract set of colours $\calC$ only depends on the associated flag variety $G/P$, and the colour map $\xi$ depends on the divisors in $\calD(G/H)=\calD^\calC(G/H)$. 

\begin{definition}[Coloured lattice]\label{def:coloured lattice}
	The triple $(N,\calC,\xi)$ constructed above is called the \define{coloured lattice associated to $G/H$}. We call $\calC$ the \define{universal colour set} and $\xi$ the \define{colour map}. We usually just denote $(N,\calC,\xi)$ by $N=N(G/H)$, and the colour structure $\xi:\calC=\calC(G/H)\to N$ is implicit. 
\end{definition} 

\begin{remark}[FAQs about coloured lattices]\label{rmk:FAQs about coloured lattices}
	\mbox{}\par 
	\begin{itemize}[leftmargin=2.5em]
		\item Can a colour point $u_\alpha$ be at the origin in $N$? The answer is \textit{yes}, e.g. see the flag variety example in \cref{ex:basic examples of coloured lattices}, or \cref{ex:orbit closure coloured fan in projective SL_3/U_3 variety}. 
		\item Can two colours $\alpha_1,\alpha_2\in\calC$ yield the same colour point $u_{\alpha_1}=u_{\alpha_2}$ in $N$? The answer is \textit{yes}, e.g. consider flag varieties, or see \cref{ex:non-factorial coloured cone with double colour point for SL_3}. This is the reason why we keep track of the abstract colours rather than just the colour points.
		\item If $u_\alpha\neq 0$, then does $u_\alpha$ have to be a primitive lattice point in $N$? The answer is \textit{no}, e.g. see \cref{ex:coloured cone of A1 singularity}. 
		\item Does the coloured lattice uniquely determine the horospherical homogeneous space? In general, the answer is \textit{no} if one allows for different full groups $G$, but if one is only looking at a specific full group $G$ and keeps track of the data appropriately, then the answer \textit{can be yes}; this is made precise in \cref{thm:uniqueness of G/H}. \qedhere
	\end{itemize}
\end{remark}

\begin{example}[Basic examples of coloured lattices]\label{ex:basic examples of coloured lattices}
	\mbox{}\par 
	\begin{itemize}[leftmargin=2.5em]
		\item For a torus $T$, the coloured lattice $N=N(T/\{1\})$ is the usual one-parameter subgroup lattice for $T$, and there are no colour points. By \cref{cor:torus if and only if colourless}, a coloured lattice $N(G/H)$ has no colour points if and only if $G/H$ is a torus. 
		\item The coloured lattice $N=N(G/P)$ for $G/P$ is trivial, i.e. rank $0$, and all colour points are at the origin. It is easy to see that a coloured lattice $N(G/H)$ is trivial if and only if $H=P$ is parabolic. 
		\item If $G_i/H_i$ are two horospherical $G_i$-homogeneous spaces for $i=1,2$, then we have $N((G_1\times G_2)/(H_1\times H_2))=N(G_1/H_1)\times N(G_2/H_2)$ and the universal colour set is $\calC(G_1/H_1)\sqcup\calC(G_2/H_2)$. \qedhere
	\end{itemize}
\end{example}

\begin{example}[Colours of $\SL_n/U_n$]\label{ex:coloured lattice for SL_n/U_n}
	Consider the horospherical homogeneous space $G/H=\SL_n/U_n$. We describe $N=N(\SL_n/U_n)$ as follows. Recall the basis $\{e_1,\ldots,e_{n-1}\}$ for $\frakX(T_n)\cong\Z^{n-1}$ from \cref{ex:characters for SL_n maximal torus}. Since $N^\vee=\frakX(T_n)$, this gives us a basis for $N^\vee$, and in turn we get a dual basis for $N$ which we also denote $\{e_1,\ldots,e_{n-1}\}$. With this basis, we have $e_i=u_{\alpha_i}$ for each $1\leq i\leq n-1$, where $\calC=\{\alpha_1,\ldots,\alpha_{n-1}\}$ is the universal colour set for $N$. 
	
	Indeed, to see why this is true, we need to check that there is a basis of $B_n^-$-eigenfunctions $\{f_1,\ldots,f_{n-1}\}$ for $k(\SL_n/U_n)^{(B_n^-)}/k^*\cong N^\vee$ such that $m_{f_j}=e_j$ (in the notation of \cref{subsec:eigenfunctions}) and $\nu_{D_{\alpha_i}}(f_j)=\delta_{i,j}$ (Kronecker delta) for each $1\leq i,j\leq n-1$. Choose $f_j$ to be the rational function on $\SL_n/U_n$ which sends a matrix to the determinant of its upper-left $j\times j$-submatrix. Then clearly $m_{f_j}=e_j$. Moreover, each $f_j$ is regular and $D_{\alpha_i}=\{f_i=0\}$, so it follows that $\nu_{D_{\alpha_i}}(f_j)=\delta_{i,j}$. 
\end{example}

The previous example shows that the colour points of $N(\SL_n/U_n)$ form a basis for the lattice. The following example generalizes this to similar coloured lattices.

\begin{example}[Primary coloured lattice]\label{ex:primary coloured lattice}
	Recall from \cref{ex:primary horospherical homogeneous space} that $G/U$ is the primary horospherical homogeneous space (with respect to $G$). The associated coloured lattice $N(G/U)$ is called the \define{primary coloured lattice} (again, with respect to $G$). Note that the universal colour set is $\calC(G/U)=S$. In this case, the colour points $u_\alpha$ with $\alpha\in\calC(G/U)$ form part of a $\Z$-basis for $N(G/U)$; this is summarized in the introduction of \cite{gagliardi2014cox} since the divisor class group $\Cl(G/U)$ is trivial. 
	
	More generally, if $P\subseteq G$ is a standard parabolic subgroup, then $G/[P,P]$ is the $P$-primary horospherical homogeneous space. The associated coloured lattice $N(G/[P,P])$ is called the \define{$P$-primary coloured lattice}, and the colour points $u_\alpha$ with $\alpha\in\calC(G/[P,P])$ form part of a $\Z$-basis for $N(G/[P,P])$; again, this uses the fact that $\Cl(G/[P,P])$ is trivial. 
\end{example}

We end this subsection by stating a uniqueness result from \cite{losev2009uniqueness} for horospherical homogeneous spaces in terms of their coloured lattices. Let $G_{D_\alpha}\subseteq G$ denote the stabilizer of a colour divisor $D_\alpha$; note that this is the maximal parabolic subgroup $P_{S\setminus\{\alpha\}}^-\supseteq B^-$. 

\begin{theorem}[Uniqueness of $G/H$]\label{thm:uniqueness of G/H}
	Let $H_1,H_2\subseteq G$ be any horospherical subgroups. Then $G/H_1$ and $G/H_2$ are $G$-equivariantly isomorphic if and only if the following holds: there is an isomorphism of lattices $\theta:N(G/H_1)\to N(G/H_2)$, and there exists a bijection $\iota:\calC(G/H_1)\to \calC(G/H_2)$ such that $G_{D_\alpha}=G_{D_{\iota(\alpha)}}$ and $\theta(u_\alpha)=u_{\iota(\alpha)}$ for each $\alpha\in\calC(G/H_1)$. 
	
	\begin{proof}
		See \cite[Theorem 1]{losev2009uniqueness}. 
	\end{proof}
\end{theorem}

\subsection{Maps of coloured lattices}\label{subsec:maps of coloured lattices}

Now let $G/H_i$ be two horospherical homogeneous spaces for $i=1,2$ and assume that $H_1\subseteq H_2$. We use the same notation developed previously for $G/H$ but now with subscripts $i=1,2$; e.g. $P_i=N_G(H_i)$ and $N_i=N(G/H_i)$. 

Consider the natural projection $\pi:G/H_1\to G/H_2$. For each colour divisor $D$ of $G/H_1$ there are two cases: either (1) $\ol{\pi(D)}$ is a colour divisor of $G/H_2$, or (2) $\ol{\pi(D)}=G/H_2$ (i.e. $D$ maps dominantly via $\pi$). We relate these possibilities to the combinatorics as follows. 

Since $P_i=TH_i$ for $i=1,2$, we have $P_1=P_{I_1}\subseteq P_2=P_{I_2}$, so $I_1\subseteq I_2$, and thus $\calC_2=S\setminus I_2\subseteq\calC_1=S\setminus I_1$. Denote the colour divisor in $G/H_i$ corresponding to $\alpha\in\calC_i$ by $D_\alpha^i$. The inclusion $\calC_2\subseteq\calC_1$ is viewed in terms of colour divisors as follows: each $D_\alpha^2\in\calD^\calC(G/H_2)$ pulls back to $D_\alpha^1=\pi^{-1}(D_\alpha^2)\in\calD^\calC(G/H_1)$. Thus, the $D_\alpha^1$ with $\alpha\in\calC_2$ are precisely the colour divisors of $G/H_1$ which map to colour divisors of $G/H_2$, i.e. these are the colour divisors in case (1) listed in the previous paragraph. On the other hand, the $D_\alpha^1$ with $\alpha\in\calC_1\setminus\calC_2$ are the colour divisors of $G/H_1$ which map dominantly to $G/H_2$ via $\pi$, i.e. these are the colour divisors in case (2) listed above.

For $i=1,2$, we have natural projection maps $G/H_i\to G/P_i$. Therefore, we have a commutative diagram of horospherical morphisms (natural projection maps):
\begin{equation*}
	\begin{tikzcd}
		{G/H_1} && {G/H_2} \\
		{G/P_1} && {G/P_2}
		\arrow["\pi", from=1-1, to=1-3]
		\arrow[from=1-3, to=2-3]
		\arrow[from=1-1, to=2-1]
		\arrow[from=2-1, to=2-3]
	\end{tikzcd}
\end{equation*}
The vertical maps are principal $P_i/H_i$-bundles, and we get an induced map of the fibres $T_{N_1}=P_1/H_1\to T_{N_2}=P_2/H_2$ which is a surjective homomorphism of tori. This induces an inclusion of character lattices and a dual map of one-parameter subgroup lattices:
\begin{align}\label{eq:associated map of coloured lattices}
	\Phi^\vee:N_2^\vee=\frakX(T_{N_2})\hookrightarrow N_1^\vee=\frakX(T_{N_1}), \qquad \Phi:N_1=\frakX^*(T_{N_1})\to N_2=\frakX^*(T_{N_2}).
\end{align}
Note that $\Phi_\R$ is surjective since $\Phi^\vee$ is injective. 

We call the map $\Phi:N_1\to N_2$ the \define{associated map of coloured lattices} (associated to $\pi:G/H_1\to G/H_2$). This map preserves the colour structure of these lattices in the following way. Recall that $\calC_2\subseteq\calC_1$. Denote $\calC_\Phi:=\calC_1\setminus\calC_2=I_2\setminus I_1$, which we call the set of \define{dominantly mapped colours} (the usage of \textit{dominant} is explained above). Note that $\calC_1=\calC_\Phi\sqcup\calC_2$. If $u_\alpha^i$ denotes the colour point in $N_i$ corresponding to $\alpha$, then it is easy to check that
\begin{align}\label{eq:colour points map to colour points}
	\Phi(u_\alpha^1)=u_\alpha^2 ~~\forall\alpha\in\calC_2 \quad\text{and}\quad \Phi(u_\alpha^1)=0 ~~\forall \alpha\in\calC_\Phi. 
\end{align}

\begin{remark}\label{rmk:alternative description of map of coloured lattices}
	Following \cite[Section 4]{knop1991luna}, there is an equivalent way of thinking about the associated map of coloured lattices $\Phi:N_1\to N_2$ constructed above. Since $\pi:G/H_1\to G/H_2$ is dominant, we get an induced inclusion $k(G/H_2)\hookrightarrow k(G/H_1)$ of function fields. Since $\pi$ is $G$-equivariant, it is easy to check that this induces an inclusion $k(G/H_2)^{(B^-)}\hookrightarrow k(G/H_1)^{(B^-)}$. Therefore by \cref{lemma:image of eigenfunctions is N^vee for horospherical}, we have an inclusion $N_2^\vee\hookrightarrow N_1^\vee$; this is the same as the inclusion in \cref{eq:associated map of coloured lattices}. Then the dual map is exactly $\Phi:N_1\to N_2$ given above. 
\end{remark}

\begin{remark}\label{rmk:mapping from primary coloured lattice}
	Recall from \cref{ex:primary coloured lattice} that we have the primary coloured lattice $N(G/U)$. For any other coloured lattice $N=N(G/H)$, we have an associated map of coloured lattices $N(G/U)\to N$ coming from the projection $G/U\to G/H$. On the other hand, if $H$ has associated parabolic $P$, then we have a map of coloured lattices $N(G/[P,P])\to N$ associated to the projection $G/[P,P]\to G/H$. Combining these two maps we get a commutative diagram of coloured lattice maps:
	\begin{equation*}
		\begin{tikzcd}
			{N(G/U)} & {N(G/[P,P])} & {N}
			\arrow[from=1-1, to=1-2]
			\arrow[from=1-2, to=1-3]
			\arrow[shift left=1, curve={height=-12pt}, from=1-1, to=1-3]
		\end{tikzcd}
	\end{equation*}
	The dual maps are the inclusions of character lattices $N^\vee\subseteq\frakX(P)\subseteq\frakX(T)$. By \cref{eq:colour points map to colour points}, the colour points of $N$ are exactly the images of the non-dominantly mapped colour points in $N(G/U)$ or $N(G/[P,P])$. 
\end{remark}

\begin{example}[Quotient coloured lattice]\label{ex:quotient coloured lattice}
	Let $N'\subseteq N=N(G/H)$ be a saturated sublattice and let $\calC'\subseteq\calC=\calC(G/H)$ be a subset of colours such that $u_\alpha\in N'$ for each $\alpha\in\calC'$; we call $N'$ a \define{coloured sublattice} of $N$. Then the quotient lattice $N/N'$ is the coloured lattice associated to some horospherical homogeneous space $G/H'$ with $H'\supseteq H$; see \cite[Theorem 4.4]{knop1991luna}. The colour structure of $N/N'$ is described as follows. The quotient map $N\to N/N'$ is a map of coloured lattices (associated to the projection $G/H\to G/H'$) whose set of dominantly mapped colours is $\calC'$. So the universal colour set of $N/N'$ is $\calC\setminus\calC'$, and the colour points of $N/N'$ are the images of the colour points of $N$ from $\calC\setminus\calC'$. 
	
	We can determine $H'$ as follows. Since $\calC=S\setminus I$ and $\calC'\subseteq\calC$, we can write $\calC'=S\setminus I'$ for a unique $I'\subseteq S$. The dual map $(N/N')^\vee\hookrightarrow N^\vee$ is an inclusion of character lattices $M':=\frakX(T_{N/N'})\hookrightarrow \frakX(T_N)$. Then $H'=H_{(I',M')}$.
\end{example}

\section{Classification of horospherical varieties}\label{sec:classification of horospherical varieties}

Now we describe the combinatorial classification of horospherical varieties via coloured fans. This section is self-contained regarding the polyhedral geometry, but we recommend \cite{cox2011toric} (mainly Section 1.2 and Section 3.1) for more on the relevant definitions and properties. 

Throughout this section, we fix a horospherical homogeneous space $G/H$ and use the notation from the previous sections. In particular, $P=P_I=N_G(H)\supseteq B$ is the associated parabolic, $N=N(G/H)$ is the associated coloured lattice with universal colour set $\calC=\calC(G/H)=S\setminus I$, and $\xi:\calC\to N$ is the colour map.

\subsection{Simple open cover}\label{subsec:simple open cover}

In toric geometry, a general toric variety admits an open cover by affine toric varieties. However, horospherical $G$-varieties usually do not admit an open cover by affine horospherical $G$-varieties. Shifting perspective, one can think of affine toric varieties as being the toric varieties with a unique closed torus orbit. This perspective proves useful for generalizing the ideas in toric geometry to horospherical geometry. In this subsection, we show that every horospherical $G$-variety admits an open cover by horospherical $G$-varieties with a unique closed $G$-orbit. 

\begin{definition}[Simple horospherical variety]\label{def:simple horospherical variety}
	A $G/H$-horospherical variety is called \define{simple} if it contains a unique closed $G$-orbit. 
\end{definition}

\begin{example}[Basic examples of simple horospherical varieties]\label{ex:basic examples of simple horospherical varieties}
	\mbox{}\par 
	\begin{itemize}[leftmargin=2.5em]
		\item The horospherical homogeneous space $G/H$ is simple because it has a single $G$-orbit. In particular, flag varieties are simple. 
		\item Any affine toric variety is simple because it has a unique closed torus orbit (this follows from the orbit-cone correspondence \cite[Theorem 3.2.6]{cox2011toric}). 
		\item If $X_i$ are two simple horospherical $G_i$-varieties for $i=1,2$, then $X_1\times X_2$ is a simple horospherical $G_1\times G_2$-variety. \qedhere
	\end{itemize}
\end{example}

Let $X$ be a $G/H$-horospherical variety and let $\calO$ be any $G$-orbit. Set 
\begin{align*}
	& \calD_\calO(X):=\{D\in\calD(X):\calO\subseteq D\},
	\\& \calD^G_\calO(X):=\{D\in\calD^G(X):\calO\subseteq D\},
	\\& \calD^\calC_\calO(X):=\{D\in\calD^\calC(X):\calO\subseteq D\}.
\end{align*}
Note that, if $X$ is simple and $\calO$ is the unique closed $G$-orbit, then $\calD^G_\calO(X)=\calD^G(X)$ since every $G$-invariant prime divisor contains $\calO$.

For any $G$-orbit $\calO\subseteq X$, let 
\begin{align*}
	X_{\calO,B^-} := X\setminus \bigcup (\calD(X)\setminus \calD_\calO(X))
\end{align*} 
(so the union is over all $B^-$-invariant prime divisors of $X$ which do \textit{not} contain $\calO$). Set $X_{\calO,G}:=G\cdot X_{\calO,B^-}$. Clearly $X_{\calO,B^-}$ is open in $X$, so $X_{\calO,G}$ is also open in $X$. 

\begin{proposition}\label{prop:describing the X_{O,B}s and X_{O,G}s}
	Let $X$ be a $G/H$-horospherical variety and let $\calO$ be any $G$-orbit. 
	\begin{enumerate}
		\item $X_{\calO,B^-}$ is a $B^-$-invariant affine open subset of $X$ and $\calO\cap X_{\calO,B^-}$ is a $B^-$-orbit. 
		\item $X_{\calO,B^-}=\{x\in X:\ol{B^-\cdot x}\supseteq\calO\}$ and $X_{\calO,G}=\{x\in X:\ol{G\cdot x}\supseteq\calO\}$. 
		\item $\calO$ is the unique closed $G$-orbit in $X_{\calO,G}$. In particular, $X_{\calO,G}$ is a simple $G/H$-horospherical variety. Moreover, if $X$ is simple and $\calO$ is the unique closed $G$-orbit, then $X=X_{\calO,G}$. 
	\end{enumerate}
	
	\begin{proof}
		See \cite[Theorem 2.1]{knop1991luna} for (1) and (3). Part (2) uses $\charac(k)=0$; see \cite[Remark 5.4]{gandini2018embeddings-spherical}. 
	\end{proof}
\end{proposition}

\begin{proposition}\label{prop:open cover by simples}
	Let $X$ be a $G/H$-horospherical variety. Then $X$ has a finite open cover consisting of simple $G/H$-horospherical varieties. Specifically, $X$ is covered by the open subvarieties $X_{\calO,G}$ as $\calO$ ranges over the $G$-orbits of $X$.
	
	\begin{proof}
		This follows from \cref{prop:describing the X_{O,B}s and X_{O,G}s} and the fact that $X$ has finitely many $G$-orbits. 
	\end{proof}
\end{proposition}

\subsection{Coloured cone of a simple horospherical variety}\label{subsec:coloured cone of a simple horospherical variety}

In the correspondence between toric varieties and fans, the simple toric varieties correspond to cones. We illustrate how this works for simple horospherical varieties using \textit{coloured} cones, where the colour structure comes from \cref{sec:colours}.

For us, a \define{cone} on $N$ is a rational polyhedral cone in $N_\R$, i.e. a set of the form $\Cone(A):=\R_{\geq 0}A$ for some finite set $A\subseteq N$. Recall that a cone is said to be \define{strongly convex} if it does not contain any positive-dimensional subspaces of $N_\R$ (see \cite[Proposition 1.2.12]{cox2011toric} for equivalent characterizations). 

\begin{definition}[Coloured cone]\label{def:coloured cone}
	A \define{coloured cone} on the coloured lattice $N$ is a pair $\sigma^c=(\sigma,\calF)$ where	$\sigma\subseteq N_\R$ is a cone and $\calF\subseteq\calC$ satisfies $\xi(\calF)\subseteq\sigma$ and $0\notin\xi(\calF)$. We say that $\sigma^c$ is \define{strongly convex} if $\sigma$ is strongly convex.
\end{definition}

Given a coloured cone $\sigma^c$ on $N$, we let $\sigma$ denote the underlying cone, and we let $\calF(\sigma^c)\subseteq\calC$ denote the \define{colour set}; so we have $\sigma^c=(\sigma,\calF(\sigma^c))$.  

For a $G/H$-horospherical variety $X$ and any $G$-orbit $\calO$, let 
\begin{align*}
	\calF_\calO(X):= \{\alpha\in\calC:D_\alpha\in\calD^\calC_\calO(X)\} 
\end{align*}
(so $\calF_\calO(X)$ consists of the $\alpha\in\calC$ such that $D_\alpha$ contains $\calO$). Note that $\calF_\calO(X)$ is in bijection with $\calD^\calC_\calO(X)$ by \cref{eq:bijection between colours and colour divisors}. This notation just allows us to explicitly view $\calD^\calC_\calO(X)$ as a subset of $\calC$. We also let 
\begin{align*}
	\calF(X):= \bigcup_{\calO\subseteq X} \calF_\calO(X)
\end{align*}
where the union is over all $G$-orbits $\calO\subseteq X$. In words, $\calF(X)$ is the set of $\alpha\in\calC$ such that $D_\alpha$ contains a $G$-orbit of $X$. Note that, when $X$ is simple and $\calO$ is the unique closed $G$-orbit, we have $\calF(X)=\calF_\calO(X)$. 

\begin{definition}[Coloured cone of a simple horospherical variety]\label{def:coloured cone of a simple horospherical variety}
	Let $X$ be a simple $G/H$-horospherical variety with unique closed $G$-orbit $\calO$. We define the \define{coloured cone associated to $X$} as $\sigma^c(X):=(\sigma(X),\calF(X))$	where $\sigma(X)$ is the cone in $N_\R$ given by
	\begin{align*}
		\sigma(X):=\Cone(u_D:D\in\calD_\calO(X))=\Cone(u_\alpha,u_D:\alpha\in\calF(X),D\in\calD^G(X)).
	\end{align*}
	This $\sigma^c(X)$ is a strongly convex coloured cone on $N$ (see \cite[Theorem 3.1]{knop1991luna}). 
\end{definition}

\begin{example}[Basic examples of associated coloured cones]\label{ex:basic examples of associated coloured cones}
	\mbox{}\par 
	\begin{itemize}[leftmargin=2.5em]
		\item The horospherical homogeneous space $G/H$ is a simple $G/H$-horospherical variety, and its associated coloured cone is the \define{trivial coloured cone} $0^c=(0,\varnothing)$ on $N$; note that this is the unique coloured cone whose underlying cone is trivial. In particular, $0^c$ is the coloured cone associated to any flag variety.
		\item If $X$ is an affine toric variety, then it corresponds to a strongly convex cone $\sigma$ by \cite[Theorem 1.3.5]{cox2011toric}. Viewing $X$ as a $T/\{1\}$-horospherical variety, the coloured cone associated to $X$ is $(\sigma,\varnothing)$. 
		\item If $X_i$ are simple $G_i/H_i$-horospherical varieties for $i=1,2$, then $\sigma^c(X_1\times X_2)=\sigma^c(X_1)\times\sigma^c(X_2)$, where the product of these coloured cones is defined componentwise as $(\sigma(X_1)\times\sigma(X_2),\calF(X_1)\sqcup\calF(X_2))$. \qedhere
	\end{itemize}
\end{example}

When drawing a coloured cone $\sigma^c=(\sigma,\calF)$, we draw the cone $\sigma$, we use ``un-filled" circular points to denote the colour points in the coloured lattice (the point $u_\alpha$ is marked with $\alpha$), and we ``fill in" (with the colour orange) the colour point $u_\alpha$ to indicate that $\alpha\in\calF$. 

\begin{example}\label{ex:coloured cone for A^2 and blow up}
	Consider $G=\SL_2$ and $H=U_2$. Then the coloured lattice $N$ is isomorphic to $\Z$, we have $\calC=\{\alpha\}$, and $u_\alpha=e_1$. 
	
	Recall from \cref{ex:A^2 and P^2 are SL_2/U_2-varieties} that $\A^2$ is a $\SL_2/U_2$-horospherical variety where $\SL_2$ acts by ordinary matrix multiplication. The unique closed $\SL_2$-orbit is the origin $\{0\}\subseteq \A^2$, so $\A^2$ is simple. The colour divisor $D_\alpha$ is $\{(0,y)\}\subseteq \A^2$, which clearly contains the closed $\SL_2$-orbit, so $\calF(\A^2)=\{\alpha\}$. One can easily check that $\A^2$ has no $\SL_2$-invariant divisors. Therefore $\sigma(\A^2)=\Cone(u_\alpha)=\Cone(e_1)$, so $\sigma^c(\A^2)=(\Cone(e_1),\{\alpha\})$; see the coloured cone on the left in the diagram below.  
	
	Now consider the blow-up $\Bl_0\A^2$ of $\A^2$ at the origin. This has an $\SL_2$-action which lifts the action on $\A^2$, and by choosing a base point which maps to the base point of $\A^2$ we can view $\Bl_0\A^2$ as a $\SL_2/U_2$-horospherical variety. Since $0\in\A^2$ is fixed by $\SL_2$, the exceptional divisor $E$ is an $\SL_2$-invariant prime divisor. One can check that the corresponding lattice point is $u_E=e_1\in N=\Z$. It is easy to see that $E$ is the unique closed $\SL_2$-orbit of $\Bl_0\A^2$, so $\Bl_0\A^2$ is simple. Now the colour divisor $D_\alpha\subseteq\Bl_0\A^2$ does \textit{not} contain $E$, so $\calF(\Bl_0\A^2)=\varnothing$. Therefore $\sigma(\Bl_0\A^2)=\Cone(u_E)=\Cone(e_1)$, so $\sigma^c(\Bl_0\A^2)=(\Cone(e_1),\varnothing)$; see the coloured cone on the right in the diagram below. 
	
	\tikzitfig{coloured-cones-of-A2-and-blow-up}
\end{example}

\begin{example}\label{ex:coloured cone of A1 singularity}
	Consider $G=\SL_2$ and $H=H_{(\varnothing,2\Z)}$ from \cref{ex:horospherical subgroups of SL_2}. Then the coloured lattice $N$ is isomorphic to $\Z$, we have $\calC=\{\alpha\}$, and $u_\alpha=2e_1$. Note that $N^\vee$ is regarded as the sublattice $2\Z$ of $\frakX(T_2)=\Z$. 
	
	Let $X$ be the affine variety $\{xz=y^2\}\subseteq\A^3$, where $x,y,z$ are the $\A^3$-coordinates. So $X$ is the affine cone over the projective variety $C\subseteq\P^2$ given by the image of $\P^1\hookrightarrow\P^2$ via $(u:v)\mapsto(u^2:uv:v^2)$. There is a natural action of $\SL_2$ on $C$: we treat $(u:v)\in\P^1$ as a vector with two coordinates and $\SL_2$ acts by matrix multiplication. Then $\SL_2$ acts on $X$ by lifting the action on $C$. The point $(1,0,0)\in X$ has stabilizer $H$ and the orbit $\SL_2/H=X\setminus\{0\}$ is open in $X$, so $X$ is a $\SL_2/H$-horospherical variety. 
	
	The unique closed $\SL_2$-orbit of $X$ is the origin $\{0\}\subseteq X$, so $X$ is simple. The colour divisor $D_\alpha\subseteq X$ is $\{x=0\}\subseteq X$. This clearly contains the closed $\SL_2$-orbit, so $\calF(X)=\{\alpha\}$. One can check that $X$ has no $\SL_2$-invariant divisors. Therefore $\sigma(X)=\Cone(u_\alpha)=\Cone(2e_1)=\Cone(e_1)$, so $\sigma^c(X)=(\Cone(e_1),\{\alpha\})$; see the diagram below.
	
	\tikzitfig{coloured-cone-of-A2-mod-Z-2}
\end{example}

\begin{example}\label{ex:SL_3/U_3 affine variety associated coloured cone}
	Consider $G=\SL_3$ and $H=U_3$. Then the coloured lattice $N$ is isomorphic to $\Z^2$, we have $\calC=\{\alpha_1,\alpha_2\}$, and $u_{\alpha_i}=e_i$ are the standard basis vectors for $i=1,2$. 
	
	Let $X$ be the affine variety $\{ux+vy+wz=0\}\subseteq\A^6$ from \cref{ex:non-toric affine SL_3/U_3-variety}; recall that $X$ is a $\SL_3/U_3$-horospherical variety. The unique closed $\SL_3$-orbit of $X$ is the origin $\{0\}\subseteq X$, which is contained in both colour divisors $D_{\alpha_1}$ and $D_{\alpha_2}$, so $\calF(X)=\{\alpha_1,\alpha_2\}$. One can check that $X$ has no $\SL_3$-invariant divisors. Therefore $\sigma(X)=\Cone(u_{\alpha_1},u_{\alpha_2})=\Cone(e_1,e_2)$, so $\sigma^c(X)=(\Cone(e_1,e_2),\{\alpha_1,\alpha_2\})$; see the diagram below. 
	
	\tikzitfig{coloured-cone-of-affine-SL-3-U-3-variety} 
\end{example}

The following is the main theorem of this subsection, which says that simple horospherical varieties correspond to strongly convex coloured cones. This generalizes the correspondence between affine toric varieties and strongly convex cones (see \cite[Theorem 1.3.5]{cox2011toric}). 

\begin{theorem}[Classification of simple horospherical varieties]\label{thm:classification of simple horospherical varieties}
	The following map is a bijection:
	\begin{align*}
		\{\text{Simple horospherical } G/H \text{-varieties}\}/\cong &\map \{\text{Strongly convex coloured cones on } N\}
		\\ X &\longmapsto \sigma^c(X).
	\end{align*}
	(Here, $\cong$ refers to $G$-equivariant isomorphisms).
	
	\begin{proof}
		See \cite[Theorem 3.1]{knop1991luna}. 
	\end{proof}
\end{theorem}

\begin{notation}\label{not:simple horospherical variety from coloured cone}
	Given a strongly convex coloured cone $\sigma^c$ on $N$, we let $X_{\sigma^c}$ denote the simple $G/H$-horospherical variety corresponding to $\sigma^c$ via \cref{thm:classification of simple horospherical varieties} (this is defined up to $G$-equivariant isomorphism). 
\end{notation}

\subsection{Coloured fan of a horospherical variety}\label{subsec:coloured fan of a horospherical variety}

The classification of general toric varieties uses fans, which are collections of cones that glue together along boundaries. The way that the cones piece together in the fan is reflected geometrically by the simple open cover of the corresponding toric variety. We illustrate how this works for horospherical varieties: we define \textit{coloured} fans, which are collections of \textit{coloured} cones, and we see how the data of a coloured fan encodes the simple open cover of a horospherical variety developed in \cref{subsec:simple open cover}. 

Recall that a \define{face} of a cone $\sigma$ is a ``boundary cone" of $\sigma$, i.e. a subset of the form $\sigma\cap m^\perp$ for some $m\in\sigma^\vee$, where $\sigma^\vee:=\{m\in N^\vee_\R:\langle m,u\rangle\geq 0 ~\forall u\in\sigma\}$ is the \define{dual cone} of $\sigma$ and $m^\perp\subseteq N^\vee_\R$ denotes the hyperplane orthogonal to $m$.

\begin{definition}[Coloured face]\label{def:coloured face}
	Let $\sigma^c$ be a coloured cone on $N$. A \define{coloured face} of $\sigma^c$ is a coloured cone $\tau^c$ on $N$ whose underlying cone $\tau\subseteq\sigma$ is a face of $\sigma$, and whose colour set is given by $\calF(\tau^c)=\{\alpha\in\calF(\sigma^c):u_\alpha\in\tau\}$. 
\end{definition}

Note that $\sigma^c$ is always a coloured face of itself. Also, $\sigma^c$ is strongly convex if and only if the trivial coloured cone $0^c$ is a coloured face of $\sigma^c$. For each $l\in\Z_{\geq 0}$, let $\sigma^c(l)$ denote the set of $l$-dimensional coloured faces of $\sigma^c$; in particular, $\sigma^c(1)$ is the set of coloured rays. Note that each coloured ray $\rho^c\in\sigma^c(1)$ has a minimal ray generator $u_\rho\in\rho\cap N$ (i.e. this is a primitive element in $N$), and $\sigma$ is generated by $\{u_\rho:\rho^c\in\sigma^c(1)\}$. 

\begin{lemma}\label{lemma:G-orbits and coloured faces}
	Let $X$ be a simple $G/H$-horospherical variety. For each $G$-orbit $\calO$, the coloured cone $\sigma^c(X_{\calO,G})$ is a coloured face of $\sigma^c(X)$. Moreover, every coloured face of $\sigma^c(X)$ is of this form. 
	
	\begin{proof}
		See \cite[Lemma 3.2]{knop1991luna}. 
	\end{proof}
\end{lemma}

Before defining coloured fans, we note that containment and intersections of coloured cones are defined componentwise, i.e. $(\sigma_1,\calF_1)\subseteq (\sigma_2,\calF_2)$ if and only if $\sigma_1\subseteq\sigma_2$ and $\calF_1\subseteq\calF_2$, and $(\sigma_1,\calF_1)\cap(\sigma_2,\calF_2)=(\sigma_1\cap\sigma_2,\calF_1\cap\calF_2)$.

\begin{definition}[Coloured fan]\label{def:coloured fan}
	A \define{coloured fan} $\Sigma^c$ on the coloured lattice $N$ is a finite collection of strongly convex coloured cones on $N$ which satisfies the following:
	\begin{enumerate}
		\item Every coloured face of any coloured cone in $\Sigma^c$ is also in $\Sigma^c$.
		\item The intersection of any two coloured cones in $\Sigma^c$ is a coloured face of both.
	\end{enumerate}
\end{definition}

To obtain the definition of an ordinary \define{fan} in polyhedral geometry, simply remove the word ``coloured" from \cref{def:coloured fan}. 

Given a coloured fan $\Sigma^c$ on $N$, we let $\Sigma:=\{\sigma:\sigma^c\in\Sigma^c\}$ denote the underlying fan, and we let $\calF(\Sigma^c):=\cup_{\sigma^c\in\Sigma^c} \calF(\sigma^c)$ denote the \define{colour set}. For each $l\in\Z_{\geq 0}$, let $\Sigma^c(l)$ denote the set of $l$-dimensional coloured cones in $\Sigma^c$; in particular, $\Sigma^c(1)$ is the set of coloured rays. We also let $\Sigma^c_{\max}$ denote the set of maximal coloured cones in $\Sigma^c$.

\begin{remark}\label{rmk:coloured cone as coloured fan}
	A coloured cone $\sigma^c$ on $N$ can be viewed as a coloured fan on $N$ by taking the (finite) collection of all coloured faces of $\sigma^c$. When we say that a coloured fan is ``generated by a single coloured cone", we mean that it has a unique maximal coloured cone in this sense. 
\end{remark}

Recall that, in a $G/H$-horospherical variety $X$, any $G$-orbit $\calO$ determines an open simple $G/H$-horospherical variety $X_{\calO,G}\subseteq X$ whose unique closed $G$-orbit is $\calO$. These $X_{\calO,G}$'s form an open cover for $X$. By \cref{thm:classification of simple horospherical varieties}, $X_{\calO,G}$ corresponds to a strongly convex coloured cone $\sigma^c(X_{\calO,G})$ on $N$.

\begin{definition}[Coloured fan of a horospherical variety]\label{def:coloured fan of a horospherical variety}
	Let $X$ be a $G/H$-horospherical variety. We define the \define{coloured fan associated to $X$} as 
	\begin{align*}
		\Sigma^c(X):= \{\sigma^c(X_{\calO,G}) : \calO\subseteq X \text{ is a $G$-orbit}\}.
	\end{align*}
	This $\Sigma^c(X)$ is a coloured fan on $N$ (see \cite[Theorem 3.3]{knop1991luna}). 
\end{definition}

Note that, for a $G/H$-horospherical variety $X$, we have $\calF(X)=\calF(\Sigma^c(X))$. 

\begin{example}[Basic examples of associated coloured fans]\label{ex:basic examples of associated coloured fans}
	\mbox{}\par
	\begin{itemize}[leftmargin=2.5em]
		\item If $X$ is a simple $G/H$-horospherical variety, then $\Sigma^c(X)$ is generated by the single coloured cone $\sigma^c(X)$ from \cref{def:coloured cone of a simple horospherical variety}; here we are using the terminology of \cref{rmk:coloured cone as coloured fan}. 
		\item If $X$ is a toric variety, then it corresponds to a fan $\Sigma$ by \cite[Corollary 3.1.8]{cox2011toric}. Viewing $X$ as a $T/\{1\}$-horospherical variety, the coloured fan associated to $X$ is $\{(\sigma,\varnothing):\sigma\in\Sigma\}$. 
		\item If $X_i$ are $G_i/H_i$-horospherical varieties for $i=1,2$, then $\Sigma^c(X_1\times X_2)=\Sigma^c(X_1)\times\Sigma^c(X_2)$, where this product coloured fan is the collection of all product coloured cones $\sigma_1^c\times\sigma_2^c$ with $\sigma_i^c\in\Sigma^c(X_i)$ (see \cref{ex:basic examples of associated coloured cones}). \qedhere
	\end{itemize}
\end{example}

The following result is the main theorem of this section, which says that horospherical varieties correspond to coloured fans. This generalizes the correspondence between toric varieties and fans (see \cite[Theorem 3.1.5, Corollary 3.1.8]{cox2011toric}). 

\begin{theorem}[Classification of horospherical varieties]\label{thm:classification of horospherical varieties}
	The following map is a bijection:
	\begin{align*}
		\{\text{$G/H$-horospherical varieties}\}/\cong &\map \{\text{Coloured fans on } N\}
		\\ X &\longmapsto \Sigma^c(X).
	\end{align*}
	(Here, $\cong$ refers to $G$-equivariant isomorphisms).
	
	\begin{proof}
		See \cite[Theorem 3.3]{knop1991luna}. 
	\end{proof}
\end{theorem}

\begin{notation}\label{not:horospherical variety from coloured fan}
	Given a coloured fan $\Sigma^c$ on $N$, we let $X_{\Sigma^c}$ denote the $G/H$-horospherical variety corresponding to $\Sigma^c$ via \cref{thm:classification of horospherical varieties} (this is defined up to $G$-equivariant isomorphism). 
\end{notation}

\begin{remark}\label{rmk:simple open cover corresponds to coloured cones of coloured fan}
	Let $\Sigma^c$ be a coloured fan on $N$. Then $X_{\Sigma^c}$ is covered by the open simple $G/H$-horospherical varieties $X_{\sigma^c}$ as $\sigma^c$ ranges through the coloured cones in $\Sigma^c$; compare with \cref{prop:open cover by simples}. In fact, $X_{\Sigma^c}$ is covered by the $X_{\sigma^c}$ for $\sigma^c\in\Sigma^c_{\max}$. 
\end{remark}

\begin{example}\label{ex:SL_2/U_2 coloured fans}
	Consider $G=\SL_2$ and $H=U_2$. Then the coloured lattice $N$ is isomorphic to $\Z$, we have $\calC=\{\alpha\}$, and $u_\alpha=e_1$. We describe all $\SL_2/U_2$-horospherical varieties by describing all coloured fans on $N$.
	
	Of course, the trivial coloured cone $0^c$ yields the horospherical variety $\SL_2/U_2\cong\A^2\setminus\{0\}$. Note that the colour divisor $D_\alpha$ of $\SL_2/U_2$ is $\{(0,y)\}\subseteq\A^2\setminus\{0\}$. It is easy to see that the following are all of the possible nontrivial coloured fans on $N$:
	
	\tikzitfig{all-SL-2-U-2-coloured-fans}
	
	The following table says what the corresponding variety is in each case:
	
	\begin{mdframed}
		\centering
		\begin{tabular}{|c||c|c|c|c|c|}
			\hline
			$\Sigma^c$ & (1) & (2) & (3) & (4) & (5)
			\\\hline &&&&& \\[-1em]
			$X_{\Sigma^c}$ & $\A^2$ & $\Bl_0\A^2$ & $\P^2\setminus\{(1:0:0)\}$ & $\P^2$ & $\Bl_{(1:0:0)}\P^2$
			\\\hline
		\end{tabular}
	\end{mdframed}

	The action of $\SL_2$ on $\A^2$ and $\P^2$ is described in \cref{ex:A^2 and P^2 are SL_2/U_2-varieties}, and the actions on the other varieties are induced by these actions. Recall that examples (1) and (2) are treated in more detail in \cref{ex:coloured cone for A^2 and blow up}. 
\end{example}

Given a coloured fan $\Sigma^c$ on $N$, the \define{support} of $\Sigma^c$ is the support of the underlying fan, i.e. $|\Sigma^c|:=|\Sigma|:=\cup_{\sigma\in\Sigma} \sigma\subseteq N_\R$. We say that $\Sigma^c$ is \define{complete} if $|\Sigma^c|=N_\R$. 

\begin{proposition}[Characterization of complete]\label{prop:complete horospherical varieties}
	A $G/H$-horospherical variety $X$ is complete (i.e. proper over $\Spec(k)$) if and only if $\Sigma^c(X)$ is complete. 
	
	\begin{proof}
		See \cite[Theorem 4.2]{knop1991luna}. 
	\end{proof}
\end{proposition}

\subsection{$G$-orbits}\label{subsec:G-orbits}

Just like the orbit-cone correspondence for toric varieties (see \cite[Theorem 3.2.6]{cox2011toric}), we have an orbit-\textit{coloured} cone correspondence for horospherical varieties. Before stating this theorem, we give a remark concerning the stabilizers of collections of colour divisors, which is relevant for calculating the dimensions of orbits in the theorem (and for the local analysis in \cref{subsec:affine local structure}). 

\begin{remark}[Stabilizers of colour divisors]\label{rmk:stabilizers of colour divisors}
	The stabilizer $G_{D_\alpha}$ of a colour divisor $D_\alpha$ is the parabolic subgroup $P_{S\setminus\{\alpha\}}^-\supseteq B^-$. Hence, if $\calF\subseteq\calC$, then the stabilizer of all colour divisors $D_\alpha$ with $\alpha\in \calC\setminus\calF=S\setminus (I\cup \calF)$ is the parabolic subgroup $\cap_{\alpha\in\calC\setminus\calF} P_{S\setminus\{\alpha\}}^-=P_{I\cup\calF}^-\supseteq B^-$. 
\end{remark}

\begin{theorem}[Orbit-coloured cone correspondence]\label{thm:orbit-coloured cone correspondence}
	Let $X$ be a $G/H$-horospherical variety. Then there is a bijection
	\begin{align*}
		\{\text{$G$-orbits in $X$}\} &\map \Sigma^c(X)
		\\ \calO &\longmapsto \sigma^c(X_{\calO,G})
	\end{align*}
	which satisfies the following properties. Let $\calO(\sigma^c)$ denote the $G$-orbit corresponding to $\sigma^c\in\Sigma^c(X)$. 
	\begin{enumerate}
		\item For each $\sigma^c=(\sigma,\calF)\in\Sigma^c(X)$, we have
		\begin{align*}
			\dim(\calO(\sigma^c)) = \dim_\R(N_\R)-\dim(\sigma)+\dim(G/P_{I\cup\calF}).
		\end{align*}
		
		\item For each $\sigma^c\in\Sigma^c(X)$, we have
		\begin{align*}
			X_{\sigma^c} = \bigcup_{\tau^c\subseteq\sigma^c} \calO(\tau^c)
		\end{align*}
		where the union is over all coloured faces $\tau^c$ of $\sigma^c$. 
		
		\item For each $\sigma^c,\tau^c\in\Sigma^c(X)$, we have that $\tau^c$ is a coloured face of $\sigma^c$ if and only if $\calO(\sigma^c)\subseteq\ol{\calO(\tau^c)}$. Moreover, for each $\tau^c\in\Sigma^c(X)$, we have
		\begin{align*}
			\ol{\calO(\tau^c)} = \bigcup_{\tau^c\subseteq\sigma^c} \calO(\sigma^c)
		\end{align*}
		where the union is over all coloured cones $\sigma^c\in\Sigma^c(X)$ which contain $\tau^c$.
	\end{enumerate}
	
	\begin{proof}
		The bijective correspondence follows from \cref{lemma:G-orbits and coloured faces}. The proof in \cite[Lemma 3.2]{knop1991luna} also gives (3). The formula in (1) is shown in \cite[Theorem 6.6]{knop1991luna}. Finally, (2) follows from the description of simple horospherical varieties and their associated coloured cones. 
	\end{proof}
\end{theorem}

In \cref{thm:orbit-coloured cone correspondence}, the trivial coloured cone $0^c\in\Sigma^c$ corresponds to the open orbit $G/H\subseteq X_{\Sigma^c}$. On the other hand, the maximal coloured cones $\sigma^c\in\Sigma^c_{\max}$ correspond to the closed $G$-orbits $\calO(\sigma^c)=\ol{\calO(\sigma^c)}$ by part (3). 

We can use \cref{thm:orbit-coloured cone correspondence} to identify the $G$-invariant prime divisors (which are closures of codimension $1$ $G$-orbits) of a horospherical variety through the associated coloured fan. Given a coloured fan $\Sigma^c$ on $N$, a coloured ray of the form $\rho^c=(\rho,\varnothing)\in\Sigma^c(1)$ with empty colour set is called a \define{non-coloured ray}. 

\begin{corollary}\label{cor:non-coloured rays and G-invariant divisors}
	Let $X$ be a $G/H$-horospherical variety. There is a bijection 
	\begin{align*}
		\{\text{Non-coloured rays of $\Sigma^c(X)$}\} &\map \calD^G(X)
		\\ \rho^c=(\rho,\varnothing) &\longmapsto D_\rho:=\ol{\calO(\rho^c)}.
	\end{align*}
	Furthermore, the minimal ray generator $u_\rho\in N$ of $\rho^c=(\rho,\varnothing)\in\Sigma^c(1)$ equals $u_{D_\rho}$. 
	
	\begin{proof}
		This follows from \cref{thm:orbit-coloured cone correspondence}, and \cite[Lemma 2.4]{knop1991luna} for the final claim.
	\end{proof}
\end{corollary}

\begin{example}\label{ex:orbits of projective SL_3/U_3 variety}
	Consider $G=\SL_3$ and $H=U_3$. The coloured lattice $N$ is isomorphic to $\Z^2$, we have $\calC=\{\alpha_1,\alpha_2\}$, and $u_{\alpha_i}=e_i$ are the standard basis vectors for $i=1,2$. Let $X=X_{\Sigma^c}$ be the $\SL_3/U_3$-horospherical variety where $\Sigma^c$ is given in the diagram below. The maximal coloured cones in $\Sigma^c$ are $\sigma_1^c:=(\Cone(e_1,e_2),\{\alpha_1\})$, $\sigma_2^c:=(\Cone(e_2,-e_1-e_2),\varnothing)$, and $\sigma_3^c:=(\Cone(e_1,-e_1-e_2),\{\alpha_1\})$. 
	
	\tikzitfig{projective-SL-3-U-3-variety-with-alpha-1}
	
	We use \cref{thm:orbit-coloured cone correspondence} to give the dimensions of the $\SL_3$-orbits of $X$. Recall that $I=\varnothing$ since $B_3=P_\varnothing$ is the parabolic associated to $H=U_3$. Note that $\dim(\SL_3/P_\varnothing)=3$ and $\dim(\SL_3/P_{\{\alpha_1\}})=\dim(\SL_3/P_{\{\alpha_2\}})=2$. There are seven coloured cones in $\Sigma^c$, so there are seven $\SL_3$-orbits. The following table shows the dimensions of the orbits.
	
	\begin{mdframed}
		\centering
		\begin{tabular}{|c||c|c|c|c|c|c|c|}
			\hline
			$\sigma^c\in\Sigma^c$ & $\sigma_1^c$ & $\sigma_2^c$ & $\sigma_3^c$ & $\sigma_1^c\cap\sigma_2^c$ & $\sigma_1^c\cap\sigma_3^c$ & $\sigma_2^c\cap\sigma_3^c$ & $0^c$
			\\\hline &&&&&&& \\[-1em]
			$\dim(\calO(\sigma^c))$ & $2$ & $3$ & $2$ & $4$ & $3$ & $4$ & $5$
			\\\hline
		\end{tabular}
	\end{mdframed}
\end{example}

Let $X$ be a $G/H$-horospherical variety with associated coloured fan $\Sigma^c$. The closure of a $G$-orbit in $X$ has the structure of a horospherical $G$-variety, and we can describe the associated coloured fan as follows. Let $\calO(\tau^c)\subseteq X$ be a $G$-orbit for some $\tau^c\in\Sigma^c$. The sublattice $\Z(\tau\cap N)$ is saturated in $N$, so using \cref{ex:quotient coloured lattice} we can form the quotient coloured lattice, which we denote $N/\tau^c$. The underlying lattice is $N/\Z(\tau\cap N)$, the universal colour set is $\calC\setminus\calF(\tau^c)$, and the quotient map $\phi:N\to N/\tau^c$ is a map of coloured lattices. For each $\sigma^c\in\Sigma^c$ such that $\sigma\supseteq\tau$, let $(\sigma/\tau)^c$ be the coloured cone on $N/\tau^c$ whose underlying cone is $\sigma/\tau$, the image of $\sigma$ under $\phi_\R$, and whose colour set is $\calF(\sigma^c)\setminus\calF(\tau^c)$. It is easy to check that $(\sigma/\tau)^c$ is a strongly convex coloured cone, and the collection 
\begin{align}\label{eq:quotient coloured fan}
	(\Sigma/\tau)^c := \{(\sigma/\tau)^c : \tau^c\subseteq\sigma^c\in\Sigma^c\}
\end{align}
is a coloured fan on $N/\tau^c$. 

\begin{proposition}[Orbit closures]\label{prop:coloured fan of orbit closure}
	Let $X_{\Sigma^c}$ be a $G/H$-horospherical variety, and let $\calO(\tau^c)\subseteq X_{\Sigma^c}$ be a $G$-orbit for some $\tau^c\in\Sigma^c$. The orbit closure $\ol{\calO(\tau^c)}\subseteq X_{\Sigma^c}$ is a horospherical $G$-variety whose associated coloured fan is $(\Sigma/\tau)^c$. 
	
	\begin{proof}
		This mostly follows from \cref{thm:orbit-coloured cone correspondence}. The fact that $\ol{\calO(\tau^c)}$ is normal uses, for example, \cite[Corollary 2.3.4]{perrin2014geometry} (and $\charac(k)=0$). 
	\end{proof}
\end{proposition}

\begin{remark}\label{rmk:orbits as horospherical homogeneous spaces}
	Let $X$ be a $G/H$-horospherical variety and take $\tau^c\in\Sigma^c(X)$. Then the $G$-orbit $\calO(\tau^c)$ is $G$-equivariantly isomorphic to the horospherical homogeneous space $G/H_{(J,M)}$ determined by the pair $(J,M)=(\calF(\tau^c),(N/\tau^c)^\vee)$. 
\end{remark}

\begin{example}\label{ex:orbit closure coloured fan in projective SL_3/U_3 variety}
	Consider $G=\SL_3$ and $H=U_3$. The coloured lattice $N$ is isomorphic to $\Z^2$, we have $\calC=\{\alpha_1,\alpha_2\}$, and $u_{\alpha_i}=e_i$ are the standard basis vectors for $i=1,2$. Let $X=X_{\Sigma^c}$ be the $\SL_3/U_3$-horospherical variety from \cref{ex:orbits of projective SL_3/U_3 variety}. 
	
	Using \cref{prop:coloured fan of orbit closure}, the coloured fans for the $\SL_3$-orbit closures $\ol{\calO(\sigma_1^c)}$ (left), $\ol{\calO(\sigma_1^c\cap\sigma_2^c)}$ (middle), and $\ol{\calO(\sigma_1^c\cap\sigma_3^c)}$ (right) are given in the diagram below. 
	
	\tikzitfig{orbit-closure-coloured-fans-for-projective-SL-3-U-3-variety}
	
	Note that $\ol{\calO(\sigma_1^c)}$ is the flag variety $\SL_3/H_{(\{\alpha_1\},\{0\})}=\SL_3/P_{\{\alpha_1\}}$, $\ol{\calO(\sigma_1^c\cap\sigma_2^c)}$ is a $\SL_3/H_{(\varnothing,0\times\Z)}$-horospherical variety, and $\ol{\calO(\sigma_1^c\cap\sigma_3^c)}$ is a $\SL_3/H_{(\varnothing,\Z\times 0)}$-horospherical variety.
\end{example}

\subsection{Horospherical morphisms and maps of coloured fans}\label{subsec:horospherical morphisms and maps of coloured fans}

In this subsection, we consider two horospherical homogeneous spaces $G/H_i$ for $i=1,2$, and we assume that $H_1\subseteq H_2$. We use the same notation developed previously for $G/H$ but now with subscripts $i=1,2$; e.g. $N_i=N(G/H_i)$ and $\calC_i=\calC(G/H_i)$. 

As in \cref{subsec:maps of coloured lattices}, the natural projection map $G/H_1\to G/H_2$ induces an associated map of coloured lattices $\Phi:N_1\to N_2$. Recall that $\calC_\Phi\subseteq\calC_1$ is the set of dominantly mapped colours. The following theorem says that the map $G/H_1\to G/H_2$ extends to a horospherical morphism of $G/H_i$-horospherical varieties precisely when $\Phi$ respects the structure of the associated coloured fans. 

\begin{definition}[Compatible coloured fan map]\label{def:compatible map of coloured fans}
	Given coloured fans $\Sigma_i^c$ on $N_i$ for $i=1,2$, we say that $\Phi:N_1\to N_2$ is \define{compatible with $\Sigma_1^c$ and $\Sigma_2^c$} if, for each $\sigma_1^c=(\sigma_1,\calF_1)\in\Sigma_1^c$, there exists $\sigma_2^c=(\sigma_2,\calF_2)\in\Sigma_2^c$ such that $\Phi_\R(\sigma_1)\subseteq\sigma_2$ and $\calF_1\setminus\calC_\Phi\subseteq\calF_2$. 
\end{definition}

\begin{theorem}[Classification of horospherical morphisms]\label{thm:classification of horospherical morphisms}
	Let $X_{\Sigma_i^c}$ be $G/H_i$-horospherical varieties for $i=1,2$. Then the natural projection $G/H_1\to G/H_2$ extends to a horospherical morphism $X_{\Sigma_1^c}\to X_{\Sigma_2^c}$ if and only if $\Phi$ is compatible with $\Sigma_1^c$ and $\Sigma_2^c$. 
	
	\begin{proof}
		See \cite[Theorem 4.1]{knop1991luna}. 
	\end{proof}
\end{theorem}

\begin{remark}\label{rmk:integral submersions}
	There is a method for studying $G/H_i$-horospherical varieties without referencing the second $G/H_2$; this uses \textit{generalized coloured fans}, which are the same as a coloured fans except that we do not require the elements to be strongly convex (compare with \textit{generalized fans} in toric geometry \cite[Definition 6.2.2]{cox2011toric}). We do not go into detail on this topic in this paper; see the discussion of integral submersions in \cite[Section 4]{knop1991luna} for details. This approach is used in \cite{brion1993mori} to combinatorially describe the minimal model program for (horo)spherical varieties.
\end{remark}

\begin{example}[Open horospherical subvarieties]\label{ex:open horospherical subvarieties and sub-coloured fans}
	Let $X=X_{\Sigma^c}$ be a $G/H$-horospherical variety. If $X'\subseteq X$ is a $G$-invariant open subvariety, then $X'$ is also a $G/H$-horospherical variety. The coloured fan $\Sigma^c(X')$ is a \define{sub-coloured fan} of $\Sigma^c$, i.e. it is a coloured fan which is a subset of $\Sigma^c$. Conversely, any sub-coloured fan $(\Sigma')^c$ of $\Sigma^c$ yields an open horospherical subvariety $X_{(\Sigma')^c}\subseteq X_{\Sigma^c}$. Indeed, the identity map $N\to N$ is compatible with $(\Sigma')^c$ and $\Sigma^c$, and the induced map $X_{(\Sigma')^c}\to X_{\Sigma^c}$ is a $G$-equivariant open embedding. 
\end{example}

Recall from \cref{rmk:equivariant automorphisms of G/H} that $T_N=P/H$ is naturally isomorphic to $\Aut^G(G/H)$, the group of $G$-equivariant automorphisms of the horospherical homogeneous space $G/H$. The following result says that all $G$-equivariant automorphisms of $G/H$ extend to any $G/H$-horospherical variety. In particular, $T_N$ acts on any $G/H$-horospherical variety by $G$-equivariant automorphisms. 

\begin{corollary}\label{cor:equivariant automorphisms of horospherical variety}
	Let $X$ be a $G/H$-horospherical variety. Then $\Aut^G(X)=\Aut^G(G/H)$. 
	
	\begin{proof}
		See \cite[Lemma 4.1]{altmann2015merging}. 
	\end{proof}
\end{corollary}

Lastly, we give a combinatorial characterization of proper morphisms, which generalizes \cref{prop:complete horospherical varieties}. 

\begin{proposition}[Characterization of proper]\label{prop:proper horospherical morphisms}
	Let $X_i=X_{\Sigma_i^c}$ be $G/H_i$-horospherical varieties for $i=1,2$, and let $\pi:X_1\to X_2$ be a horospherical morphism. Then $\pi$ is proper if and only if $\Phi_\R^{-1}(|\Sigma_2^c|)=|\Sigma_1^c|$. 
	
	\begin{proof}
		See \cite[Theorem 4.2]{knop1991luna}. 
	\end{proof}
\end{proposition}

\subsection{Torus factors}\label{subsec:torus factors--overview}

To finish this section, we examine the notion of ``torus factors" for horospherical varieties, extending \cite[Proposition 3.3.9]{cox2011toric} for toric varieties.

\begin{definition}[Torus factor]\label{def:torus factor}
	Let $X$ be a $G/H$-horospherical variety with coloured fan $\Sigma^c$. We say that $X$ has a \define{torus factor} if $X$ is $G$-equivariantly isomorphic to $X_1\times T_1$ where $X_1$ is a horospherical $G$-variety of smaller dimension than $X$ and $T_1$ is a nontrivial torus. 
\end{definition}

The following proposition gives a nice combinatorial interpretation of torus factors. Note that, for a $G/H$-horospherical variety $X$, the set $\{u_D:D\in\calD(X)\}$ is the set of all colour points in $N$ together with all the minimal generators for the non-coloured rays of $\Sigma^c(X)$. 

\begin{proposition}[Characterizations of torus factors]\label{prop:characterizations of torus factors}
	Let $X$ be a $G/H$-horospherical variety with coloured fan $\Sigma^c$. The following conditions are equivalent:
	\begin{enumerate}
		\item $X$ has a torus factor. 
		\item $\scrO_X^*(X)\neq k^*$, i.e. $X$ has a nonconstant function $X\to k^*$. 
		\item $\{u_D:D\in\calD(X)\}$ does not span $N_\R$. 
	\end{enumerate}
	
	\begin{proof}
		(1)$\Rightarrow$(2) If $X\cong X_1\times T_1$ in the notation of \cref{def:torus factor}, where $T_1$ is a nontrivial torus, then $k^*\neq \scrO_{T_1}^*(T_1)\hookrightarrow \scrO_X^*(X)$. 
		
		(2)$\Rightarrow$(3) See \cite[Proposition A.1]{monahan2025horospherical-stack} for an equivalence of (2) and (3). 
		
		(3)$\Rightarrow$(1) Suppose that $\{u_D:D\in\calD(X)\}$ does not span $N_\R$, so $\Span_\R\{u_D:D\in\calD(X)\}\cap N$ is a proper saturated sublattice $N'$ of $N$. Let $N''$ be a nontrivial complementary sublattice, i.e. $N=N'\times N''$. Then we have a projection $N\to N'$, which induces an inclusion $(N')^\vee\subseteq N^\vee$, which in turn induces a containment of horospherical subgroups $H=H_{(I,N^\vee)}\subseteq H_{(I,(N')^\vee)}=:H'$ (using \cref{prop:classification of horospherical subgroups}). Notice that $\Sigma^c$ can be viewed as a coloured fan on either $N$ or $N'$ because $\Sigma^c$ is contained in $N'_\R$. In the former case, $\Sigma^c$ corresponds the the $G/H$-horospherical variety $X=X_{\Sigma^c,N}$, but in the latter case, it corresponds to a $G/H'$-horospherical variety $X_{\Sigma^c,N'}$. Finally, we have 
		\begin{align*}
			X = X_{\Sigma^c,N} = X_{\Sigma^c,N'}\times T_{N''}
		\end{align*}
		where $T_{N''}$ is a nontrivial torus. Thus, $X$ has a torus factor. 
	\end{proof}
\end{proposition}

\begin{remark}[Removing torus factors]\label{rmk:removing torus factors}
	In the notation of \cref{prop:characterizations of torus factors}, let $N(\Sigma^c):=\Span_\R\{u_D:D\in\calD(X)\}\cap N$, which is a saturated sublattice of $N$; note that $N(\Sigma^c)_\R=\Span_\R\{u_D:D\in\calD(X)\}$ is the subspace of $N_\R$ spanned by all rays of $\Sigma$ together with all colour points in $N$. Then \cref{prop:characterizations of torus factors} says that we have $N(\Sigma^c)=N$ if and only if $X=X_{\Sigma^c}$ has no torus factors. 
	
	Assuming that $N(\Sigma^c)\neq N$, the proof of \cref{prop:characterizations of torus factors} shows us that
	\begin{align}\label{eq:removing torus factor}
		X=X_{\Sigma^c,N} = X_{\Sigma^c,N(\Sigma^c)}\times T_{N/N(\Sigma^c)}
	\end{align}
	where $T_{N/N(\Sigma^c)}$ is the nontrivial torus with one-parameter subgroup lattice $N/N(\Sigma^c)$. Now $X_{\Sigma^c,N(\Sigma^c)}$ is a $G/H_{(I,N(\Sigma^c)^\vee)}$-horospherical variety which has no torus factors.
\end{remark}

\begin{example}\label{ex:torus factor example}
	Consider $G=\SL_2\times\G_m$ and $H=U_2\times\{1\}$. Then the coloured lattice $N$ is isomorphic to $\Z^2$, we have $\calC=\{\alpha\}$, and $u_\alpha=e_1$. Let $X=X_{\Sigma^c}$ be the $G/H$-horospherical variety where $\Sigma^c$ is given in the diagram below. The unique maximal coloured cone in $\Sigma^c$ is $(\Cone(e_1),\{\alpha\})$. 
	
	\tikzitfig{coloured-fan-of-torus-factor-SL-2xG-m-variety}
	
	As in \cref{rmk:removing torus factors}, $N(\Sigma^c)=\Z$ is the sublattice of $N=\Z^2$ spanned by $e_1$. It is easy to see that $N(\Sigma^c)$ is the coloured lattice associated to $\SL_2/U_2$. It follows that
	\begin{align*}
		X &= X_{\Sigma^c,N} = X_{\Sigma^c,N(\Sigma^c)}\times T_{N/N(\Sigma^c)} = \A^2\times\G_m.
		\qedhere
	\end{align*}
\end{example}

\section{Affine horospherical varieties and local structure}\label{sec:affine horospherical varieties and local structure}

In this section we characterize affine horospherical varieties, give a way of computing the coordinate ring of an affine horospherical variety, illustrate a certain local description for a general horospherical variety using affine horospherical varieties, and state combinatorial characterizations of factoriality and smoothness. 

Throughout this section, we fix a horospherical homogeneous space $G/H$ and use the notation from the previous sections. In particular, $P=P_I=N_G(H)\supseteq B$ is the associated parabolic, and $N=N(G/H)$ is the associated coloured lattice with universal colour set $\calC=\calC(G/H)=S\setminus I$.

\subsection{Affine horospherical varieties}\label{subsec:affine horospherical varieties}

As mentioned in \cref{subsec:simple open cover}, every simple \textit{toric} variety is affine, but not every simple \textit{horospherical} variety is affine. The following proposition tells us exactly which horospherical varieties are affine. 

\begin{proposition}[Characterization of affine]\label{prop:characterization of affine}
	Let $X$ be a $G/H$-horospherical variety. Then the following are equivalent:
	\begin{enumerate}
		\item $X$ is affine.
		\item $X$ is simple and $X=X_{\calO,B^-}$, where $\calO$ denotes the unique closed $G$-orbit of $X$. 
		\item $X$ is simple and $\calF(X)=\calC$. 
	\end{enumerate}
	
	\begin{proof}
		See \cite[Theorem 6.7]{knop1991luna} for the equivalence of (1) and (3). Recall that $\calV(X)=N_\R$ since $X$ is horospherical, so the only functional which is nonpositive on $\calV(X)$ is the zero functional. The equivalence between (2) and (3) follows from the definition of $X_{\calO,B^-}$.
	\end{proof}
\end{proposition}

Consider an affine $G/H$-horospherical variety $X$; so $X$ corresponds to a coloured cone $\sigma^c=(\sigma,\calF)$ on $N$ such that $\calF=\calC$. Due to the condition on the colour set for an affine horospherical variety, we see that $X$ is completely determined by the cone $\sigma$. We can use $\sigma$ to describe the coordinate ring of $X$ as follows. Similar to \cref{subsec:eigenfunctions}, we can consider the regular $B^-$-eigenfunctions $k[G]^{(B^-)}$ of $k[G]$, which has a $\frakX(T)$-grading:
\begin{align*}
	k[G]^{(B^-)}\cup\{0\} = \bigoplus_{m\in\frakX(T)} k[G]^{(B^-)}_{-m}
\end{align*}
where
\begin{align*}
	k[G]^{(B^-)}_{-m} &:= \{f\in k[G] : b\cdot f=\chi^{-m}(b)f ~\forall b\in B^-\}
	\\&= \{f\in k[G] : f(b\cdot x)=\chi^m(b)f(x) ~\forall b\in B^-, x\in G\}.
\end{align*}
Since $N^\vee\subseteq\frakX(T)$, we can consider the graded pieces $k[G]^{(B^-)}_{-m}$ for $m\in N^\vee$. 

For example, if $G=T$ is a torus, then $k[G]=k[T]$ is the ring of Laurent monomials (in the number of variables equal to $\dim(T)$) and $k[G]^{(B^-)}_{-m}=k[T]^{(T)}_{-m}$ consists of Laurent monomials of degree $m$. 

\begin{proposition}[Coordinate ring]\label{prop:coordinate ring}
	Let $X=X_{\sigma^c}$ be an affine $G/H$-horospherical variety. The coordinate ring of $X$ is given by 
	\begin{align*}
		k[X] = \bigoplus_{m\in \sigma^\vee\cap N^\vee} k[G]^{(B^-)}_{-m}.
	\end{align*}

	\begin{proof}
		See \cite[Theorem 28.3]{timashev2011homogeneous}. 
	\end{proof}
\end{proposition}

\begin{example}\label{ex:coordinate ring of toric variety}
	Let $X$ be an affine toric variety with open torus $T$, i.e. $X$ is an affine $T/\{1\}$-horospherical variety. Then $X$ corresponds to a (strongly convex) cone $\sigma$ and the associated coloured cone is $\sigma^c(X)=(\sigma,\varnothing)$. As in \cite[Theorem 1.3.5]{cox2011toric}, the coordinate ring of $X$ is the monoid algebra $k[X]=k[\sigma^\vee\cap N^\vee]$. This agrees with \cref{prop:coordinate ring} since $k[T]^{(T)}_{-m}$ is the space of Laurent monomials of degree $m\in N^\vee$. 
\end{example}

\begin{example}\label{ex:coordinate ring of affine SL_3/U_3-variety}
	Consider $G=\SL_3$ and $H=U_3$. Then the coloured lattice $N$ is isomorphic to $\Z^2$, we have $\calC=\{\alpha_1,\alpha_2\}$, and $u_{\alpha_i}=e_i$ are the standard basis vectors for $i=1,2$. Let $X$ be the affine $\SL_3/U_3$-horospherical variety whose associated coloured cone is $(\Cone(e_1,e_2),\calC)$; compare with \cref{ex:SL_3/U_3 affine variety associated coloured cone}. We compute the coordinate ring $k[X]$ using \cref{prop:coordinate ring} as follows.
	
	First, we can view $\SL_3$ as a closed subvariety of $\A^9$ where the nine coordinates represent the matrix entries. Thus $k[\SL_3]$ is a quotient of $k[r,s,t,u,w,v,x,y,z]$; here $(r,s,t)$ represents the first row of entries in $\SL_3$, and so on. 
	
	To compute the direct sum in \cref{prop:coordinate ring}, we mainly need to compute $k[\SL_3]^{(B_3^-)}_{-m}$ for $m$ in a generating set of the monoid $\sigma^\vee\cap N^\vee$. Note that $\sigma^\vee\cap N^\vee=\Z_{\geq 0}^2$ where $N^\vee=\Z^2$, so $\{e_1,e_2\}$ is a generating set. The action of $B_3^-$ on $\SL_3$ tells us that
	\begin{align*}
		k[\SL_3]^{(B_3^-)}_{-e_1} = k\cdot\{r,s,t\} \quad\text{and}\quad k[\SL_3]^{(B_3^-)}_{-e_2} = k\cdot\{rv-su,rw-tu,sw-tv\}
	\end{align*}
	(note that these quantities $rv-su$, $rw-tu$, and $sw-tv$ are the upper $2\times 2$-minors of the matrices in $\SL_3$). Therefore, we have
	\begin{align*}
		k[X] = k[r,s,t,rv-su,rw-tu,sw-tv] \cong k[a,b,c,d,e,f]/\langle ad+be+cf\rangle.
		\qedhere
	\end{align*}
\end{example}

\subsection{Affine local structure}\label{subsec:affine local structure}

Consider a $G/H$-horospherical variety $X_{\Sigma^c}$. In this subsection we present a local description for $X$ in terms of affine horospherical varieties. 
First, we know that $X_{\Sigma^c}$ is covered by open simple horospherical varieties for each $\sigma^c\in\Sigma^c$, so to provide this local description we first reduce to one of these simple pieces. Let $X=X_{\sigma^c}\subseteq X_{\Sigma^c}$ be an open simple subvariety. As in the last subsection, we know that $X$ may not be affine, so we cannot simply use this simple open cover as an affine open cover of $X_{\Sigma^c}$. 

Let $Q$ be the parabolic subgroup $P_{I\cup\calF(X)}$; recall from \cref{rmk:stabilizers of colour divisors} that $Q^-$ is the stabilizer of the colour divisors coming from $\calC\setminus\calF(X)$. Then $P=P_I\subseteq Q$, so we have a natural projection $G/H\to G/P\to G/Q$. The universal colour set $\calC(G/Q)$ for $G/Q$ is $\calC\setminus\calF(X)$, so this projection $G/H\to G/Q$ extends to a horospherical morphism $X\to G/Q$ by \cref{thm:classification of horospherical morphisms}. Let $Z$ denote the fibre of the base point $eQ$; so $Z$ is a closed subvariety of $X$. Note that every fibre of $X\to G/Q$ is isomorphic to $Z$ since the map is $G$-equivariant and the base is a single $G$-orbit. In particular, we have $X=G\times^Q Z$; recall that the associated bundle $G\times^Q Z$ is $(G\times Z)/Q$ where $q\cdot (g,z):=(gq^{-1},q\cdot z)$. Since $eQ$ is a $Q$-fixed point in $G/Q$, we see that $Z$ is $Q$-invariant. Moreover, $\Rad_u(Q)$ acts on $Z$ trivially, so $Z$ is an $L$-variety where $L=L_{I\cup\calF(X)}$ is the standard Levi subgroup of $Q$. 

Notice that the coloured lattice $N(L/(L\cap H))$ satisfies the following properties. As a lattice, it is the same as $N(G/H)$ because $L/(L\cap H)$ and $G/H$ have the same associated torus. However, the colour structure is different: the universal colour set is $\calC(L/(L\cap H))=\calF(X)$. That is, using the terminology of \cref{ex:quotient coloured lattice}, $N(L/(L\cap H))$ is a coloured sublattice of $N(G/H)$, and the coloured quotient lattice $N(G/H)/N(L/(L\cap H))$ is precisely $N(G/Q)$. 

\begin{theorem}[Affine local structure]\label{thm:affine local structure}
	Let $X=X_{\sigma^c}$ be a simple $G/H$-horospherical variety. In the notations immediately above, we have the following:
	\begin{enumerate}
		\item $Z$ is an affine $L/(L\cap H)$-horospherical variety whose associated coloured cone is $\sigma^c$ on the coloured lattice $N(L/(L\cap H))$.
		\item If $\calO\subseteq X$ denotes the unique closed $G$-orbit, then we have $Z\subseteq X_{\calO,B^-}$ and the map $\Rad_u(Q^-)\times Z\to X_{\calO,B^-}$ via $(q,z)\mapsto q\cdot z$ is a $Q^-$-equivariant isomorphism. 
	\end{enumerate}
	
	\begin{proof}
		See \cite[Theorem 28.2]{timashev2011homogeneous} for (1) and \cite[Theorem 3.2.2]{perrin2018sanya} for (2). 
	\end{proof}
\end{theorem}

\begin{remark}\label{rmk:comment about affine local structure fixed point}
	One can get a similar result to \cref{thm:affine local structure} where the affine closed subvariety $Z\subseteq X$ has a fixed point; see \cite[Theorem 28.2]{timashev2011homogeneous} for details.
\end{remark}

\begin{example}\label{ex:affine local structure for SL_3 variety}
	Consider $G=\SL_3$ and $H=U_3$. Then the coloured lattice $N$ is isomorphic to $\Z^2$, we have $\calC=\{\alpha_1,\alpha_2\}$, and $u_{\alpha_i}=e_i$ are the standard basis vector for $i=1,2$. Let $X$ be the simple $\SL_3/U_3$-horospherical variety whose associated coloured cone is $\sigma^c(X)=(\Cone(e_1,e_2),\{\alpha_1\})$; see the coloured cone on the left in the diagram below. 
	
	In the notation at the beginning of this subsection, we have $Q=P_{\{\alpha_1\}}$ since $I=\varnothing$ and $\calF(X)=\{\alpha_1\}$, so $L=L_{\{\alpha_1\}}$. Note that, using \cref{rmk:ss and torus isogeny horospherical varieties}, we can think of $L/(L\cap U_3)$ as $(\SL_2\times\Gm)/(U_2\times\{1\})$. Then $Z$ is the $L/(L\cap U_3)$-horospherical variety whose associated coloured cone is $\sigma^c(Z)=(\Cone(e_1,e_2),\{\alpha_1\})$ on the coloured lattice $\Z^2$ with universal colour set $\{\alpha_1\}$; see the coloured cone on the right in the diagram below. 
	
	\tikzitfig{example-of-affine-reduction-1}
	
	We can describe $Z$ as $\A^2\times\A^1=\A^3$ as a $(\SL_2\times\Gm)/(U_2\times\{1\})$-horospherical variety, where $\SL_2$ acts on the $\A^2$ by matrix multiplication and $\G_m$ acts on the $\A^1$ by multiplication. Then $X=\SL_3\times^Q \A^3$. 
\end{example}

\subsection{Toroidal horospherical varieties}\label{subsec:toroidal horospherical varieties}

Now we look at a special class of horospherical varieties whose affine local structure (from \cref{thm:affine local structure}) uses toric varieties. These varieties are particularly nice because one can study all of their local properties (e.g. singularities) using toric geometry.

\begin{definition}[Toroidal]\label{def:toroidal}
	A $G/H$-horospherical variety $X$ is called \define{toroidal} if $\calF(X)=\varnothing$, or equivalently, if the projection $G/H\to G/P$ extends to a horospherical morphism $X\to G/P$ (by \cref{thm:classification of horospherical morphisms}). 
\end{definition}

For example, toric varieties and flag varieties are toroidal. Notice that being toroidal means that, for each open simple subvariety, we have $Q=P$ and $L/(L\cap H)=T_N$, where $Q$ and $L$ are defined in \cref{subsec:affine local structure}.
Therefore, for a toroidal $G/H$-horospherical variety $X$, we can proceed as at the beginning of \cref{subsec:affine local structure} to get a \textit{global} fibre $Z$ of the map $X\to G/P$ which has an action by $T_N$. Importantly, we do not need to assume that $X$ is simple. 

\begin{corollary}\label{cor:affine local structure of toroidal}
	Let $X=X_{\Sigma^c}$ be a toroidal $G/H$-horospherical variety. If $Z\subseteq X$ denotes the fibre over $eP$ of the map $X\to G/P$, then $X\cong G\times^P Z$ and $Z$ is a toric variety with open torus $T_N$ whose associated fan is $\Sigma$. In particular, if $X=X_{\sigma^c}$ is simple, then $Z$ is an affine toric variety. 
\end{corollary}

\begin{example}\label{ex:open toroidal subvariety}
	Let $X=X_{\Sigma^c}$ be a $G/H$-horospherical variety. We describe a canonical open toroidal subvariety of $X$ as follows. Let $X'\subseteq X$ denote the open subset where we remove all $G$-orbits in $X$ of codimension greater than $1$. Based on \cref{thm:orbit-coloured cone correspondence}, $X'$ is a $G/H$-horospherical variety whose coloured fan consists of the trivial coloured cone and the non-coloured rays in $\Sigma^c$. In particular, $X'$ is toroidal. 
\end{example}

\begin{example}[Decolouration]\label{ex:decolouration}
	Let $X=X_{\Sigma^c}$ be a $G/H$-horospherical variety. We describe a canonical toroidal resolution of $X$ as follows. Consider the coloured fan $\wt\Sigma^c$ whose underlying fan is exactly $\Sigma$ and whose colour set is empty, i.e. $\calF(\wt\Sigma^c)=\varnothing$. This yields a toroidal $G/H$-horospherical variety $\wt X=X_{\wt\Sigma^c}$, which we call the \define{decolouration} of $X$. Moreover, the identity map of coloured lattices $N\to N$ is compatible with $\wt\Sigma^c$ and $\Sigma^c$, so we have a horospherical morphism $\wt X\to X$; note that this map is surjective and proper. 
	
	For example, $\Bl_0\A^2$ is the decolouration of $\A^2$ where both are viewed as $\SL_2/U_2$-horospherical varieties; see \cref{ex:coloured cone for A^2 and blow up}. 
\end{example}

\subsection{Factoriality and smoothness}\label{subsec:factoriality and smoothness}

We finish this section with combinatorial descriptions of two local properties: factoriality and smoothness. Given a coloured fan $\Sigma^c$ on $N$ and $\sigma^c\in\Sigma^c$, consider the multiset
\begin{align}\label{eq:multiset for cone regularity}
	\{u_\rho:(\rho,\varnothing)\in\sigma^c(1)\}\cup\{u_\alpha:\alpha\in\calF(\sigma^c)\}
\end{align}
of lattice points in $N$; recall that these $u_\rho$ are the minimal ray generators for the non-coloured rays $(\rho,\varnothing)\in\sigma^c(1)$. Note that we say \textit{multiset} here because we want to count colour points with multiplicity. We say that $\sigma^c$ is \define{simplicial} (resp. \define{regular}) if the multiset \cref{eq:multiset for cone regularity} is $\R$-linearly independent (resp. is part of a $\Z$-basis for $N$). We say that $\Sigma^c$ is \define{simplicial} (resp. \define{regular}) if each coloured cone in $\Sigma^c$ is simplicial (resp. regular). Clearly regular implies simplicial.

\begin{proposition}[Factoriality criteria]\label{prop:factoriality}
	Let $X_{\Sigma^c}$ be a $G/H$-horospherical variety. 
	\begin{enumerate}
		\item $X_{\Sigma^c}$ is factorial if and only if $\Sigma^c$ is regular.
		\item $X_{\Sigma^c}$ is $\Q$-factorial if and only if $\Sigma^c$ is simplicial. 
	\end{enumerate}

	\begin{proof}
		See \cite[Theorem 4.2.3]{perrin2018sanya}. 
	\end{proof}
\end{proposition}

\begin{example}\label{ex:non-factorial coloured cone with double colour point for SL_3}
	Consider $G=\SL_3$ and $H=H_{(\varnothing,M)}$ where $M=\Z$ is the diagonal sublattice inside $\frakX(T_3)=\Z^2$. Then the coloured lattice $N$ is isomorphic to $\Z$, we have $\calC=\{\alpha_1,\alpha_2\}$, and $u_{\alpha_i}=e_1$ for both $i=1,2$. 
	
	Let $\sigma^c=(\sigma,\calF)=(\Cone(e_1),\{\alpha_1,\alpha_2\})$; see the left part of the diagram below. The multiset \cref{eq:multiset for cone regularity} consists of two copies of $e_1$, so $\sigma^c$ is neither regular nor simplicial. Therefore, the corresponding simple $\SL_3/H$-horospherical variety $X_{\sigma^c}$ is not $\Q$-factorial by \cref{prop:factoriality}. Note that, if $\calF$ consisted of only one colour, say $\alpha_1$, then $\sigma^c$ would be regular; see the right part of the diagram below. 
	
	\tikzitfig{non-factorial-coloured-cone-with-double-colour-point-for-SL-3}
	
	Since $X_{\sigma^c}$ is affine (by \cref{prop:characterization of affine}), one can use \cref{prop:coordinate ring} to check that 
	\begin{align*}
		k[X_{\sigma^c}] &= \frac{k[ad,ae,af,bd,be,bf,cd,ce,cf]}{\langle ad+be+cf\rangle} 
		\\&\cong \frac{k[r,s,t,u,v,w,x,y,z]}{\Big\langle\begin{matrix}rv-su, rw-tu, ry-sx, rz-tx, sw-tv,\\ sz-ty, uy-vx, vz-wx, uz-wx, r+v+z\end{matrix}\Big\rangle}
	\end{align*}
	where the $a,b,c,d,e,f$ coordinates are the ones used in \cref{ex:coordinate ring of affine SL_3/U_3-variety}.
\end{example}

For those familiar with toric geometry, they will recognize that a regular fan (not a \textit{coloured} fan) is a smooth fan; see \cite[Definition 3.1.18]{cox2011toric}. In \cite[Theorem 3.1.19]{cox2011toric}, smooth fans are used to characterize smooth toric varieties. It is known that a toric variety is smooth if and only if it is factorial; so smooth fans really characterize \textit{factorial} toric varieties, which is more in-line with \cref{prop:factoriality}. For horospherical varieties, we need the following combinatorial condition, originally due to \cite[Definition 2.4]{pasquier2006thesis},\footnote{Originally, the term \textit{lisse} (translating to \textit{smooth}) is used in \cite[Definition 2.4]{pasquier2006thesis}. We choose to avoid \textit{lisse} in this context because this alone does not correspond to smoothness of horospherical varieties, and we want to avoid confusion with the notion of a smooth fan in toric geometry.} on the Dynkin diagram of $G$ to bridge the gap between factorial and smooth.

\begin{definition}[Vivid]\label{def:vivid}
	Let $X$ be a $G/H$-horospherical variety with coloured fan $\Sigma^c$. Recall that $P=N_G(H)$ corresponds to the subset $I$ of simple roots of $G$. For $J\subseteq S$, let $\Gamma_J$ be the sub-Dynkin diagram of $G$ whose vertices are the simple roots in $J$. We say that $\Sigma^c$ (or $X$ itself) is \define{vivid} if, for each $\sigma^c=(\sigma,\calF)\in\Sigma^c$ and each $\alpha\in\calF$, the following conditions hold:
	\begin{enumerate}
		\item The connected component $\Gamma^\alpha$ of $\Gamma_{I\cup\calF}$ containing $\alpha$ contains no other element of $\calF$. 
		\item $\Gamma^\alpha=\Gamma_{\{\alpha\}\cup I_\alpha}$ for some $I_\alpha\subseteq I$, and $\Gamma^\alpha$ has type $\mathbf{A}_n$ or $\mathbf{C}_n$ (for some $n$) in which $\alpha$ is the first simple root (for $\mathbf{C}_n$ this means the opposite end to the longest simple root). 
	\end{enumerate} 
\end{definition}

\begin{remark}\label{ex:vivid Dynkin diagrams}
	If $\Gamma^\alpha$ is one of the connected components from \cref{def:vivid}, then $\Gamma^\alpha=\Gamma_{\{\alpha\}\cup I_\alpha}$ has one of the two forms below:
	\begin{equation*}
		\begin{tikzpicture}
			\def\a{1}
			\tikzset{dynkin/.style={circle,draw,inner sep=0pt,minimum size=2mm}}
			\path
			(0,0)    node{$\mathbf{A}_n$:}  
			++(0:\a)	node[dynkin] (N1) {} +(-90:.5) node{$\alpha$}
			++(0:\a)   node[dynkin] (N2) {} 
			++(0:\a) node[dynkin] (N3) {}
			++(0:0.5*\a) coordinate (A) ++(0:\a) coordinate (B)
			++(0:0.5*\a) node[dynkin] (N4) {} ;
			
			\draw[dashed] (A)--(B);
			\draw (N1)--(N2)--(N3)--(A) (B)--(N4);
			\draw[decorate,decoration={brace,raise=3mm},thick]
			(N4.center)--(N2.center) node[midway,below=4mm]{$I_\alpha$};
			
			\path
			(7*\a,0)      node{$\mathbf{C}_n$:}
			++(0:\a)	node[dynkin] (M1) {} +(-90:.5) node{$\alpha$}
			++(0:\a)   node[dynkin] (M2) {} 
			++(0:\a) node[dynkin] (M3) {}
			++(0:0.5*\a) coordinate (C) ++(0:\a) coordinate (D)
			++(0:0.5*\a) node[dynkin] (M4) {} 
			++(0:0.5*\a) coordinate (E)
			++(0:0.5*\a) node[dynkin] (M5) {} ;
			
			\draw[dashed] (C)--(D);
			\draw (M1)--(M2)--(M3)--(C) (D)--(M4);
			\draw[double,double distance=1mm] (M4)--(M5);
			\draw[-{Classical TikZ Rightarrow[length=2mm]},double,double distance=1mm] (M5)--(E.west);
			\draw[decorate,decoration={brace,raise=3mm},thick]
			(M5.center)--(M2.center) node[midway,below=4mm]{$I_\alpha$};
		\end{tikzpicture}
	\end{equation*}
	If $G^\alpha$ is the semisimple simply connected group associated to this Dynkin diagram $\Gamma^\alpha$, then $G^\alpha=\SL_{n+1}$ or $G^\alpha=\opn{Sp}_{2n}$, respectively. Having one of the two forms above is equivalent to $G^\alpha$ acting transitively on projective space (of dimension $n$ or $2n-1$, respectively) where the stabilizer (up to conjugation) is the parabolic subgroup of $G^\alpha$ corresponding to $I_\alpha$. 
\end{remark}

\begin{example}[{cf. \cite[Example 2.5]{pasquier2006thesis}}]\label{ex:vivid examples}	
	Consider $G=\opn{Sp}_{14}\times T_0$ where $T_0$ is any torus, so $G$ has Dynkin type $\mathbf{C}_7$.
	\begin{equation*}
		\begin{tikzpicture}
			\def\a{1}
			\tikzset{dynkin/.style={circle,draw,inner sep=0pt,minimum size=2mm}}
			\path
			(0,0)   node[dynkin] (M1) {} +(-90:.5) node{$\alpha_1$}
			++(0:\a)   node[dynkin] (M2) {} +(-90:.5) node{$\alpha_2$}
			++(0:\a) node[dynkin] (M3) {} +(-90:.5) node{$\alpha_3$}
			++(0:\a) node[dynkin] (M4) {} +(-90:.5) node{$\alpha_4$}
			++(0:\a) node[dynkin] (M5) {} +(-90:.5) node{$\alpha_5$}
			++(0:\a) node[dynkin] (M6) {} +(-90:.5) node{$\alpha_6$}
			++(0:0.5*\a) coordinate (E)
			++(0:0.5*\a) node[dynkin] (M7) {} +(-90:.5) node{$\alpha_7$} ;
			
			\draw (M1)--(M2)--(M3)--(M4)--(M5)--(M6);
			\draw[double,double distance=1mm] (M6)--(M7);
			\draw[-{Classical TikZ Rightarrow[length=2mm]},double,double distance=1mm] (M7)--(E.west);
		\end{tikzpicture}
	\end{equation*}
	
	Let $\sigma^c=(\sigma,\calF)$ be a coloured cone on $N$ with $\calF\subseteq \calC=S\setminus I$. We give the following examples with $I$ and $\calF$:
	\begin{enumerate}[label=(\alph*)]
		\item If $I=\{\alpha_2,\alpha_6,\alpha_7\}$ and $\calF=\{\alpha_1,\alpha_5\}$, then $\sigma^c$ is vivid; note that $\Gamma^{\alpha_1}$ has type $\mathbf{A}_2$ and $\Gamma^{\alpha_5}$ has type $\mathbf{C}_3$. 
		\item If $I=\{\alpha_4,\alpha_6,\alpha_7\}$ and $\calF=\{\alpha_5\}$, then $\sigma^c$ is \textit{not} vivid because $\Gamma^{\alpha_5}$ has type $\mathbf{C}_4$ but $\alpha_5$ is not the first simple root. 
		\item If $I=\{\alpha_2,\alpha_5\}$ and $\calF=\{\alpha_1,\alpha_3\}$, then $\sigma^c$ is \textit{not} vivid because $\alpha_1$ and $\alpha_3$ are in the same connected component of $\Gamma_{I\cup\calF}$.  \qedhere
	\end{enumerate}
\end{example}

\begin{theorem}[Smoothness criterion]\label{thm:smoothness}
	Let $X_{\Sigma^c}$ be a $G/H$-horospherical variety. Then $X$ is smooth if and only if $\Sigma^c$ is both regular and vivid. 

	\begin{proof}
		See \cite[Theorem 28.10]{timashev2011homogeneous}. 
	\end{proof}
\end{theorem}

\begin{example}\label{ex:smoothness of some SL_4 varieties}
	Consider $G=\SL_5$ and $H$ the minimal horospherical subgroup with associated parabolic $P_{\{\alpha_2,\alpha_4\}}$ (see \cref{rmk:minimal horospherical subgroups}). Then the coloured lattice $N$ is isomorphic to $\Z^2$, we have $\calC=\{\alpha_1,\alpha_3\}$, and $u_{\alpha_1}=e_1$ and $u_{\alpha_3}=e_2$ are the standard basis vectors. In the table below we consider three simple $\SL_5/H$-horospherical varieties given by coloured cones $\sigma_i^c$ for $i=1,2,3$. We test smoothness of each variety $X_{\sigma_i^c}$ using \cref{thm:smoothness}.
	
	\begin{mdframed}
		\centering
		\begin{tabular}{|c||c|c|}
			\hline
			$i$ & $\sigma_i^c$ & smoothness
			\\\hline\hline && \\[-1em]
			$1$ & $(\Cone(e_1,-e_1+e_2),\calC)$ & NOT smooth
			\\\hline && \\[-1em]
			$2$ & $(\Cone(e_1,e_2),\{\alpha_1\})$ & smooth
			\\\hline && \\[-1em]
			$3$ & $(\Cone(e_1,e_2),\{\alpha_3\})$ & NOT smooth
			\\\hline
		\end{tabular}
	\end{mdframed}

	For $i=1$, the variety is not smooth because $\sigma_1^c$ is not regular. Indeed, the multiset \cref{eq:multiset for cone regularity} is $\{e_1,e_2,-e_1+e_2\}$, which clearly cannot be extended to a $\Z$-basis for $N$. Moreover, $\sigma_1^c$ is not vivid since the elements of $\calF(\sigma_1^c)=\calC$ are in the same connected component of the Dynkin diagram of $\SL_5$. 
	
	For $i=2,3$, it is clear that each $\sigma_i^c$ is regular since the multiset \cref{eq:multiset for cone regularity} is $\{e_1,e_2\}$. It remains to check vividness for each of these $\sigma_i^c$. 
	
	For $i=2$, it is vivid. Indeed, $\alpha_1$ is the only colour in $\calF(\sigma_2^c)$ and it is connected with the component $\{\alpha_1,\alpha_2\}$ of $I\cup\calF(\sigma_2^c)=\{\alpha_1,\alpha_2,\alpha_4\}$ so that $\{\alpha_1\}\cup\{\alpha_2\}$ has type $\mathbf{A}_2$ where $\alpha_1$ is the first simple root. 
	
	For $i=3$, it is not vivid. Indeed, $\alpha_3$ is not the first simple root in $I\cup\calF(\sigma_3^c)=\{\alpha_2,\alpha_3,\alpha_4\}$. 
\end{example}

For toroidal horospherical varieties, vividness in \cref{thm:smoothness} is vacuous, so we recover the smoothness criterion that one would expect from toric geometry. 

\begin{corollary}\label{cor:smoothness of toroidal}
	Let $X_{\Sigma^c}$ be a toroidal $G/H$-horospherical variety. Then $X_{\Sigma^c}$ is smooth if and only if $\Sigma^c$ is regular. 
\end{corollary}

\section{Weil and Cartier divisors}\label{sec:Weil and Cartier divisors}

To end this paper, we discuss Weil and Cartier divisors on horospherical varieties. In particular, we show how to compute the class group and Picard group using the $B^-$-invariant divisors and the coloured fan combinatorics.

Throughout this section, we fix a horospherical homogeneous space $G/H$ and use the notation from the previous sections. In particular, $P=P_I=N_G(H)\supseteq B$ is the associated parabolic, and $N=N(G/H)$ is the associated coloured lattice with universal colour set $\calC=\calC(G/H)=S\setminus I$.

\subsection{The class group}\label{subsec:the class group}

Given a $G/H$-horospherical variety $X$, let $\Div(X)$ denote the group of Weil divisors on $X$, i.e. the free abelian group generated by the prime divisors of $X$. Recall that the class group of $X$, denoted $\Cl(X)$, is the quotient of $\Div(X)$ by the subgroup generated by principal divisors, i.e. divisors of the form $\prindiv(f)$ for $f\in k(X)^*$.

Let $\Div_{B^-}(X)$ denote the subgroup of $\Div(X)$ consisting of $B^-$-invariant Weil divisors. It is easy to see that $\Div_{B^-}(X)$ is freely generated by $\calD(X)$, i.e. the $B^-$-invariant prime divisors of $X$. Thus, any $\delta\in\Div_{B^-}(X)$ can be written as
\begin{align}\label{eq:Weil divisor form}
	\delta = \sum_{D\in\calD(X)} a_D D = \sum_{(\rho,\varnothing)} a_\rho D_\rho + \sum_{\alpha\in\calC} a_\alpha D_\alpha
\end{align}
for some coefficients $a_D\in\Z$ or $a_\rho,a_\alpha\in\Z$; here the first sum in the last expression is implicitly over all non-coloured rays $(\rho,\varnothing)\in\Sigma^c(1)$, these $D_\rho$ are the $G$-invariant prime divisors of $X$, and these $D_\alpha$ are the colour divisors of $X$. 

Since $N^\vee\cong k(X)^{(B^-)}/k^*$ (see \cref{lemma:image of eigenfunctions is N^vee for horospherical}), an element $m\in N^\vee$ determines a $B^-$-eigenfunction $f_m$ which is unique up to a nonzero constant multiple. The principal divisors of the form $\prindiv(f_m)$ are precisely the $B^-$-invariant principal divisors. One can check that
\begin{align}\label{eq:principal divisor expression}
	\prindiv(f_m) = \sum_{D\in\calD(X)} \langle m,u_D\rangle D = \sum_{(\rho,\varnothing)} \langle m,u_\rho\rangle D_\rho + \sum_{\alpha\in\calC} \langle m,u_\alpha\rangle D_\alpha.
\end{align}

The following Theorem tells us that $\Cl(X)$ is generated by the $B^-$-invariant prime divisors, and the relations come from the principal divisors of the form $\prindiv(f_m)$. 

\begin{theorem}[Class group]\label{thm:class group}
	Let $X=X_{\Sigma^c}$ be a $G/H$-horospherical variety. There is an exact sequence of groups
	\begin{align*}
		N^\vee &\map \Div_{B^-}(X) \map \Cl(X) \map 0
		\\ m &\longmapsto \prindiv(f_m)
	\end{align*}
	where the surjection sends $\delta\in\Div_{B^-}(X)$ to its class in $\Cl(X)$. Moreover, the following are equivalent:
	\begin{enumerate}
		\item The above exact sequence is left exact, i.e. $N^\vee\to\Div_{B^-}(X)$ is injective.
		\item $X$ has no torus factors. 
		\item $\scrO_X^*(X)=k^*$.
		\item $\{u_D:D\in\calD(X)\}$ spans $N_\R$ (e.g. $\Sigma^c$ has a coloured cone of full dimension). 
	\end{enumerate}
	
	\begin{proof}
		See \cite[Theorem 4.2.1]{perrin2018sanya} for the main claim and the equivalence between (1) and (4). The equivalence between (2), (3), and (4) is contained in \cref{prop:characterizations of torus factors}. 
	\end{proof}
\end{theorem}

\begin{example}\label{ex:class group of non-Q-factorial projective SL_3 U_3 variety}
	Consider $G=\SL_3$ and $H=U_3$. Then the coloured lattice $N$ is isomorphic to $\Z^2$, we have $\calC=\{\alpha_1,\alpha_2\}$, and $u_{\alpha_i}=e_i$ are the standard basis vectors for $i=1,2$. Let $X=X_{\Sigma^c}$ be the $\SL_3/U_3$-horospherical variety whose coloured fan is given in the diagram below. The maximal coloured cones in $\Sigma^c$ are $\sigma_1^c:=(\Cone(e_1+e_2,e_1-e_2),\{\alpha_1\})$, $\sigma_2^c:=(\Cone(-e_1,e_1+e_2),\varnothing)$, and $\sigma_3^c:=(\Cone(-e_1,e_1-e_2),\varnothing)$. We use \cref{thm:class group} to compute $\Cl(X)$. 
	
	\tikzitfig{coloured-fan-of-non-toric-non-Q-factorial-SL-3-U-3-variety}
	
	Let $D_{e_1+e_2}$, $D_{-e_1}$, and $D_{e_1-e_2}$ denote the $\SL_3$-invariant prime divisors corresponding to the non-coloured rays generated by $e_1+e_2$, $-e_1$, and $e_1-e_2$, respectively. Then $\Cl(X)$ is generated by the elements of $\calD(X)$: $D_{e_1+e_2}$, $D_{-e_1}$, $D_{e_1-e_2}$, $D_{\alpha_1}$, and $D_{\alpha_2}$. The relations are given by
	\begin{align*}
		0 &\sim \sum_{D\in\calD(X)} \langle e_1,u_D\rangle D = 1\cdot D_{e_1+e_2} + (-1)\cdot D_{-e_1} + 1\cdot D_{e_1-e_2} + 1\cdot D_{\alpha_1} + 0\cdot D_{\alpha_2},
		\\ 0 &\sim \sum_{D\in\calD(X)} \langle e_2,u_D\rangle D = 1\cdot D_{e_1+e_2} + 0\cdot D_{-e_1} + (-1)\cdot D_{e_1-e_2} + 0\cdot D_{\alpha_1} + 1\cdot D_{\alpha_2}.
	\end{align*}
	Hence, $\Cl(X)$ is isomorphic to $\Z^3$ and can be generated by $D_{e_1+e_2}$, $D_{\alpha_1}$, and $D_{\alpha_2}$. 
\end{example}

\subsection{The Picard group}\label{subsec:the picard group}

Given a $G/H$-horospherical variety $X$, we let $\CDiv(X)$ denote the subgroup of $\Div(X)$ consisting of Cartier divisors on $X$, and we let $\CDiv_{B^-}(X)$ denote the subgroup of $\CDiv(X)$ consisting of $B^-$-invariant Cartier divisors on $X$. Recall that the Picard group of $X$, denoted $\Pic(X)$, is the quotient of $\CDiv(X)$ by the subgroup generated by principal divisors. As a consequence of \cref{thm:class group}, we have an exact sequence
\begin{align}\label{eq:Picard group exact sequence}
\begin{split}
	N^\vee &\map \CDiv_{B^-}(X) \map \Pic(X) \map 0
	\\ m &\longmapsto \prindiv(f_m)
\end{split}
\end{align}
(and the same conditions from \cref{thm:class group} hold for left exactness). In particular, the Picard group $\Pic(X)$ is generated by $B^-$-invariant prime divisors. 

Similar to \cite[Section 4.2]{cox2011toric}, we can define piecewise linear functions (or ``support functions") on a coloured fan $\Sigma^c$. A \define{piecewise linear function} on $\Sigma^c$ is a function $\varphi:|\Sigma^c|\to\R$ which is linear on each cone of $\Sigma$ and satisfies $\varphi(|\Sigma^c|\cap N)\subseteq\Z$; let $\opn{PLF}(\Sigma^c)$ be the abelian group of piecewise linear functions on $\Sigma^c$. 
Since $\varphi$ is linear on each cone $\sigma\in\Sigma$, we can write $\varphi(u)=\langle m_\sigma,u\rangle$ for some $m_\sigma\in N^\vee$. So we can view $\varphi$ as a collection of piecewise linear data $\{m_\sigma\}_{\sigma\in\Sigma}$ where these $m_\sigma\in N^\vee$ satisfy $\langle m_\sigma,u\rangle=\langle m_{\sigma'},u\rangle$ for all $u\in\sigma\cap\sigma'$. Note that $\varphi$ is completely determined by the collection $\{m_\sigma\}_{\sigma\in\Sigma_{\max}}$. 

If a piecewise linear function $\varphi$ on $\Sigma^c$ is given by two sets of piecewise linear data $\{m_\sigma\}_\sigma$ and $\{m_\sigma'\}_\sigma$ then, for each $\sigma\in\Sigma$, the points $m_\sigma$ and $m_\sigma'$ have the same image in $N^\vee/(\sigma^\perp\cap N^\vee)$. In particular, if $\Sigma^c$ is complete, then the collection $\{m_\sigma\}_{\sigma\in\Sigma_{\max}}$ is unique to $\varphi$. 

\begin{proposition}[Cartier divisors]\label{prop:Cartier divisors}
	Let $X=X_{\Sigma^c}$ be a $G/H$-horospherical variety. Then 
	\begin{align*}
		\delta = \sum_{(\rho,\varnothing)} a_\rho D_\rho + \sum_{\alpha\in\calC} a_\alpha D_\alpha \in\Div_{B^-}(X)
	\end{align*}
	is Cartier if and only if there exists a piecewise linear function $\varphi$ on $\Sigma^c$ such that, for each $\sigma^c=(\sigma,\calF)\in\Sigma^c$, we have $\varphi(u_\rho)=a_\rho$ for all non-coloured rays $(\rho,\varnothing)\in\sigma^c(1)$ and $\varphi(u_\alpha)=a_\alpha$ for all $\alpha\in\calF$. 
	
	\begin{proof}
		See \cite[Lemma 4.3.1]{perrin2018sanya}. 
	\end{proof}
\end{proposition}

Given $\delta\in\CDiv_{B^-}(X)$ written in the form \cref{eq:Weil divisor form}, we denote the associated piecewise linear function from \cref{prop:Cartier divisors} by $\varphi_\delta$; note that the assignment $\delta\mapsto\varphi_\delta$ is indeed a well-defined group homomorphism $\CDiv_{B^-}\to \opn{PLF}(\Sigma^c)$. We call piecewise linear data $\{m_\sigma\}_\sigma$ for $\varphi_\delta$ \define{Cartier data} for $\delta$. Note that, if $\delta=\prindiv(f_m)$ is principal for some $m\in N^\vee$, then $\varphi_\delta$ is globally linear because $\varphi_\delta=\langle m,-\rangle$, i.e. $\varphi_\delta(u)=\langle m,u\rangle$ for all $u\in |\Sigma^c|$.

\begin{remark}\label{rmk:non-included-colour divisors are Cartier}
	Let $X=X_{\Sigma^c}$ be a $G/H$-horospherical variety and denote $\calF:=\calF(\Sigma^c)$. As a result of \cref{prop:Cartier divisors}, the Weil divisors of the form $\delta=\sum_{\alpha\in\calC\setminus\calF} a_\alpha D_\alpha$ are always Cartier and are precisely the Cartier divisors with trivial Cartier data $\varphi_\delta=0$. In particular, each colour divisor $D_\alpha\subseteq X$ with $\alpha\in\calC\setminus\calF$ is Cartier. 
\end{remark}

\begin{remark}\label{rmk:splitting of Cartier divisors group}
	Let $X=X_{\Sigma^c}$ be a $G/H$-horospherical variety and denote $\calF:=\calF(\Sigma^c)$. We have an exact sequence
	\begin{align*}
		0 \map \Z^{\calC\setminus \calF} &\map \CDiv_{B^-}(X) \map \opn{PLF}(\Sigma^c) \map 0
		\\ (a_\alpha)_{\alpha\in\calC\setminus\calF} &\longmapsto \sum_{\alpha\in\calC\setminus\calF} a_\alpha D_\alpha
	\end{align*}	
	(injectivity from \cref{rmk:non-included-colour divisors are Cartier}, and surjectivity uses the argument in \cite[Theorem 4.2.12]{cox2011toric}). Since $\opn{PLF}(\Sigma^c)$ is a free abelian group, the exact sequence splits. Combining this with \cref{eq:Picard group exact sequence}, we have a surjection
	\begin{align}\label{eq:surjection of PLF onto Picard}
		\begin{split}
			\Z^{\calC\setminus\calF} \oplus \opn{PLF}(\Sigma^c) &\surjmap \Pic(X)
			\\ \left( (a_\alpha)_{\alpha\in\calC\setminus\calF} , \varphi \right) &\mapsto \sum_{\alpha\in\calC\setminus\calF} a_\alpha [D_\alpha] + \sum_{(\rho,\varnothing)} \varphi(u_\rho) [D_\rho] + \sum_{\alpha\in\calF} \varphi(u_\alpha) [D_\alpha].
		\end{split}
	\end{align}
\end{remark}

The following theorem tells us how to compute $\Pic(X_{\Sigma^c})$ in a way that is similar to \cite[Theorem 4.2.12]{cox2011toric} for toric varieties.

\begin{theorem}[Picard group]\label{thm:Picard group}
	Let $X=X_{\Sigma^c}$ be a $G/H$-horospherical variety and denote $\calF:=\calF(\Sigma^c)$. There is an exact sequence 
	\begin{alignat*}{2}
		N^\vee &\map \Z^{\calC\setminus\calF} \oplus \opn{PLF}(\Sigma^c) \map \Pic(X) \map 0
		\\ m &\longmapsto \left((\langle m,u_\alpha\rangle)_{\alpha\in\calC\setminus\calF},\langle m,-\rangle\right)
%
	\end{alignat*}
	where the right map is \cref{eq:surjection of PLF onto Picard}. 
	Moreover, the sequence is left exact, i.e. the left-most map is injective, if and only if any of the equivalent conditions (1)-(4) in \cref{thm:class group} are satisfied.

	\begin{proof}
		This follows from \cref{eq:surjection of PLF onto Picard} and \cref{eq:Picard group exact sequence}.
	\end{proof}
\end{theorem}

\begin{remark}\label{rmk:other exact Picard sequence}
	Let $X=X_{\Sigma^c}$ be a $G/H$-horospherical variety and denote $\calF:=\calF(\Sigma^c)$. We assume that $|\Sigma^c|$ spans $N_\R$, i.e. the equivalent conditions (1)-(4) of \cref{thm:class group} are satisfied for left exactness. The following exact sequence is presented in \cite[Theorem 4.3.4]{perrin2018sanya} instead of the one we give above in \cref{thm:Picard group}:
	\begin{align*}
		0 \map \Z^{\calC\setminus\calF} \map \Pic(X) \map \opn{PLF}(\Sigma^c)/N^\vee \map 0
	\end{align*}
	where the first map is the first map $\Z^{\calC\setminus\calF}\to\CDiv_{B^-}(X)$ from \cref{rmk:splitting of Cartier divisors group} followed by the projection to $\Pic(X)$, and the second map sends $[\delta]\in\Pic(X)$ to its associated $\varphi_\delta$, but only defined up to globally linear functions $\langle m,-\rangle$ with $m\in N^\vee$. 
	
	It is important to note that $\opn{PLF}(\Sigma^c)/N^\vee$ is not always torsion free, e.g. \cite[Example 4.2.3]{cox2011toric} is an example of a toric variety with fan $\Sigma_0$ depending on an arbitrary $d\in\Z_{>1}$, and in this case $\opn{PLF}(\Sigma_0)/N^\vee\cong \Z/d\Z$. However, it is torsion free when $\Sigma$ contains a cone of full dimension (see below). 
\end{remark}

\begin{corollary}\label{cor:consequences of Picard group theorem}
	Let $X=X_{\Sigma^c}$ be a $G/H$-horospherical variety and denote $\calF:=\calF(\Sigma^c)$.
	\begin{enumerate}
		\item If $\Sigma$ contains a cone of full dimension, then $\Pic(X)$ is a free abelian group. In this case, $$\Pic(X)\cong \Z^{\calC\setminus\calF}\oplus(\opn{PLF}(\Sigma^c)/N^\vee).$$ 
		\item If $X$ is simple, i.e. $\Sigma^c=\sigma^c$ is generated by a single coloured cone, and if $\sigma$ has full dimension, then $\Pic(X)=\Z^{\calC\setminus\calF}$.
		\item If $X$ is affine, then $\Pic(X)=0$. 
	\end{enumerate}

	\begin{proof}
		If $X$ is affine, then $\calC\setminus\calF=\varnothing$ and $\Sigma^c=\sigma^c$ is generated by a single coloured cone, so \cref{thm:Picard group} gives an exact sequence $N^\vee\to\opn{PLF}(\sigma^c)\to\Pic(X)\to 0$. Since $\sigma^c$ is a single coloured cone, the map $N^\vee\to\opn{PLF}(\sigma^c)$ is surjective, hence $\Pic(X)=0$, establishing (3).
		
		For the remainder, we may assume that $\Sigma$ contains a cone of full dimension. In this case, $\opn{PLF}(\Sigma^c)/N^\vee$ is a free abelian group, so $\Pic(X)$ is free by \cref{rmk:other exact Picard sequence} (alternatively, $N^\vee$ is saturated in $\Z^{\calC\setminus\calF}\oplus\opn{PLF}(\Sigma^c)$, so we can use \cref{thm:Picard group}; compare with the proof of \cite[Proposition 4.2.5]{cox2011toric}). This establishes (1) by, for example, splitting the exact sequence in \cref{rmk:other exact Picard sequence}. 
		
		Now we further assume that $\Sigma^c=\sigma^c$ is generated by a single coloured cone. Then $\opn{PLF}(\sigma^c)/N^\vee=0$, so (2) follows from (1). 
	\end{proof}
\end{corollary}

\begin{example}\label{ex:Picard group of non-Q-factorial SL_3/U_3 variety}
	Consider $G=\SL_3$ and $H=U_3$. Then the coloured lattice $N$ is isomorphic to $\Z^2$, we have $\calC=\{\alpha_1,\alpha_2\}$, and $u_{\alpha_i}=e_i$ are the standard basis vectors for $i=1,2$. Let $X=X_{\Sigma^c}$ be the $\SL_3/U_3$-horospherical variety from \cref{ex:class group of non-Q-factorial projective SL_3 U_3 variety}; we recall the diagram for $\Sigma^c$ below. We use \cref{thm:Picard group} (and \cref{cor:consequences of Picard group theorem}) to compute $\Pic(X)$. 
	
	\tikzitfig{coloured-fan-of-non-toric-non-Q-factorial-SL-3-U-3-variety}
	
	We first compute $\opn{PLF}(\Sigma^c)/N^\vee$. Since we are looking at the piecewise linear function $\{m_\sigma\}_\sigma$ up to adding a globally linear function, we may assume that $m_{\sigma_1}=0$. Then we require
	\begin{alignat*}{3}
		\begin{cases}
			\langle m_{\sigma_2},e_1+e_2\rangle=0
			\\ \langle m_{\sigma_3},e_1-e_2\rangle=0
			\\ \langle m_{\sigma_2}-m_{\sigma_3},-e_1\rangle=0
		\end{cases} \quad\implies 
		\begin{split}
			& m_{\sigma_2}=(a,-a)
			\\& m_{\sigma_3}=(a,a)
		\end{split} \quad\text{for some } a\in\Z.
	\end{alignat*}
	In particular, $\opn{PLF}(\Sigma^c)/N^\vee\cong\Z$. Since $\Z(\calC\setminus\calF(\Sigma^c))=\Z D_{\alpha_2}\cong\Z$, we find that $\Pic(X)$ is isomorphic to $\Z^2$ and can be generated by $D_{-e_1}$ and $D_{\alpha_2}$. Note that when $a=-1$ above, $\{m_\sigma\}_\sigma$ is Cartier data for $D_{-e_1}$. 
\end{example}

\subsection{Basepoint free and ample}\label{subsec:basepoint free and ample}

Similar to \cite[Section 6.1]{cox2011toric}, we can consider convex and strictly convex piecewise linear functions, and see how they correspond to basepoint free (i.e. globally generated) and ample Cartier divisors. 

Let $\varphi:|\Sigma^c|\to\R$ be a piecewise linear function on $\Sigma^c$, and let $\{m_\sigma\}_\sigma$ be piecewise linear data for $\varphi$. We say that $\varphi$ is \define{convex} (resp. \define{strictly convex}) if, for each $\sigma\in\Sigma$, we have $\langle m_{\sigma'},u\rangle\leq \langle m_\sigma,u\rangle$ (resp. $\langle m_{\sigma'},u\rangle< \langle m_\sigma,u\rangle$) for all $\sigma'\in\Sigma\setminus\{\sigma\}$ and all $u\in\sigma\setminus\sigma'$. 

\begin{proposition}[Basepoint free and ample criteria]\label{prop:basepoint free and ample}
	Let $X=X_{\Sigma^c}$ be a complete $G/H$-horospherical variety, and let 
	\begin{align*}
		\delta = \sum_{(\rho,\varnothing)} a_\rho D_\rho + \sum_{\alpha\in\calC} a_\alpha D_\alpha \in\CDiv_{B^-}(X).
	\end{align*}
	Then $\delta$ is basepoint free (resp. ample) if and only if $\varphi_\delta$ is convex (resp. strictly convex) and, for each $\alpha\in\calC\setminus\calF(\Sigma^c)$, we have $\varphi_\delta(u_\alpha)\leq a_\alpha$ (resp. $\varphi_\delta(u_\alpha)<a_\alpha$). 
	
	\begin{proof}
		See \cite[Corollary 4.3.8]{perrin2018sanya}. 
	\end{proof}
\end{proposition}

\begin{remark}\label{rmk:non-included-colour divisors are basepoint free}
	Let $X=X_{\Sigma^c}$ be a $G/H$-horospherical variety. Any Weil divisor of the form $\delta=\sum_{\alpha\in\calC\setminus\calF(\Sigma^c)} a_\alpha D_\alpha$ is basepoint free by \cref{prop:basepoint free and ample}. In particular, each colour divisor $D_\alpha\subseteq X$ with $\alpha\in\calC\setminus\calF(\Sigma^c)$ is basepoint free. 
\end{remark}

\begin{corollary}[Characterization of projective]\label{cor:characterization of projective}
	Let $X_{\Sigma^c}$ be a complete $G/H$-horospherical variety. Then the following are equivalent:
	\begin{enumerate}
		\item $X_{\Sigma^c}$ is projective.
		\item There exists a strictly convex piecewise linear function on $\Sigma^c$. 
		\item The toric variety $X_{\Sigma}$ is projective.
	\end{enumerate}
	In particular, if $\dim_\R(N_\R)\leq 2$, then $X_{\Sigma^c}$ is projective. 

	\begin{proof}
		(1) and (2) are equivalent by \cref{prop:basepoint free and ample}. Since (2) only depends on the underlying fan $\Sigma$, it follows that (2) and (3) are equivalent. 
	\end{proof}
\end{corollary}

\begin{example}\label{ex:basepoint free and ample divisors for non-Q-factorial SL_3 U_3 variety}
	Consider $G=\SL_3$ and $H=U_3$. Then the coloured lattice $N$ is isomorphic to $\Z^2$, we have $\calC=\{\alpha_1,\alpha_2\}$, and $u_{\alpha_i}=e_i$ are the standard basis vectors for $i=1,2$. Let $X=X_{\Sigma^c}$ be the $\SL_3/U_3$-horospherical variety from \cref{ex:class group of non-Q-factorial projective SL_3 U_3 variety}; we recall the diagram for $\Sigma^c$ below. We use \cref{prop:basepoint free and ample} to check which Cartier divisors on $X$ are basepoint free or ample. 
	
	\tikzitfig{coloured-fan-of-non-toric-non-Q-factorial-SL-3-U-3-variety} 
	
	Recall from \cref{ex:Picard group of non-Q-factorial SL_3/U_3 variety} that $\Pic(X)$ is generated by $D_{-e_1}$ and $D_{\alpha_2}$. For a general Cartier divisor of the form $\delta=aD_{-e_1}+bD_{\alpha_2}$, with $a,b\in\Z$, we showed that the Cartier data $\{m_\sigma\}_\sigma$ is given by $m_{\sigma_1}=0$, $m_{\sigma_2}=(-a,a)$, and $m_{\sigma_3}=(-a,-a)$. One can check that $\varphi_\delta$ is convex (resp. strictly convex) if and only if $a\geq 0$ (resp. $a>0$). Therefore, $\delta$ is basepoint free (resp. ample) if and only if $a\geq 0$ and $a\leq b$ (resp. $a>0$ and $a<b$). 
\end{example}

\subsection{The canonical divisor}\label{subsec:the canonical divisor}

To end this section we examine the canonical divisor $K_X$ of a horospherical $G$-variety $X$. The coefficients are given in terms of the pairing between roots and dual roots via the Dynkin diagram of $G$. For convenience, we state the following proposition for the \textit{anti}canonical divisor $-K_X$. 

\begin{proposition}[Anticanonical divisor]\label{prop:anticanonical divisor}
	Let $X$ be a $G/H$-horospherical variety. The anticanonical divisor of $X$ is given by 
	\begin{align*}
		-K_X = \sum_{(\rho,\varnothing)} D_\rho + \sum_{\alpha\in\calC} b_\alpha D_\alpha
	\end{align*}
	where the coefficient $b_\alpha$ in front of the colour divisor $D_\alpha\subseteq X$ is given by 
	\begin{align*}
		b_\alpha = \sum_{\gamma\in R^+\setminus R_I} \left\langle \gamma, \alpha^\vee \right\rangle_R
	\end{align*}
	where $\langle \gamma,\alpha^\vee\rangle_R$ is the pairing between roots and dual roots in $R$.\footnote{See \cite[Page 6]{pasquier2009introduction} for a quick reference on this root pairing.} Moreover, $b_\alpha\in\Z_{\geq 2}$ for each $\alpha\in\calC$. 

	\begin{proof}
		This follows from Theorem \cite[Theorem 6.1.2]{perrin2018sanya} and \cite[Proposition 6.1.1]{perrin2018sanya} (see the following remark). 
	\end{proof}
\end{proposition}

\begin{remark}\label{rmk:anticanonical divisor parts}
	 One should think of the anticanonical divisor $-K_X$ in \cref{prop:anticanonical divisor} as follows. First, $-K_X=-K_{X'}$ where $X'\subseteq X$ is the open toroidal subvariety from \cref{ex:open toroidal subvariety}. By \cref{subsec:affine local structure}, the fibre of the map $X'\to G/P$ is the toric variety whose fan is the underlying fan associated to $X'$. Thus, $\sum_{(\rho,\varnothing)} D_\rho$ is the anticanonical divisor of the fibre, and $\sum_{\alpha\in\calC} b_\alpha D_\alpha$ is the anticanonical divisor of the base $G/P$. 
\end{remark}

\begin{example}\label{ex:anticanonical divisor of SL_n/B_n}
	Consider $G=\SL_n$. If $H$ is any (standard) horospherical subgroup of $\SL_n$ with associated parabolic $B_n$ (e.g. $H=U_n$), then $I=\varnothing$ so $\calC=S\setminus I=S=\{\alpha_1,\ldots,\alpha_{n-1}\}$. Using the notation of \cref{prop:anticanonical divisor}, one can check that $b_\alpha=2$ for each $\alpha\in\calC$; for example, in $\SL_3$ we have
	\begin{align*}
		b_{\alpha_1} = \langle \alpha_1,\alpha_1^\vee\rangle_R + \langle \alpha_2,\alpha_1^\vee\rangle_R + \langle \alpha_1+\alpha_2,\alpha_1^\vee\rangle_R = 2+(-1)+(2+(-1)) = 2.
	\end{align*}
	Therefore, if $X=X_{\Sigma^c}$ is any $\SL_n/H$-horospherical variety, then
	\begin{align*}
		-K_X = \sum_{(\rho,\varnothing)} D_\rho + \sum_{\alpha\in\calC} 2D_\alpha.
	\end{align*}
	
	For example, if $X$ is the $\SL_3/U_3$-horospherical variety from \cref{ex:class group of non-Q-factorial projective SL_3 U_3 variety}, then
	\begin{align*}
		-K_X =& D_{e_1+e_2}+D_{-e_1}+D_{e_1-e_2}+2D_{\alpha_1}+2D_{\alpha_2}.
		\qedhere
	\end{align*}
\end{example}

\section*{Acknowledgements}

I want to thank my PhD advisor Matthew Satriano for originally suggesting that I study horospherical geometry, and for being very helpful and encouraging throughout my time studying these objects. Also, thank you to Xuemiao Chen and Changho Han for helpful conversations. 

I am grateful for the support from NSERC via a PGS-D scholarship (reference number: PGSD3-558713-2021). 

The arXiv V2 of this article was submitted to a journal (and was rejected), so I would like to thank the anonymous referee for providing very detailed and helpful feedback which contributed to several changes that appear in this version of the article.

\printbibliography

\end{document}

%% file: Figures/sample.tikzstyles
\tikzstyle{orange dot}=[fill={rgb,255: red,255; green,128; blue,0}, draw=black, shape=circle, scale=0.75]
\tikzstyle{empty dot}=[fill=none, draw=black, shape=circle, scale=0.75]
\tikzstyle{vertical marker}=[fill=none, draw=none, shape=circle, label={center:{\footnotesize $|$}}]
\tikzstyle{origin}=[fill=none, draw=none, shape=circle, label={center:{\Large $|$}}]
\tikzstyle{dot}=[fill=black, draw=black, shape=circle, scale=0.5]
\tikzstyle{fantastack square}=[fill=none, draw=black, shape=rectangle, scale=1.25]
\tikzstyle{big orange dot}=[fill={rgb,255: red,255; green,128; blue,0}, draw=black, shape=circle, scale=1.00]
\tikzstyle{big empty dot}=[fill=none, draw=black, shape=circle, scale=1.00]

\tikzstyle{arrow}=[->, line width=1.0]
\tikzstyle{dashed line}=[-, dashed, fill=none]
\tikzstyle{cone fill}=[-, draw=none, fill={rgb,255: red,191; green,191; blue,191}]
\tikzstyle{backed arrow}=[{|->}, line width=1]
\tikzstyle{filled line}=[-, line width=1]